\documentclass{amsart}

\usepackage{amsmath}
\usepackage[latin1]{inputenc}
\usepackage{amsfonts}
\usepackage{amssymb}
\usepackage{graphicx,epsfig}
\usepackage{amsthm}
\usepackage{amscd}
\usepackage[all]{xy}
\usepackage[latin1]{inputenc}

\usepackage[T1]{fontenc}

\usepackage{amsmath}

\usepackage{amsfonts}

\usepackage{amssymb}

\usepackage{amsthm}

\usepackage{graphics}

\newcommand{\mk}{\medskip}

\newcommand{\ZZ}{\mathbb{Z}}

\newcommand{\CC}{\mathbb{C}}

\newcommand{\NN}{\mathbb{N}}

\newcommand{\QQ}{\mathbb{Q}}

\newcommand{\Glie}{\mathfrak{g}}

\newcommand{\Yim}{\mathcal{Y}}

\newcommand{\Hlie}{\mathfrak{h}}

\newcommand{\demo}{\noindent {\it \small Proof:}\quad}

\renewcommand{\NN}{\ensuremath{\mathbb{N}}}

\renewcommand{\CC}{\ensuremath{\mathbb{C}}}

\renewcommand{\QQ}{\ensuremath{\mathbb{Q}}}

\newcommand{\U}{\mathcal{U}}

\newtheorem{thm}{Theorem}[section]

\newtheorem{defi}[thm]{Definition}

\newtheorem{cor}[thm]{Corollary}

\newtheorem{prop}[thm]{Proposition}

\newtheorem{lem}[thm]{Lemma}

\newtheorem{conj}[thm]{Conjecture}

\author{David Hernandez}\thanks{\'Ecole Normale Sup\'erieure - DMA, 45, Rue d'Ulm F-75230 PARIS, Cedex 05  FRANCE
 \newline email: David.Hernandez@ens.fr, URL: http://www.dma.ens.fr/$\sim$dhernand}
\address{David Hernandez: \'Ecole Normale Sup\'erieure - DMA, 45, Rue d'Ulm F-75230 PARIS,
  Cedex 05  FRANCE}
\email{David.Hernandez@ens.fr\\ URL: http://www.dma.ens.fr/$\sim$dhernand}

\usepackage[all]{xy}
\title{Algebraic Approach to $q,t$-Characters}

\begin{document}

\begin{abstract} Frenkel and Reshetikhin \cite{Fre} introduced $q$-characters to study finite dimensional representations of the quantum affine algebra $\U_q(\hat{\Glie})$. In the simply laced case Nakajima \cite{Naa}\cite{Nab} defined deformations of $q$-characters called $q,t$-characters. The definition is combinatorial but the proof of the existence uses the geometric theory of quiver varieties which holds only in the simply laced case. In this article we propose an algebraic general (non necessarily simply laced) new approach to $q,t$-characters motivated by the deformed screening operators \cite{Her01}. The $t$-deformations are naturally deduced from the structure of $\U_q(\hat{\Glie})$: the parameter $t$ is analog to the central charge $c\in\U_q(\hat{\Glie})$. The $q,t$-characters lead to the construction of a quantization of the Grothendieck ring and to general analogues of Kazhdan-Lusztig polynomials in the same spirit as Nakajima did for the simply laced case.\end{abstract}

\maketitle

\tableofcontents

\section{Introduction}
 
 We suppose $q\in\CC^*$ is not a root of unity. In the case of a semi-simple Lie algebra $\Glie$, the structure of the Grothendieck ring $\text{Rep}(\U_q(\Glie))$ of finite dimensional representations of the quantum algebra $\U_q(\Glie)$ is well understood. It is analogous to the classical case $q=1$. In particular we have ring isomorphisms:
$$\text{Rep}(\U_q(\Glie))\simeq \text{Rep}(\Glie)\simeq \ZZ [\Lambda]^W\simeq \ZZ[T_1,...,T_n]$$ 
deduced from the injective homomorphism of characters $\chi$:
$$\chi(V)=\underset{\lambda\in\Lambda}{\sum}\text{dim}(V_{\lambda})\lambda$$
where $V_{\lambda}$ are weight spaces of a representation $V$ and $\Lambda$ is the weight lattice.

\noindent For the general case of Kac-Moody algebras the picture is less clear. In the affine case $\U_q(\hat{\Glie})$, Frenkel and Reshetikhin \cite{Fre} introduced an injective ring homomorphism of $q$-characters:
$$\chi_q:\text{Rep}(\U_q(\hat{\Glie}))\rightarrow \ZZ[Y_{i,a}^{\pm}]_{1\leq i\leq n,a\in\CC^*}=\Yim$$

\noindent The homomorphism $\chi_q$ allows to describe the ring $\text{Rep}(\U_q(\hat{\Glie}))\simeq\ZZ[X_{i,a}]_{i\in I,a\in\CC^*}$, where the $X_{i,a}$ are fundamental representations. It particular $\text{Rep}(\U_q(\hat{\Glie}))$ is commutative. 

\noindent The morphism of $q$-characters has a symmetry property analogous to the classical action of the Weyl group $\text{Im}(\chi)=\ZZ[\Lambda]^W$: Frenkel and Reshetikhin defined $n$ screening operators $S_i$ such that $\text{Im}(\chi_q)=\underset{i\in I}{\bigcap}\text{Ker}(S_i)$ (the result was proved by Frenkel and Mukhin for the general case in \cite{Fre2}).

\noindent In the simply laced case Nakajima introduced $t$-analogues of $q$-characters (\cite{Naa}, \cite{Nab}): it is a $\ZZ[t^{\pm}]$-linear map 
$$\chi_{q,t}:\text{Rep}(\U_q(\hat{\Glie}))\otimes_{\ZZ}\ZZ[t^{\pm}]\rightarrow\Yim_t=\ZZ[Y_{i,a}^{\pm},t^{\pm}]_{i\in I,a\in\CC^*}$$ 
which is a deformation of $\chi_q$ and multiplicative in a certain sense. A combinatorial axiomatic definition of $q,t$-characters is given. But the existence is non-trivial and is proved with the geometric theory of quiver varieties which holds only in the simply laced case. 

\noindent In \cite{Her01} we introduced $t$-analogues of screening operators $S_{i,t}$ such that in the simply laced case:
$$\underset{i\in I}{\bigcap}\text{Ker}(S_{i,t})=\text{Im}(\chi_{q,t})$$
\noindent It is a first step in the algebraic approach to $q,t$-characters proposed in this article: we define and construct $q,t$-characters in the general (non necessarily simply laced) case. The motivation of the construction appears in the non-commutative structure of the Cartan subalgebra $\U_q(\hat{\Hlie})\subset\U_q(\hat{\Glie})$, the study of screening currents and of deformed screening operators. 

\noindent As an application we construct a deformed algebra structure and an involution of the Grothendieck ring, and analogues of Kazhdan-Lusztig polynomials in the general case in the same spirit as Nakajima did for the simply laced case. In particular this article proves a conjecture that Nakajima made for the simply laced case (remark 3.10 in \cite{Nab}): there exists a purely combinatorial proof of the existence of $q,t$-characters.

This article is organized as follows: after some backgrounds in section \ref{back}, we define a deformed non-commutative algebra structure on $\Yim_t=\ZZ[Y_{i,a}^{\pm},t^{\pm}]_{i\in I,a\in\CC^*}$ (section \ref{defoal}): it is naturally deduced from the relations of $\U_q(\hat{\Hlie})\subset\U_q(\hat{\Glie})$ (theorem \ref{dessus}) by using the quantization in the direction of the central element $c$. In particular in the simply laced case it can be used to construct the deformed multiplication of Nakajima \cite{Nab} (proposition \ref{form}) and of Varagnolo-Vasserot \cite{Vas} (section \ref{varva}). 

\noindent This picture allows us to introduce the deformed screening operators of \cite{Her01} as commutators of Frenkel-Reshetikhin's screening currents of \cite{Freb} (section \ref{scr}). In \cite{Her01} we gave explicitly the kernel of each deformed screening operator (theorem \ref{her}). 

\noindent In analogy to the classic case where $\text{Im}(\chi_q)=\underset{i\in I}{\bigcap}\text{Ker}(S_i)$, we have to describe the intersection of the kernels of deformed screening operators. We introduce a completion of this intersection (section \ref{complesection}) and give its structure in proposition \ref{thth}. It is easy to see that it is not too big (lemma \ref{leasto}); but the point is to prove that it contains enough elements: it is the main result of our construction in theorem \ref{con} which is crucial for us. It is proved by induction on the rank $n$ of $\Glie$. 

\noindent We define a $t$-deformed algorithm (section \ref{defialgo}) analog to the Frenkel-Mukhin's algorithm \cite{Fre2} to construct $q,t$-characters in the completion of $\Yim_t$. An algorithm was also used by Nakajima in the simply laced case in order to compute the $q,t$-characters for some examples (\cite{Naa}) assuming they exist (which was geometrically proved). Our aim is different : we do not know {\it a priori} the existence in the general case. That is why we have to show the algorithm is well defined, never fails (lemma \ref{nfail}) and gives a convenient element (lemma \ref{conv}). 

\noindent This construction gives $q,t$-characters for fundamental representations; we deduce from them the injective morphism of $q,t$-characters $\chi_{q,t}$ (definition \ref{mqt}). We study the properties of $\chi_{q,t}$ (theorem \ref{axiomes}). Some of them are generalization of the axioms that Nakajima defined in the simply laced case (\cite{Nab}); in particular we have constructed the morphism of \cite{Nab}. 

\noindent We have some applications: the morphism gives a deformation of the Grothendieck ring  because the image of $\chi_{q,t}$ is a subalgebra for the deformed multiplication (section \ref{quanta}). Moreover we define an antimultiplicative involution of the deformed Grothendieck ring (section \ref{invo}); the construction of this involution is motivated by the new point view adopted in this paper : it is just replacing $c$ by $-c$ in $\U_q(\hat{\Glie})$. In particular we define constructively analogues of Kazhdan-Lusztig polynomials and a canonical basis (theorem \ref{expol}) motivated by the introduction of \cite{Nab}. We compute explicitly the polynomials for some examples.

\noindent In section \ref{quest} we raise some questions : we conjecture that the coefficients of $q,t$-characters are in $\NN[t^{\pm}]\subset\ZZ[t^{\pm}]$. In the $ADE$-case it a result of Nakajima;  we give an alternative elementary proof for the $A$-cases in section \ref{acase}. The cases $G_2, B_2, C_2$ are also checked in section \ref{fin}. The cases $F_4, B_n, C_n$ ($n\leq 10$) have been checked on a computer.

\noindent We also conjecture that the generalized analogues to Kazhdan-Lusztig polynomials give at $t=1$ the multiplicity of simple modules in standard modules. We propose some generalizations and further applications which will be studied elsewhere.

\noindent In the appendix (section \ref{fin}) we give explicit computations of $q,t$-characters for semi-simple Lie algebras of rank 2. They are used in the proof of theorem \ref{con}.

\noindent For convenience of the reader we give at the end of this article an index of notations defined in the main body of the text.

\noindent {\bf Acknowledgments.} The author would like to thank M. Rosso for encouragements and precious comments on a previous version of this paper, I. B. Frenkel for having encouraged him in this direction, E. Frenkel for encouragements, useful discussions and references,  E. Vasserot for very interesting explanations about \cite{Vas}, O. Schiffmann for valuable comments and his kind hospitality in Yale university, and T. Schedler for help on programming.

\section{Background}\label{back}

\subsection{Cartan matrix}\label{recalu}

A generalized Cartan matrix of rank $n$ is a matrix $C=(C_{i,j})_{1\leq i,j\leq n}$\label{carmat} such that $C_{i,j}\in\ZZ$ and:
$$C_{i,i}=2$$
$$i\neq j\Rightarrow C_{i,j}\leq 0$$
$$C_{i,j}=0\Leftrightarrow C_{j,i}=0$$
Let $I = \{1,...,n\}$.

\noindent We say that $C$ is symmetrizable if there is a matrix $D=\text{diag}(r_1,...,r_n)$ ($r_i\in\NN^*$) such that $B=DC$\label{symcar} is symmetric. 

\noindent Let $q\in\CC^*$ be the parameter of quantization. In the following we suppose it is not a root of unity. $z$ is an indeterminate.\label{qz}

\noindent If $C$ is symmetrizable, let $q_i=q^{r_i}$, $z_i=z^{r_i}$ and $C(z)=(C(z)_{i,j})_{1\leq i,j\leq n}$ the matrix with coefficients in $\ZZ[z^{\pm}]$ such that:
$$C(z)_{i,j}=[C_{i,j}]_z\text{ if $i\neq j$}$$
$$C(z)_{i,i}=[C_{i,i}]_{z_i}=z_i+z_i^{-1}$$
where for $l\in\ZZ$ we use the notation: 
$$[l]_z=\frac{z^l-z^{-l}}{z-z^{-1}}\text{ ($=z^{-l+1}+z^{-l+3}+...+z^{l-1}$ for $l\geq 1$)}$$
In particular, the coefficients of $C(z)$ are symmetric Laurent polynomials (invariant under $z\mapsto z^{-1}$). We define the diagonal matrix $D_{i,j}(z)=\delta_{i,j}[r_i]_z$ and the matrix $B(z)=D(z)C(z)$.

\noindent In the following we suppose that $C$ is of finite type, in particular $\text{det}(C)\neq 0$. In this case $C$ is symmetrizable; if $C$ is indecomposable there is a unique choice of $r_i\in\NN^*$ such that $r_1\wedge...\wedge r_n=1$. We have $B_{i,j}(z)=[B_{i,j}]_z$ and $B(z)$ is symmetric. See \cite{bou} or \cite{Kac} for a classification of those finite Cartan matrices.

\noindent We say that $C$ is simply-laced if $r_1=...=r_n=1$. In this case $C$ is symmetric, $C(z)=B(z)$ is symmetric. In the classification those matrices are of type $ADE$.

\noindent  Denote by $\mathfrak{U}\subset\QQ(z)$\label{mathu} the subgroup $\ZZ$-linearly spanned by the $\frac{P(z)}{Q(z^{-1})}$ such that $P(z)\in\ZZ[z^{\pm}]$, $Q(z)\in\ZZ[z]$, the zeros of $Q(z)$ are roots of unity and $Q(0)=1$. It is a subring of $\QQ(z)$, and for $R(z)\in\mathfrak{U},m\in\ZZ$ we have $R(q^m)\in\mathfrak{U}$ and $R(q^m)\in\CC$ makes sense. 

\noindent It follows from lemma 1.1 of \cite{Fre2} that $C(z)$ has inverse $\tilde{C}(z)$\label{invcar} with coefficients of the form $R(z)\in\mathfrak{U}$. 

\subsection{Finite quantum algebras} We refer to \cite{Ro} for the definition of the finite quantum algebra $\U_q(\Glie)$ associated to a finite Cartan matrix, the definition and properties of the type $1$-representations of $\U_q(\Glie)$, the Grothendieck ring $\text{Rep}(\U_q(\Glie))$ and the injective ring morphism of characters $\chi:\text{Rep}(\U_q(\Glie))\rightarrow \ZZ[y_i^{\pm}]$.

\subsection{Quantum affine algebras} The quantum affine algebra associated to a finite Cartan matrix $C$ is the $\CC$-algebra $\U_q(\hat{\Glie})$\label{qadefi} defined (Drinfeld new realization) by generators $x_{i,m}^{\pm}$ ($i\in I$, $m\in\ZZ$), $k_i^{\pm}$ ($i\in I$), $h_{i,m}$ ($i\in I$, $m\in\ZZ^*$), central elements $c^{\pm\frac{1}{2}}$, and relations:
$$k_ik_j=k_jk_i$$
$$k_ih_{j,m}=h_{j,m}k_i$$
$$k_ix_{j,m}^{\pm}k_i^{-1}=q^{\pm B_{ij}}x_{j,m}^{\pm}$$
$$[h_{i,m},x_{j,m'}^{\pm}]=\pm \frac{1}{m}[mB_{ij}]_qc^{\mp\frac{\mid m\mid}{2}} x_{j,m+m'}^{\pm}$$
$$x_{i,m+1}^{\pm}x_{j,m'}^{\pm}-q^{\pm B_{ij}}x_{j,m'}^{\pm}x_{i,m+1}^{\pm}=q^{\pm B_{ij}}x_{i,m}^{\pm}x_{j,m'+1}^{\pm}-x_{j,m'+1}^{\pm}x_{i,m}^{\pm}$$
$$[h_{i,m},h_{j,m'}]=\delta_{m,-m'}\frac{1}{m}[mB_{ij}]_q\frac{c^m-c^{-m}}{q-q^{-1}}$$
$$[x_{i,m}^+,x_{j,m'}^-]= \delta_{ij}\frac{c^{\frac{m-m'}{2}}\phi^+_{i,m+m'}-c^{-\frac{m-m'}{2}}\phi^-_{i,m+m'}}{q_i-q_i^{-1}}$$
$$\underset{\pi\in \Sigma_s}{\sum}\underset{k=0..s}{\sum}(-1)^k\begin{bmatrix}s\\k\end{bmatrix}_{q_i}x_{i,m_{\pi(1)}}^{\pm}...x_{i,m_{\pi(k)}}^{\pm}x_{j,m'}^{\pm}x_{i,m_{\pi(k+1)}}^{\pm}...x_{i,m_{\pi(s)}}^{\pm}=0$$
where the last relation holds for all $i\neq j$, $s=1-C_{ij}$, all sequences of integers $m_1,...,m_s$. $\Sigma_s$ is the symmetric group on $s$ letters. For $i\in I$ and $m\in\ZZ$, $\phi_{i,m}^{\pm}\in \U_q(\hat{\Glie})$ is determined by the formal power series in $\U_q(\hat{\Glie})[[u]]$ (resp. in $\U_q(\hat{\Glie})[[u^{-1}]]$):
$$\underset{m=0..\infty}{\sum}\phi_{i,\pm m}^{\pm}u^{\pm m}=k_i^{\pm}\text{exp}(\pm(q-q^{-1})\underset{m'=1..\infty}{\sum}h_{i,\pm m'}u^{\pm m'})$$
and $\phi_{i,m}^+=0$ for $m<0$, $\phi_{i,m}^-=0$ for $m>0$.

\noindent One has an embedding $\U_q(\Glie)\subset\U_q(\hat{\Glie})$ and a Hopf algebra structure on $\U_q(\hat{\Glie})$ (see \cite{Fre} for example).

\noindent The Cartan algebra $\U_q(\hat{\Hlie})\subset\U_q(\hat{\Glie})$\label{qhdefi} is the $\CC$-subalgebra of $\U_q(\hat{\Glie})$ generated by the $h_{i,m},c^{\pm}$ ($i\in I, m\in\ZZ-\{0\}$).

\subsection{Finite dimensional representations of $\U_q(\hat{\Glie})$}

A finite dimensional representation $V$ of $\U_q(\hat{\Glie})$ is called of type $1$ if $c$ acts as $\text{Id}$ and $V$ is of type $1$ as a representation of $\U_q(\Glie)$. Denote by $\text{Rep}(\U_q(\hat{\Glie}))$ the Grothendieck ring of finite dimensional representations of type $1$.

\noindent The operators $\{\phi^{\pm}_{i,\pm m},i\in I, m\in\ZZ\}$ commute on $V$. So we have a pseudo-weight space decomposition: 
$$V=\underset{\gamma\in \CC^{I\times \ZZ}\times\CC^{I\times\ZZ}}{\bigoplus} V_{\gamma}$$ 
where for $\gamma=(\gamma^+,\gamma^-)$, $V_{\gamma}$ is a simultaneous generalized eigenspace:
$$V_{\gamma}=\{x\in V/\exists p\in\NN,\forall i\in\{1,...,n\},\forall m\in\ZZ,(\phi_{i,m}^{\pm}-\gamma_{i,m}^{\pm})^p.x=0\}$$
The $\gamma_{i,m}^{\pm}$ are called pseudo-eigen values of $V$.
\begin{thm}({\bf Chari, Pressley} \cite{Cha},\cite{Cha2}) Every simple representation $V\in\text{Rep}(\U_q(\hat{\Glie}))$ is a highest weight representation $V$, that is to say there is $v_0\in V$ (highest weight vector) $\gamma_{i,m}^{\pm}\in\CC$ (highest weight) such that:
$$V=\U_q(\hat{\Glie}).v_0\text{ , }c^{\frac{1}{2}}.v_0=v_0$$
$$\forall i\in I,m\in\ZZ,x_{i,m}^+.v_0=0\text{ , }\text{ , }\phi_{i,m}^{\pm}.v_0=\gamma_{i,m}^{\pm}v_0$$
Moreover we have an $I$-uplet $(P_i(u))_{i\in I}$ of (Drinfeld-)polynomials such that $P_i(0)=1$ and:
$$\gamma_i^{\pm}(u)=\underset{m\in\NN}{\sum}\gamma_{i,\pm m}^{\pm}u^{\pm}=q_i^{\deg(P_i)}\frac{P_i(uq_i^{-1})}{P_i(uq_i)}\in\CC[[u^{\pm}]]$$
and $(P_i)_{i\in I}$ parameterizes simple modules in $\text{Rep}(\U_q(\hat{\Glie}))$.
\end{thm}

\begin{thm}({\bf Frenkel, Reshetikhin} \cite{Fre}) The eigenvalues $\gamma_i(u)^{\pm}\in\CC[[u]]$ of a representation $V\in\text{Rep}(\U_q(\hat{\Glie}))$ have the form:
$$\gamma_i^{\pm}(u)=q_i^{deg(Q_i)-deg(R_i)}\frac{Q_i(uq_i^{-1})R_i(uq_i)}{Q_i(uq_i)R_i(uq_i^{-1})}$$
where $Q_i(u),R_i(u)\in\CC[u]$ and $Q_i(0)=R_i(0)=1$.
\end{thm}

\noindent Note that the polynomials $Q_i,R_i$ are uniquely defined by $\gamma$. Denote by $Q_{\gamma,i}$, $R_{\gamma,i}$ the polynomials associated to $\gamma$.

\subsection{q-characters}\label{qcar}

Let $\Yim$ be the commutative ring $\Yim=\ZZ[Y_{i,a}^{\pm}]_{i\in I,a\in\CC^*}$. 

\begin{defi} For $V\in\text{Rep}(\U_q(\hat{\Glie}))$ a representation, the $q$-character $\chi_q(V)$\label{chiqdefi} of $V$ is:
$$\chi_q(V)=\underset{\gamma}{\sum}\text{dim}(V_{\gamma})\underset{i\in I,a\in\CC^*}{\prod}Y_{i,a}^{\lambda_{\gamma,i,a}-\mu_{\gamma,i,a}}\in \Yim$$   
where for $\gamma\in\CC^{I\times\ZZ}\times\CC^{I\times\ZZ}$, $i\in I$, $a\in\CC^*$ the $\lambda_{\gamma,i,a},\mu_{\gamma,i,a}\in\ZZ$ are defined by:
$$Q_{\gamma, i}(z)=\underset{a\in\CC^*}{\prod}(1-za)^{\lambda_{\gamma,i,a}}\text{ , }R_{\gamma, i}(z)=\underset{a\in\CC^*}{\prod}(1-za)^{\mu_{\gamma,i,a}}$$
\end{defi}

\begin{thm}({\bf Frenkel, Reshetikhin} \cite{Fre}) The map 
$$\chi_q:\text{Rep}(\U_q(\hat{\Glie}))\rightarrow \Yim$$ 
is an injective ring homomorphism and the following diagram is commutative:
$$\begin{array}{rcccl}
\text{Rep} (\U_q(\hat{\Glie}))&\stackrel{\chi_q}{\longrightarrow}&\ZZ[Y^{\pm}_{i,a}]_{i\in I,a\in\CC^*}\\
\downarrow res&&\downarrow\beta\\
\text{Rep}(\U_q(\Glie))&\stackrel{\chi}{\longrightarrow}&\ZZ[y^{\pm}_i]_{i\in I}\\\end{array}$$ 
where $\beta$ is the ring homomorphism such that $\beta(Y_{i,a})=y_i$ ($i\in I,a\in\CC^*$).
\end{thm}

\noindent For $m\in\Yim$ of the form $m=\underset{i\in I,a\in\CC^*}{\prod}Y_{i,a}^{u_{i,a}(m)}$ ($u_{i,a}(m)\geq 0$), denote $V_m\in\text{Rep}(\U_q(\hat{\Glie}))$ the simple module with Drinfeld polynomials $P_i(u)=\underset{a\in\CC^*}{\prod}(1-ua)^{u_{i,a}(m)}$. In particular for $i\in I,a\in\CC^*$ denote $V_{i,a}=V_{Y_{i,a}}$ and $X_{i,a}=\chi_q(V_{i,a})$. The simple modules $V_{i,a}$ are called fundamental representations.

\noindent Denote by $M_m\in\text{Rep}(\U_q(\hat{\Glie}))$ the module $M_m=\underset{i\in I,a\in\CC^*}{\bigotimes}V_{i,a}^{\otimes u_{i,a}(m)}$. It is called a standard module and his $q$-character is $\underset{i\in I,a\in\CC^*}{\prod}X_{i,a}^{u_{i,a}(m)}$.

\begin{cor}({\bf Frenkel, Reshetikhin} \cite{Fre}) The ring $\text{Rep}(\U_q(\hat{\Glie}))$ is commutative and isomorphic to $\ZZ[X_{i,a}]_{i\in I,a\in\CC^*}$.
\end{cor}

\begin{prop}({\bf Frenkel, Mukhin} \cite{Fre2})\label{aidafm} For $i\in I,a\in\CC^*$, we have $X_{i,a}\in\ZZ[Y_{j,aq^l}^{\pm}]_{j\in I,l\geq 0}$.\end{prop}

\noindent In particular for $a\in\CC^*$ we have an injective ring homomorphism:
$$\chi_q^a:\text{Rep}_a=\ZZ[X_{i,aq^l}]_{i\in I,l\in\ZZ}\rightarrow\Yim_a=\ZZ[Y_{i,aq^l}^{\pm}]_{i\in I,l\in\ZZ}$$
For $a,b\in\CC^*$ denote $\alpha_{b,a}:\text{Rep}_a\rightarrow\text{Rep}_b$ and $\beta_{b,a}:\Yim_a\rightarrow\Yim_b$ the canonical ring homomorphism. 
\begin{lem} We have a commutative diagram:
$$\begin{array}{rcccl}
\text{Rep}_a&\stackrel{\chi_q^a}{\longrightarrow}&\Yim_a\\
\alpha_{b,a}\downarrow &&\downarrow\beta_{b,a}\\
\text{Rep}_b&\stackrel{\chi_q^b}{\longrightarrow}&\Yim_b\\\end{array}$$\end{lem}
\noindent This result is a consequence of theorem \ref{simme} (or see \cite{Fre}, \cite{Fre2}). In particular it suffices to study $\chi_q^1$. In the following denote $\text{Rep}=\text{Rep}_1$, $X_{i,l}=X_{i,q^l}$\label{xil},  $\Yim=\Yim_1$ and $\chi_q=\chi_q^1:\text{Rep}\rightarrow \Yim$.\label{rep}

\section{Twisted polynomial algebras related to quantum affine algebras}\label{defoal} The aim of this section is to define the $t$-deformed algebra $\Yim_t$ and to describe its structure (theorem \ref{dessus}). We define the Heisenberg algebra $\mathcal{H}$, the subalgebra $\Yim_u\subset \mathcal{H}[[h]]$ and eventually $\Yim_t$ as a quotient of $\Yim_u$.

\subsection{Heisenberg algebras related to quantum affine algebras}

\subsubsection{The Heisenberg algebra $\mathcal{H}$} 
\begin{defi} $\mathcal{H}$\label{zq} is the $\CC$-algebra defined by generators $a_i[m]$\label{aim} ($i\in I, m\in\ZZ-\{0\}$), central elements $c_r$\label{cr} ($r>0$) and relations ($i,j\in I,m,r\in\ZZ-\{0\}$):
$$[a_i[m],a_j[r]]=\delta_{m,-r}(q^m-q^{-m})B_{i,j}(q^m)c_{|m|}$$
\end{defi}
\noindent This definition is motivated by the structure of $\U_q(\hat{\Glie})$: in $\mathcal{H}$ the $c_r$ are algebraically independent, but we have a surjective homomorphism from $\mathcal{H}$ to $\U_q(\hat{\Hlie})$ such that $a_i[m]\mapsto (q-q^{-1})h_{i,m}$ and $c_r\mapsto\frac{c^{r}-c^{-r}}{r}$.

\subsubsection{Properties of $\mathcal{H}$} For $j\in I,m\in\ZZ$ we set\label{yim}:
$$y_j[m]=\underset{i\in I}{\sum} \tilde{C}_{i,j}(q^m)a_i[m]\in \mathcal{H}$$

\begin{lem}\label{calczq} We have the Lie brackets in $\mathcal{H}$ ($i,j\in I,m,r\in\ZZ$):
$$[a_i[m],y_j[r]]=(q^{mr_i}-q^{-r_im})\delta_{m,-r}\delta_{i,j}c_{|m|}$$
$$[y_i[m],y_j[r]]=\delta_{m,-r}\tilde{C}_{j,i}(q^m)(q^{mr_j}-q^{-mr_j})c_{|m|}$$
\end{lem}

\demo We compute in $\mathcal{H}$:
$$[a_i[m],y_j[r]]
=[a_i[m],\underset{k\in I}{\sum} \tilde{C}_{k,j}(q^r)a_k[r]]
=\delta_{m,-r}c_{|m|}\underset{k\in I}{\sum}\tilde{C}_{k,j}(q^{-m})[r_i]_{q^m}C_{i,k}(q^m)(q^m-q^{-m})$$
$$=\delta_{i,j}\delta_{m,-r}(q^{mr_i}-q^{-mr_i})c_{|m|}$$
$$[y_i[m],y_j[r]]
=[\underset{k\in I}{\sum} \tilde{C}_{k,i}(q^m)a_k[m],y_j[r]]
=\delta_{m,-r}\tilde{C}_{j,i}(q^m)(q^{mr_j}-q^{-mr_j})c_{|m|}$$\qed

\noindent Let $\pi_+$\label{piplus} and $\pi_-$ be the $\CC$-algebra endomorphisms of $\mathcal{H}$ such that ($i\in I$, $m>0$, $r<0$):
$$\pi_+(a_i[m])=a_i[m]\text{ , }\pi_+(a_i[r])=0\text{ , }\pi_+(c_m)=0$$
$$\pi_-(a_i[m])=0\text{ , }\pi_-(a_i[r])=a_i[r]\text{ , }\pi_-(c_m)=0$$
They are well-defined because the relations are preserved. We set $\mathcal{H}^+=\text{Im}(\pi_+)\subset \mathcal{H}$\label{zqplus} and $\mathcal{H}^-=\text{Im}(\pi_-)\subset \mathcal{H}$.

\noindent Note that $\mathcal{H}^+$ (resp. $\mathcal{H}^-$) is the subalgebra of $\mathcal{H}$ generated by the $a_i[m]$, $i\in I,m>0$ (resp. $m<0$). So $\mathcal{H}^+$ and $\mathcal{H}^-$ are commutative algebras, and:
$$\mathcal{H}^+\simeq \mathcal{H}^-\simeq \CC[a_i[m]]_{i\in I,m>0}$$

\noindent We say that $m\in \mathcal{H}$ is a $\mathcal{H}$-monomial if it is a product of the generators $a_i[m],c_r$.

\begin{lem} There is a unique $\CC$-linear endomorphism $::$ of $\mathcal{H}$ such that for all $\mathcal{H}$-monomials $m$ we have:
$$:m:=\pi_+(m)\pi_-(m)$$
\end{lem}

\noindent In particular there is a vector space triangular decomposition $\mathcal{H}\simeq \mathcal{H}^+\otimes \CC[c_r]_{r>0}\otimes \mathcal{H}^-$.

\demo The $\mathcal{H}$-monomials span the $\CC$-vector space $\mathcal{H}$, so the map is unique. But there are non trivial linear combinations between them because of the relations of $\mathcal{H}$: it suffices to show that for $m_1$, $m_2$ $\mathcal{H}$-monomials the definition of $::$ is compatible with the relations ($i,j\in I$, $l,k\in\ZZ-\{0\}$):
$$m_1a_i[k]a_j[l]m_2-m_1a_j[l]a_i[k]m_2=\delta_{k,-l}(q^k-q^{-k})B_{i,j}(q^k)m_1c_{|k|}m_2$$
As $\mathcal{H}^+$ and $\mathcal{H}^-$ are commutative, we have:
$$\pi_+(m_1a_i[k]a_j[l]m_2)\pi_-(m_1a_i[k]a_j[l]m_2)=\pi_+(m_1a_j[l]a_i[k]m_2)\pi_-(m_1a_j[l]a_i[k]m_2)$$
and we can conclude because $\pi_+(m_1c_{|k|}m_2)=\pi_-(m_1c_{|k|}m_2)=0$.\qed

\subsection{The deformed algebra $\Yim_u$}

\subsubsection{Construction of $\Yim_u$} Consider the $\CC$-algebra $\mathcal{H}_h=\mathcal{H}[[h]]$\label{zqh}.  The application $\text{exp}$ is well-defined on the subalgebra $h\mathcal{H}_h$:
$$\text{exp}:h\mathcal{H}_h\rightarrow \mathcal{H}_h$$
For $l\in\ZZ$, $i\in I$, introduce $\tilde{A}_{i,l},\tilde{Y}_{i,l}\in \mathcal{H}_h$\label{tail} such that:
$$\tilde{A}_{i,l}=\text{exp}(\underset{m>0}{\sum}h^m a_i[m]q^{lm})\text{exp}(\underset{m>0}{\sum}h^m a_i[-m]q^{-lm})$$
$$\tilde{Y}_{i,l}=\text{exp}(\underset{m>0}{\sum}h^m y_{i}[m]q^{lm})\text{exp}(\underset{m>0}{\sum}h^m y_i[-m]q^{-lm})$$
Note that $\tilde{A}_{i,l}$ and $\tilde{Y}_{i,l}$ are invertible in $\mathcal{H}_h$ and that:
$$\tilde{A}_{i,l}^{-1}=\text{exp}(-\underset{m>0}{\sum}h^m a_i[-m]q^{-lm})\text{exp}(-\underset{m>0}{\sum}h^m a_i[m]q^{lm})$$
$$\tilde{Y}_{i,l}^{-1}=\text{exp}(-\underset{m>0}{\sum}h^m y_i[-m]q^{-lm})\text{exp}(-\underset{m>0}{\sum}h^m y_{i}[m]q^{lm})$$
Recall the definition $\mathfrak{U}\subset\QQ(z)$ of section \ref{recalu}. For $R\in \mathfrak{U}$, introduce $t_{R}\in \mathcal{H}_h$\label{tr}:
$$t_R=\text{exp}(\underset{m>0}{\sum}h^{2m}R(q^m)c_m)$$

\begin{defi}\label{yu} $\Yim_u$ is the $\ZZ$-subalgebra of $\mathcal{H}_h$ generated by the $\tilde{Y}_{i,l}^{\pm},\tilde{A}_{i,l}^{\pm},t_R$ ($i\in I,l\in\ZZ,R\in\mathfrak{U}$).\end{defi}

\noindent In this section we give properties of $\Yim_u$ and subalgebras of $\Yim_u$ which will be useful in section \ref{studyyt}.

\subsubsection{Relations in $\Yim_u$}\label{defij}

\begin{lem}\label{relu} We have the following relations in $\Yim_u$ ($i,j\in I$ $l,k\in\ZZ$):
\begin{equation}\label{ay}\tilde{A}_{i,l}\tilde{Y}_{j,k}\tilde{A}_{i,l}^{-1}\tilde{Y}_{j,k}^{-1}
=t_{\delta_{i,j}(z^{-r_i}-z^{r_i})(-z^{(l-k)}+z^{(k-l)})}\end{equation}
\begin{equation}\label{yy}\tilde{Y}_{i,l}\tilde{Y}_{j,k}\tilde{Y}_{i,l}^{-1}\tilde{Y}_{j,k}^{-1}
=t_{\tilde{C}_{j,i}(z)(z^{r_j}-z^{-r_j})(-z^{(l-k)}+z^{(k-l)})}\end{equation}
\begin{equation}\label{aa}\tilde{A}_{i,l}\tilde{A}_{j,k}\tilde{A}_{i,l}^{-1}\tilde{A}_{j,k}^{-1}=t_{B_{i,j}(z)(z^{-1}-z)(-z^{(l-k)}+z^{(k-l)})}\end{equation}
\end{lem}

\demo  For $A,B\in h\mathcal{H}_h$ such that $[A,B]\in h\CC[c_r]_{r>0}$, we have:
$$\text{exp}(A)\text{exp}(B)=\text{exp}(B)\text{exp}(A)\text{exp}([A,B])$$ 
So we can compute (see lemma \ref{calczq}):

$\tilde{A}_{i,l}\tilde{A}_{j,k}
\\=\text{exp}(\underset{m>0}{\sum}h^m a_{i}[m]q^{lm})(\text{exp}(\underset{m>0}{\sum}h^ma_i[-m]q^{-lm})\text{exp}(\underset{m>0}{\sum}h^ma_{j}[m]q^{km}))\text{exp}(\underset{m>0}{\sum}h^ma_j[-m]q^{-km})
\\=\text{exp}(\underset{m>0}{\sum}h^{2m}B_{i,j}(q^m)(q^{-m}-q^{m})q^{m(k-l)}c_m)
\\\text{exp}(\underset{m>0}{\sum}h^m a_{i}[m]q^{lm})\text{exp}(\underset{m>0}{\sum}h^m a_{j}[m]q^{km})\text{exp}(\underset{m>0}{\sum}h^ma_i[-m]q^{-lm})\text{exp}(\underset{m>0}{\sum}h^m a_j[-m]q^{-km})
\\=\text{exp}(\underset{m>0}{\sum}h^{2m}B_{i,j}(q^m)(q^{-m}-q^{m})(-q^{m(l-k)}+q^{m(k-l)})c_m)\tilde{A}_{j,k}\tilde{A}_{i,l}$

\mk

$\tilde{A}_{i,l}\tilde{Y}_{j,k}
\\=\text{exp}(\underset{m>0}{\sum}h^ma_{i}[m]q^{lm})(\text{exp}(\underset{m>0}{\sum}h^ma_i[-m]q^{-lm})\text{exp}(\underset{m>0}{\sum}h^my_j[m]q^{km}))\text{exp}(\underset{m>0}{\sum}h^my_j[-m]q^{-km})
\\=\text{exp}(\underset{m>0}{\sum}h^{2m}\delta_{i,j}(q^{-mr_i}-q^{mr_i})q^{m(k-l)}c_m)\text{exp}(\underset{m>0}{\sum}h^ma_{i}[m]q^{ml})
\\\text{exp}(\underset{m>0}{\sum}h^my_j[m]q^{mk})\text{exp}(\underset{m>0}{\sum}h^ma_i[-m]q^{-ml})\text{exp}(\underset{m>0}{\sum}h^my_j[-m]q^{-mk})
\\=\text{exp}(\underset{m>0}{\sum}h^{2m}\delta_{i,j}(q^{-mr_i}-q^{mr_i})(-q^{m(l-k)}+q^{m(k-l)})c_m)\tilde{Y}_{j,k}\tilde{A}_{i,l}$

\mk

$\tilde{Y}_{i,l}\tilde{Y}_{j,k}
\\=\text{exp}(\underset{m>0}{\sum}h^my_i[m]q^{ml})(\text{exp}(\underset{m>0}{\sum}h^my_i[-m]q^{-ml})\text{exp}(\underset{m>0}{\sum}h^my_j[m]q^{mk}))\text{exp}(\underset{m>0}{\sum}h^my_j[-m]q^{-mk})
\\=\text{exp}(\underset{m>0}{\sum}h^{2m}q^{m(k-l)}\tilde{C}_{j,i}(q^m)(q^{mr_j}-q^{-mr_j})c_m)\text{exp}(\underset{m>0}{\sum}h^my_{i}[m]q^{ml})
\\\text{exp}(\underset{m>0}{\sum}h^my_j[m]q^{mk})\text{exp}(\underset{m>0}{\sum}h^ma_i[-m]q^{-ml})\text{exp}(\underset{m>0}{\sum}h^my_j[-m]q^{-mk})
\\=\text{exp}(\underset{m>0}{\sum}h^{2m}\tilde{C}_{j,i}(q^m)(q^{mr_j}-q^{-mr_j})(-q^{m(l-k)}+q^{m(k-l)})c_m)\tilde{Y}_{j,k}\tilde{Y}_{i,l}$\qed

\subsubsection{Commutative subalgebras of $\mathcal{H}_h$}\label{defipip} The $\CC$-algebra endomorphisms $\pi_+,\pi_-$ of $\mathcal{H}$ are naturally extended to $\CC$-algebra endomorphisms of $\mathcal{H}_h$. As $\Yim_u\subset \mathcal{H}_h$, we have by restriction the $\ZZ$-algebra morphisms $\pi_{\pm}:\Yim_u\rightarrow \mathcal{H}_h$.

\noindent Introduce $\Yim=\pi_+(\Yim_u)\subset \mathcal{H}^+[[h]]$\label{y}. In this section \ref{defipip} we study $\Yim$. In particular we will see in proposition \ref{circ} that the notation $\Yim$ is consistent with the notation of section \ref{qcar}.

\noindent For $i\in I,l\in\ZZ$, denote\label{ail}:
$$Y_{i,l}^{\pm}=\pi_+(\tilde{Y}_{i,l}^{\pm})=\text{exp}(\pm\underset{m>0}{\sum}h^m y_{i}[m]q^{lm})$$ 
$$A_{i,l}^{\pm}=\pi_+(\tilde{A}_{i,l}^{\pm})=\text{exp}(\pm\underset{m>0}{\sum}h^m a_i[m]q^{lm})$$

\begin{lem}\label{gen} For $i\in I,l\in\ZZ$, we have:
$$A_{i,l}=Y_{i,l-r_i}Y_{i,l+r_i}(\underset{j/C_{j,i}=-1}{\prod}Y_{j,l}^{-1})(\underset{j/C_{j,i}=-2}{\prod}Y_{j,l+1}^{-1}Y_{j,l-1}^{-1})(\underset{j/C_{j,i}=-3}{\prod}Y_{j,l+2}^{-1}Y_{j,l}^{-1}Y_{j,l-2}^{-1})$$
In particular $\Yim$ is generated by the $Y_{i,l}^{\pm}$ ($i\in I,l\in\ZZ$).
\end{lem}

\demo

We have $a_i[m]=\underset{j\in I}{\sum}C_{j,i}(q^m)y_j[m]$, and:
$$\pi_+(\tilde{A}_{i,l})=\text{exp}(\underset{m>0}{\sum}h^m a_i[m]q^{lm})=\underset{j\in I}{\prod}\text{exp}(\underset{m>0}{\sum}h^m C_{j,i}(q^m)y_j[m]q^{lm})$$
As $C_{i,i}(q)=q^{r_i}+q^{-r_i}$, we have:
$$\text{exp}(\underset{m>0}{\sum}h^m C_{i,i}(q^m)y_i[m]q^{lm})=\text{exp}(\underset{m>0}{\sum}h^m y_i[m]q^{(l-r_i)m})\text{exp}(\underset{m>0}{\sum}h^m y_i[m]q^{(l+r_i)m})=Y_{i,l-r_i}Y_{i,l+r_i}$$
If $C_{j,i}<0$, we have $C_{j,i}(q)=-\underset{k=C_{j,i}+1, C_{j,i}+3...-C_{j,i}-1}{\sum}q^{k}$ and:
$$\text{exp}(-\underset{m>0}{\sum}h^m C_{j,i}(q^m)y_j[m]q^{lm})=\underset{k=C_{j,i}+1, C_{j,i}+3...-C_{j,i}-1}{\prod}\text{exp}(-\underset{m>0}{\sum}h^m y_j[m]q^{(l+k)m})$$
As $\Yim_u$ is generated by the $\tilde{Y}_{i,l}^{\pm},\tilde{A}_{i,l}^{\pm},t_R$ we get the last point.\qed

\noindent Note that the formula of lemma \ref{gen} already appeared in \cite{Fre}.

\noindent We need a general technical lemma to describe $\Yim$:

\begin{lem}\label{indgene} Let $J=\{1,...,r\}$ and let $\Lambda$ be the polynomial commutative algebra 
\\$\Lambda=\CC[\lambda_{j,m}]_{j\in J,m\geq 0}$. For $R=(R_1,...,R_r)\in\mathfrak{U}^{r}$, consider:
$$\Lambda_R=\text{exp}(\underset{j\in J,m>0}{\sum}h^mR_j(q^m)\lambda_{j,m})\in\Lambda[[h]]$$
Then the $(\Lambda_R)_{R\in\mathfrak{U}^r}$ are $\CC$-linearly independent. In particular the $\Lambda_{j,l}=\Lambda_{(0,...,0,z^l,0,...,0)}$ ($j\in J$, $l\in\ZZ$) are $\CC$-algebraically independent.\end{lem}

\demo Suppose we have a linear combination ($\mu_R\in\CC$, only a finite number of $\mu_R\neq 0$):
$$\underset{R\in\mathfrak{U}^r}{\sum}\mu_R\Lambda_R=0$$
The coefficients of $h^L$ in $\Lambda_R$ are of the form $R_{j_1}(q^{l_1})^{L_1}R_{j_2}(q^{l_2})^{L_2}...R_{j_N}(q^{l_N})^{L_N}\lambda_{j_1,l_1}^{L_1}\lambda_{j_2,l_2}^{L_2}...\lambda_{j_N,l_N}^{L_N}$ where $l_1L_1+...+l_NL_N=L$. So for $N\geq 0$, $j_1,...,j_N\in J$, $l_1,...,l_N>0$, $L_1,...,L_N\geq 0$ we have:
$$\underset{R\in\mathfrak{U}^r}{\sum}\mu_RR_{j_1}(q^{l_1})^{L_1}R_{j_2}(q^{l_2})^{L_2}...R_{j_N}(q^{l_N})^{L_N}=0$$
If we fix $L_2,...,L_{N}$, we have for all $L_1=l\geq 0$:
$$\underset{\alpha_1\in\CC}{\sum}\alpha_1^l\underset{R\in\mathfrak{U}^r/R_{j_1}(q^{l_1})=\alpha_1}{\sum}\mu_RR_{j_2}(q^{l_2})^{L_2}...R_{j_N}(q^{l_N})^{L_N}=0$$
We get a Van der Monde system which is invertible, so for all $\alpha_1\in\CC$:
$$\underset{R\in\mathfrak{U}^r/R_{j_1}(q^{l_1})=\alpha_1}{\sum}\mu_RR_{j_2}(q^{l_2})^{L_2}...R_{j_N}(q^{l_N})^{l_N}=0$$
By induction we get for $r'\leq N$ and all $\alpha_1,...,\alpha_{r'}\in\CC$: 
$$\underset{R\in\mathfrak{U}^r/R_{j_1}(q^{l_1})=\alpha_1,...,R_{j_{r'}}(q^{l_{r'}})=\alpha_{r'}}{\sum}\mu_{R}R_{j_{r'+1}}(q^{l_{r'+1}})^{L_{r'+1}}...R_{j_N}(q^{l_N})^{L_N}=0$$
And so for $r'=N$:
$$\underset{R\in\mathfrak{U}^r/R_{j_1}(q^{l_1})=\alpha_1,...,R_{j_{N}}(q^{l_{N}})=\alpha_{N}}{\sum}\mu_{R}=0$$
Let be $S\geq 0$ such that for all $\mu_R,\mu_{R'}\neq 0$, $j\in J$ we have $R_j-R_j'=0$ or $R_j-R_j'$ has at most $S-1$ roots. We set $N=Sr$ and $((j_1,l_1),...,(j_S,l_S))=((1,1),(1,2),...,(1,S),(2,1),...,(2,S),(3,1),...,(r,S))$. We get for all $\alpha_{j,l}\in\CC$ ($j\in J,1\leq l\leq S$):
$$\underset{R\in\mathfrak{U}^r/\forall j\in J,1\leq l\leq S, R_{j}(q^{l})=\alpha_{j,l}}{\sum}\mu_{R}=0$$
It suffices to show that there is at most one term is this sum. But consider $P,Q\in\mathfrak{U}$ such that for all $1\leq l\leq S$, $P(q^l)=P'(q^l)$. As $q$ is not a root of unity the $q^l$ are different and $P-P'$ has $S$ roots, so is $0$.

\noindent For the last assertion, we can write a monomial $\underset{j\in J,l\in\ZZ}{\prod}\Lambda_{j,l}^{u_{j,l}}=\Lambda_{\underset{l\in\ZZ}{\sum}u_{1,l}z^l,...,\underset{l\in\ZZ}{\sum}u_{r,l}z^l}$. In particular there is no trivial linear combination between those monomials.\qed

\noindent It follows from lemma \ref{gen} and lemma \ref{indgene}:

\begin{prop}\label{circ} The $Y_{i,l}\in\Yim$ are $\ZZ$-algebraically independent and generate the $\ZZ$-algebra $\Yim$. In particular, $\Yim$ is the commutative polynomial algebra $\ZZ[Y_{i,l}^{\pm}]_{i\in I,l\in\ZZ}$.

\noindent The $A_{i,l}^{-1}\in\Yim$ are $\ZZ$-algebraically independent. In particular the subalgebra of $\Yim$ generated by the $A_{i,l}^{-1}$ is the commutative polynomial algebra $\ZZ[A_{i,l}^{-1}]_{i\in I,l\in\ZZ}$.\end{prop}

\subsubsection{Generators of $\Yim_u$}\label{ptpt} 

\noindent The $\CC$-linear endomorphism $::$ of $\mathcal{H}$ is naturally extended to a $\CC$-linear endomorphism of $\mathcal{H}_h$. As $\Yim_u\subset \mathcal{H}_h$, we have by restriction a $\ZZ$-linear morphism $::$ from $\Yim_u$ to $\mathcal{H}_h$.

\noindent We say that $m\in\Yim_u$ is a $\Yim_u$-monomial if it is a product of generators $\tilde{A}_{i,l}^{\pm},\tilde{Y}_{i,l}^{\pm},t_R$.

\noindent In the following, for a product of non commuting terms, denote $\overset{\rightarrow}{\underset{s=1..S}{\prod}}U_s=U_1U_2...U_S$.

\begin{lem}\label{rel} The algebra $\Yim_u$ is generated by the $\tilde{Y}_{i,l}^{\pm},t_R$ ($i\in I,l\in\ZZ,R\in\mathfrak{U}$).
\end{lem}

\demo Let be $i\in I$, $l\in\ZZ$. It follows from proposition \ref{circ} that $\pi_+(\tilde{A}_{i,l})$ is of the form $\pi_+(\tilde{A}_{i,l})=\underset{i\in I,l\in\ZZ}{\prod}Y_{i,l}^{u_{i,l}}$ and that $:m:=:\overset{\rightarrow}{\underset{l\in\ZZ}{\prod}}\underset{i\in I}{\prod}\tilde{Y}_{i,l}^{u_{i,l}}:$. So it suffices to show that for $m$ a $\Yim_u$-monomial, there is a unique $R_m\in \mathfrak{U}$ such that $m=t_{R_m}:m:$. Let us write $m=t_R\overset{\rightarrow}{\underset{s=1.. S}{\prod}}U_s$ where $U_s\in\{\tilde{A}_{i,l}^{\pm},\tilde{Y}_{i,l}^{\pm}\}_{i\in I,l\in\ZZ}$ are generators. Then:
$$:m:=(\underset{s=1..S}{\prod}\pi_+(U_s))(\underset{s=1..S}{\prod}\pi_-(U_s))$$
And we can conclude because it follows from the proof of lemma \ref{relu} that for $1\leq s,s'\leq S$, there is $R_{s,s'}\in\mathfrak{U}$ such that $\pi_+(U_s)\pi_-(U_{s'})=t_{R_{s,s'}}\pi_-(U_{s'})\pi_+(U_s)$.\qed

\noindent In particular it follows from this proof that $:\Yim_u:\subset \Yim_u$.

\subsection{The deformed algebra $\Yim_t$}\label{studyyt}

\subsubsection{Construction of $\Yim_t$}\label{defipire} Denote by $\ZZ((z^{-1}))$ the ring of series of the form $P=\underset{r\leq R_P}{\sum}P_rz^r$ where $R_P\in\ZZ$ and the coefficients $P_r\in\ZZ$. Recall the definition $\mathfrak{U}$ of section \ref{recalu}. We have an embedding $\mathfrak{U}\subset\ZZ((z^{-1}))$ by expanding $\frac{1}{Q(z^{-1})}$ in $\ZZ[[z^{-1}]]$ for $Q(z)\in\ZZ[z]$ such that $Q(0)=1$. So we can introduce maps:\label{pir} 
$$\pi_r:\mathfrak{U}\rightarrow \ZZ\text{ , }P=\underset{k\leq R_P}{\sum}P_k z^k\mapsto P_r$$
Note that we could have consider the expansion in $\ZZ((z))$ and that the maps $\pi_r$ are not independent of our choice.

\begin{defi} We define  $\Yim_t$\label{tyt} (resp. $\mathcal{H}_t$)\label{zqt} as the algebra quotient of $\Yim_u$ (resp. $\mathcal{H}_h$) by relations: 
$$t_R=t_{R'}\text{ if $\pi_0(R)=\pi_0(R')$}$$\end{defi}
\noindent We keep the notations $\tilde{Y}_{i,l}^{\pm},\tilde{A}_{i,l}^{\pm}$ for their image in $\Yim_t$. Denote by $t$ the image of $t_1=\text{exp}(\underset{m>0}{\sum}h^{2m}c_m)$ in $\Yim_t$. As $\pi_0$ is additive, the image of $t_R$ in $\Yim_t$ is $t^{\pi_o(R)}$\label{t}. In particular $\Yim_t$ is generated by the $\tilde{Y}_{i,l}^{\pm},\tilde{A}_{i,l}^{\pm},t^{\pm}$.

\noindent As the defining relations of $\mathcal{H}_t$ involve only the $c_l$ and $\pi_+(c_l)=\pi_-(c_l)=0$, the algebra endomorphisms $\pi_+,\pi_-$ of $\mathcal{H}_t$ are well-defined. So we can define\label{zqtplus} $\mathcal{H}_t^+,\mathcal{H}_t^-,\Yim_t^+,\Yim_t^-$\label{tytplus} in the same way as in section \ref{defipip} and $::$ a $\CC$-linear endomorphism of $\mathcal{H}_t$ as in section \ref{ptpt}. The $\ZZ[t^{\pm}]$-subalgebra $\Yim_t\subset \mathcal{H}_t$ verifies $:\Yim_t:\subset\Yim_t$ (proof of lemma \ref{rel}). We have $\Yim_t^+\simeq\Yim$.

\noindent We say that $m\in\Yim_t$ (resp. $m\in\Yim$) is a $\Yim_t$-monomial (resp. a $\Yim$-monomial) if it is a product of the generators $\tilde{Y}_{i,m}^{\pm},t^{\pm}$ (resp. $Y_{i,m}^{\pm}$).

\subsubsection{Structure of $\Yim_t$}\label{dessusdeux}

The following theorem gives the structure of $\Yim_t$:

\begin{thm}\label{dessus} The algebra $\Yim_t$ is defined by generators $\tilde{Y}_{i,l}^{\pm}$ $(i\in I,l\in\ZZ)$, central elements $t^{\pm}$ and relations ($i,j\in I, k,l\in\ZZ$):
$$\tilde{Y}_{i,l}\tilde{Y}_{j,k}=t^{\gamma(i,l,j,k)}\tilde{Y}_{j,k}\tilde{Y}_{i,l}$$
where $\gamma: (I\times\ZZ)^2\rightarrow\ZZ$ is given by (recall the maps $\pi_r$ of section \ref{defipire})\label{gamma}:
$$\gamma(i,l,j,k)=\underset{r\in\ZZ}{\sum}\pi_r(\tilde{C}_{j,i}(z))(-\delta_{l-k,-r_j-r}-\delta_{l-k,r-r_j}+\delta_{l-k,r_j-r}+\delta_{l-k,r_j+r})$$
\end{thm}

\demo As the image of $t_R$ in $\Yim_t$ is $t^{\pi_o(R)}$, we can deduce the relations from lemma \ref{relu}. For example formula \ref{yy} (p. \pageref{yy}) gives:
$$\tilde{Y}_{i,l}\tilde{Y}_{j,k}\tilde{Y}_{i,l}^{-1}\tilde{Y}_{j,k}^{-1}=t^{\pi_0((\tilde{C}_{j,i}(z)(z^{r_j}-z^{-r_j})(-z^{(l-k)}+z^{(k-l)}))}$$
where:
$$\pi_0(\tilde{C}_{j,i}(z)(z^{r_j}-z^{-r_j})(-z^{(l-k)}+z^{(k-l)}))$$
$$=\underset{r\in\ZZ}{\sum}\pi_r(\tilde{C}_{j,i}(z))(\delta_{r_j+r+k-l,0}+\delta_{-r_j+r+l-k,0}-\delta_{r_j+r+l-k,0}-\delta_{-r_j+r+k-l,0})=\gamma(i,l,j,k)$$
It follows from lemma \ref{gen} that $\Yim_t$ is generated by the $\tilde{Y}_{i,l}^{\pm},t^{\pm}$. 

\noindent It follows from lemma \ref{indgene} that the $t_R\in\Yim_u$ ($R\in\mathfrak{U}$) are $\ZZ$-linearly independent. So the $\ZZ$-algebra $\ZZ[t_R]_{R\in\mathfrak{U}}$ is defined by generators $(t_R)_{R\in\mathfrak{U}}$ and relations $t_{R+R'}=t_Rt_{R'}$ for $R,R'\in\mathfrak{U}$. In particular the image of $\ZZ[t_R]_{R\in\mathfrak{U}}$ in $\Yim_t$ is $\ZZ[t^{\pm}]$. 

\noindent Let $A$ be the classes of $\Yim_t$-monomials modulo $t^{\ZZ}$. So we have:
$$\underset{m\in A}{\sum}\ZZ[t^{\pm}].m=\Yim_t$$
We prove the sum is direct: suppose we have a linear combination $\underset{m\in A}{\sum}\lambda_m(t)m=0$ where $\lambda_m(t)\in\ZZ[t^{\pm}]$. We saw in proposition \ref{circ} that $\Yim\simeq \ZZ[{Y}_{i,l}^{\pm}]_{i\in I,l\in\ZZ}$. So $\lambda_m(1)=0$ and $\lambda_m(t)=(t-1)\lambda_m^{(1)}(t)$ where $\lambda_m^{(1)}(t)\in\ZZ[t^{\pm}]$. In particular $\underset{m\in A}{\sum}\lambda_m(t)^{(1)}(t)m=0$ and we get by induction $\lambda_m(t)\in (t-1)^r\ZZ[t^{\pm}]$ for all $r\geq 0$. This is possible if and only if all $\lambda_m(t)=0$.\qed

\noindent In the same way using the last assertion of proposition \ref{circ}, we have:

\begin{prop}\label{yenga} The sub $\ZZ[t^{\pm}]$-algebra of $\Yim_t$ generated by the $\tilde{A}_{i,l}^{-1}$ is defined by generators $\tilde{A}_{i,l}^{-1},t^{\pm}$ $(i\in I,l\in\ZZ)$ and relations\label{alpha}:
$$\tilde{A}_{i,l}^{-1}\tilde{A}_{j,k}^{-1}=t^{\alpha(i,l,j,k)}\tilde{A}_{j,k}^{-1}\tilde{A}_{i,l}^{-1}$$
where $\alpha: (I\times\ZZ)^2\rightarrow\ZZ$ is given by:
$$\alpha(i,l,i,k)=2(-\delta_{l-k,2r_i}+\delta_{l-k,-2r_i})$$
$$\alpha(i,l,j,k)=2\underset{r=C_{i,j}+1,C_{i,j}+3,...,-C_{i,j}-1}{\sum}(-\delta_{l-k,-r_i+r}+\delta_{l-k,r_i+r})\text{ (if $i\neq j$)}$$\end{prop}

\noindent Moreover we have the following relations in $\Yim_t$:
$$\tilde{A}_{i,l}\tilde{Y}_{j,k}=t^{\beta(i,l,j,k)}\tilde{Y}_{j,k}\tilde{A}_{i,l}$$
where $\beta: (I\times\ZZ)^2\rightarrow\ZZ$ is given by\label{beta}:
$$\beta(i,l,j,k)=2\delta_{i,j}(-\delta_{l-k,r_i}+\delta_{l-k,-r_i})$$

\subsection{Notations and properties related to monomials} In this section we study some technical properties of the $\Yim$-monomials and the $\Yim_t$-monomials which will be used in the following.

\subsubsection{Basis} Denote by $A$\label{a} the set of $\Yim$-monomials. It is a $\ZZ$-basis of $\Yim$ (proposition \ref{circ}). Let us define an analog $\ZZ[t^{\pm}]$-basis of $\Yim_t$: denote $A'$\label{ap} the set of $\Yim_t$-monomials of the form $m=:m:$. It follows from theorem \ref{dessus} that:
$$\Yim_t=\underset{m\in A'}{\bigoplus}\ZZ[t^{\pm}]m$$
The map $\pi:A'\rightarrow A$ defined by $\pi(m)=\pi_+(m)$\label{pi} is a bijection. In the following we identify $A$ and $A'$. In particular we have an embedding $\Yim\subset\Yim_t$ and an isomorphism of $\ZZ[t^{\pm}]$-modules $\Yim\otimes_{\ZZ}\ZZ[t^{\pm}]\simeq\Yim_t$. Note that it depends on the choice of the $\ZZ[t^{\pm}]$-basis of $\Yim_t$. 

\noindent We say that $\chi_1\in\Yim_t$ has the same monomials as $\chi_2\in \Yim$ if in the decompositions $\chi_1=\underset{m\in A}{\sum}\lambda_m(t)m$, $\chi_2=\underset{m\in A}{\sum}\mu_mm$ we have $\lambda_m(t)=0\Leftrightarrow\mu_m=0$.

\subsubsection{The notation $u_{i,l}$}\label{notuil} For $m$ a $\Yim$-monomial we set $u_{i,l}(m)\in\ZZ$\label{uil} such that $m=\underset{i\in I,l\in\ZZ}{\prod}{Y}_{i,l}^{u_{i,l}(m)}$ and $u_i(m)=\underset{l\in\ZZ}{\sum}u_{i,l}(m)$\label{ui}. For $m$ a $\Yim_t$-monomial, we set $u_{i,l}(m)=u_{i,l}(\pi_+(m))$ and $u_i(m)=u_i(\pi_+(m))$. Note that $u_{i,l}$ is invariant by multiplication by $t$ and compatible with the identification of $A$ and $A'$.

\noindent Note that section \ref{dessusdeux} implies that for $i\in I, l\in \ZZ$ and $m$ a $\Yim_t$-monomial we have:
$$\tilde{A}_{i,l}m=t^{-2u_{i,l-r_i}(m)+2u_{i,l+r_i}(m)}m\tilde{A}_{i,l}$$

\noindent Denote by $B_i\subset A$ the set of $i$-dominant $\Yim$-monomials, that is to say $m\in B_i$\label{bi} if $\forall l\in\ZZ$, $u_{i,l}(m)\geq 0$. For $J\subset I$ denote $B_J=\underset{i\in J}{\bigcap}B_i$\label{bj} the set of $J$-dominant $\Yim$-monomials. In particular, $B=B_I$ is the set of dominant $\Yim$-monomials.\label{b}

\noindent We recall we can define a partial ordering on $A$ by putting $m\leq m'$ if there is a $\Yim$-monomial $M$ which is a product of $A_{i,l}^{\pm}$ ($i\in I,l\in\ZZ$) such that $m=Mm'$ (see for example \cite{Her01}). A maximal (resp. lowest, higher...) weight $\Yim$-monomial is a maximal (resp. minimal, higher...) element of $A$ for this ordering. We deduce from $\pi_+$ a partial ordering on the $\Yim_t$-monomials.

\noindent Following \cite{Fre2}, a $\Yim$-monomial $m$ is said to be right negative if the factors $Y_{j,l}$ appearing in $m$, for which $l$ is maximal, have negative powers. A product of right negative $\Yim$-monomials is right negative. It follows from lemma \ref{gen} that the $A_{i,l}^{-1}$ are right negative. A $\Yim_t$-monomial is said to be right negative if $\pi_+(m)$ is right negative.

\subsubsection{Some technical properties}

\begin{lem}\label{genea} Let $(i_1,l_1),...,(i_K,l_K)$ be in $(I\times\ZZ)^K$. For $U\geq 0$, the set of the $m=\underset{k=1...K}{\prod}A_{i_k,l_k}^{-v_{i_k,l_k}(m)}$ ($v_{i_k,l_k}(m)\geq 0$) such that $\underset{i\in I,k\in\ZZ}{\text{min}}u_{i,k}(m)\geq -U$ is finite.\end{lem}

\demo Suppose it is not the case: let be $(m_p)_{p\geq 0}$ such that $\underset{i\in I,k\in\ZZ}{\text{min}}u_{i,k}(m_p)\geq -U$ but 
\\$\underset{k=1...K}{\sum}v_{i_k,l_k}(m_p)\underset{p\rightarrow\infty}{\rightarrow}+\infty$. So there is at least one $k$ such that $v_{i_k,l_k}(m_p)\underset{p\rightarrow\infty}{\rightarrow}+\infty$. Denote by $\mathfrak{R}$ the set of such $k$. Among those $k\in\mathfrak{R}$, such that $l_k$ is maximal suppose that $r_{i_k}$ is maximal (recall the definition of $r_i$ in section \ref{recalu}). In particular, we have $u_{i_k,l_k+r_{i_k}}(m_p)=-v_{i_k,l_k}(m_p)+f(p)$ where $f(p)$ depends only of the $v_{i_{k'},l_{k'}}(m_p)$, $k'\notin\mathfrak{R}$. In particular, $f(p)$ is bounded and $u_{i_k,l_k+r_{i_k}}(m_p)\underset{p\rightarrow\infty}{\rightarrow}-\infty$.\qed

\begin{lem}\label{fini} For $M\in B$, $K\geq 0$ the set of $\Yim$-monomials $\{MA_{i_1,l_1}^{-1}...A_{i_R,l_R}^{-1}/R\geq 0,l_1,...,l_R\geq K\}\cap B$ is finite. \end{lem}

\demo Let us write $M=Y_{i_1,l_1}...Y_{i_R,l_R}$ such that $l_1=\underset{r=1...R}{\text{min}}l_r$, $l_R=\underset{r=1...R}{\text{max}}l_r$ and consider $m$ in the set. It is of the form $m=MM'$ where $M'=\underset{i\in I,l\geq K}{\prod}A_{i,l}^{-v_{i,l}}$ ($v_{i,l}\geq 0$). Let $L=\text{max}\{l\in\ZZ/\exists i\in I, u_{i,l}(M')<0\}$. $M'$ is right negative so for all $i\in I$, $l>L\Rightarrow v_{i,l}=0$. But $m$ is dominant, so $L\leq l_R$. In particular $M'=\underset{i\in I,K\leq l\leq l_R}{\prod}A_{i,l}^{-v_{i,l}}$. It suffices to prove that the $v_{i,l}(m_r)$ are bounded under the condition $m$ dominant. This follows from lemma \ref{genea}.\qed

\subsection{Presentations of deformed algebras} Our construction of $\Yim_t$ using $\mathcal{H}_h$ (section \ref{studyyt}) is a ``concrete'' presentation of the deformed structure. Let us look at another approach: in this section we define two bicharacters $\mathcal{N},\mathcal{N}_t$ related to basis of $\Yim_t$. All the information of the multiplication of $\Yim_t$ is contained is those bicharacters because we can construct a deformed $*$ multiplication on the ``abstract'' $\ZZ[t^{\pm}]$-module $\Yim\otimes_{\ZZ}\ZZ[t^{\pm}]$ by putting for $m_1,m_2\in A$ $\Yim$-monomials:
$$m_1*m_2=t^{\mathcal{N}(m_1,m_2)-\mathcal{N}(m_2,m_1)}m_2*m_1$$
or
$$m_1*m_2=t^{\mathcal{N}_t(m_1,m_2)-\mathcal{N}_t(m_2,m_1)}m_2*m_1$$
Those presentations appeared earlier in the literature \cite{Nab}, \cite{Vas} for the simply laced case. In particular this section identifies our approach with those articles and gives an algebraic motivation of the deformed structures of \cite{Nab}, \cite{Vas} related to the structure of $\U_q(\hat{\Glie})$.

\subsubsection{The bicharacter $\mathcal{N}$} It follows from the proof of lemma \ref{rel} that for $m$ a $\Yim_t$-monomial, there is $N(m)\in\ZZ$\label{n} such that $m=t^{N(m)}:m:$. For $m_1,m_2$ $\Yim_t$-monomials we define $\mathcal{N}(m_1,m_2)=N(m_1m_2)-N(m_1)-N(m_2)$. We have $N(Y_{i,l})=N(A_{i,l})=0$. Note that for $\alpha,\beta\in\ZZ$ we have:
$$N(t^{\alpha}m)=\alpha+N(m)\text{ , }\mathcal{N}(t^{\alpha}m_1,t^{\beta}m_2)=\mathcal{N}(m_1,m_2)$$
In particular the map $\mathcal{N}:A\times A\rightarrow\ZZ$ is well-defined and independent of the choice of a representant in $\pi_+^{-1}(A)$.
\begin{lem}\label{aaa} For $m_1,m_2$ $\Yim_t$-monomials, we have in $\mathcal{H}_t$:
$$\pi_-(m_1)\pi_+(m_2)=t^{\mathcal{N}(m_1,m_2)}\pi_+(m_2)\pi_-(m_1)$$
\end{lem}

\demo We have:
$$m_1=t^{N(m_1)}\pi_+(m_1)\pi_-(m_1)\text{ , }m_2=t^{N(m_2)}\pi_+(m_2)\pi_-(m_2)$$
and so:
$$m_1m_2=t^{N(m_1m_2)}\pi_+(m_1)\pi_+(m_2)\pi_-(m_1)\pi_-(m_2)=t^{N(m_1)+N(m_2)}\pi_+(m_1)\pi_-(m_1)\pi_+(m_2)\pi_-(m_2)$$
\qed

\begin{lem}\label{bibi} The map $\mathcal{N}:A\times A\rightarrow \ZZ$ is a bicharacter, that is to say for $m_1,m_2,m_3\in A$, we have:
$$\mathcal{N}(m_1m_2,m_3)=\mathcal{N}(m_1,m_3)+\mathcal{N}(m_2,m_3)\text{ and }\mathcal{N}(m_1,m_2m_3)=\mathcal{N}(m_1,m_2)+\mathcal{N}(m_1,m_3)$$
Moreover for $m_1,...,m_k$ $\Yim_t$-monomials, we have:
$$N(m_1m_2...m_k)=N(m_1)+N(m_2)+...+N(m_k)+\underset{1\leq i<j\leq k}{\sum}\mathcal{N}(m_i,m_j)$$
\end{lem}

\demo For the first point it follows from lemma \ref{aaa}:
$$\pi_-(m_1m_2)\pi_+(m_3)=t^{\mathcal{N}(m_1m_2,m_3)}\pi_+(m_3)\pi_-(m_1m_2)=t^{\mathcal{N}(m_2,m_3)}\pi_-(m_1)\pi_+(m_3)\pi_-(m_2)$$
$$=t^{\mathcal{N}(m_1,m_3)+\mathcal{N}(m_2,m_3)}\pi_+(m_3)\pi_-(m_1m_2)$$
For the second point we have first:
$$N(m_1m_2)=N(m_1)+N(m_2)+\mathcal{N}(m_1,m_2)$$
and by induction:
$$N(m_1m_2...m_k)=N(m_1)+N(m_2...m_k)+\mathcal{N}(m_1,m_2...m_k)$$
$$=N(m_1)+N(m_2)+...+N(m_k)+\underset{1<i<j\leq k}{\sum}\mathcal{N}(m_i,m_j)+\mathcal{N}(m_1,m_2)+...+\mathcal{N}(m_1,m_k)$$
\qed

\subsubsection{The bicharacter $\mathcal{N}_t$}\label{tildeat}

For $m$ a $\Yim_t$-monomial and $l\in\ZZ$, denote $\pi_l(m)=\underset{j\in I}{\prod}\tilde{Y}_{j,l}^{u_{j,l}(m)}$. It is well defined because for $i,j\in I$ and $l\in\ZZ$ we have $\tilde{Y}_{i,l}\tilde{Y}_{j,l}=\tilde{Y}_{j,l}\tilde{Y}_{i,l}$ (theorem \ref{dessus}). Moreover for $m_1,m_2$ $\Yim_t$-monomials we have $\pi_l(m_1m_2)=\pi_l(m_1)\pi_l(m_2)=\underset{i\in I}{\prod}\tilde{Y}_{i,l}^{u_{i,l}(m_1)+u_{i,l}(m_2)}$.

\noindent For $m$ a $\Yim_t$-monomial denote $\tilde{m}={\underset{l\in\ZZ}{\overset{\rightarrow}{\prod}}}\pi_l(m)$\label{tm}, and $A_t$\label{tat} the set of $\Yim_t$-monomials of the form $\tilde{m}$. From theorem \ref{dessus} there is a unique $N_t(m)\in\ZZ$\label{nt} such that $m=t^{N_t(m)}\tilde{m}$, and:
$$\Yim_t=\underset{m\in A_t}{\bigoplus}\ZZ[t^{\pm}]m$$
For $m_1,m_2$ $\Yim_t$-monomials we define $\mathcal{N}_t(m_1,m_2)=N_t(m_1m_2)-N_t(m_1)-N_t(m_2)$. We have $N_t(Y_{i,l})=0$. Note that for $\alpha,\beta\in\ZZ$ we have:
$$N_t(t^{\alpha}m)=\alpha+N_t(m)\text{ , }\mathcal{N}_t(t^{\alpha}m_1,t^{\beta}m_2)=\mathcal{N}_t(m_1,m_2)$$
In particular the map $\mathcal{N}_t:A\times A\rightarrow\ZZ$ is well-defined and independent of the choice of $A$.
\begin{lem} For $m_1,m_2$ $\Yim_t$-monomials, we have:
$$\mathcal{N}_t(m_1,m_2)=\underset{l>l'}{\sum}(\mathcal{N}(\pi_l(m_1),\pi_{l'}(m_2))-\mathcal{N}(\pi_{l'}(m_2),\pi_{l}(m_1)))$$
In particular, $\mathcal{N}_t$ is a bicharacter and for $m_1,...,m_k$ $\Yim_t$-monomials, we have:
$$N_t(m_1m_2...m_k)=N_t(m_1)+N_t(m_2)+...+N_t(m_k)+\underset{1\leq i<j\leq k}{\sum}\mathcal{N}_t(m_i,m_j)$$
\end{lem}

\demo For the first point, it follows from the definition that $\tilde{(m_1m_2)}=t^{\mathcal{N}_t(m_1,m_2)}\tilde{m_1}\tilde{m_2}$. But:
$$\tilde{(m_1m_2)}=\overset{\rightarrow}{\underset{l\in\ZZ}{\prod}}\pi_l(m_1)\pi_l(m_2)\text{ , }\tilde{m_1}\tilde{m_2}=(\overset{\rightarrow}{\underset{l\in\ZZ}{\prod}}\pi_l(m_1))(\overset{\rightarrow}{\underset{l\in\ZZ}{\prod}}\pi_l(m_2))$$
So we have to commute $\pi_l(m_1)$ and $\pi_{l'}(m_2)$ for $l>l'$. The last assertion is proved as in lemma \ref{bibi}.\qed

\subsubsection{Presentation related to the basis $A_t$ and identification with \cite{Nab}}

We suppose we are in the $ADE$-case. 

\noindent Let be $m_1=:\underset{i\in I,l\in\ZZ}{\prod}\tilde{Y}_{i,l}^{y_{i,l}}\tilde{A}_{i,l}^{-v_{i,l}}:,m_2=:\underset{i\in I,l\in\ZZ}{\prod}\tilde{Y}_{i,l}^{y_{i,l}'}\tilde{A}_{i,l}^{-v_{i,l}'}:\in\Yim_t$. We set $m_1^y=:\underset{i\in I,l\in\ZZ}{\prod}\tilde{Y}_{i,l}^{y_{i,l}}:$ and $m_2^y=:\underset{i\in I,l\in\ZZ}{\prod}\tilde{Y}_{i,l}^{y_{i,l}'}:$.

\begin{prop}\label{form} We have $\mathcal{N}_t(m_1,m_2)=\mathcal{N}_t(m_1^y,m_2^y)+2d(m_1,m_2)$, where:\label{d}
$$d(m_1,m_2)=\underset{i\in I,l\in\ZZ}{\sum}v_{i,l+1}u_{i,l}'+y_{i,l+1}v_{i,l}'=\underset{i\in I,l\in\ZZ}{\sum}u_{i,l+1}v_{i,l}'+v_{i,l+1}y_{i,l}'$$
where $u_{i,l}=y_{i,l}-v_{i,l-1}-v_{i,l+1}+\underset{j/C_{i,j}=-1}{\sum}v_{j,l}$ and $u_{i,l}'=y_{i,l}'-v_{i,l-1}'-v_{i,l+1}'+\underset{j/C_{i,j}=-1}{\sum}v_{j,l}'$.
\end{prop}

\demo

First notice that we have ($i\in I,l\in\ZZ$):
$$\mathcal{N}_t(Y_{i,l},A_{i,l-1}^{-1})=2\text{ , }\mathcal{N}_t(A_{i,l+1}^{-1},Y_{i,l})=2\text{ , }\mathcal{N}_t(A_{i,l+1}^{-1},A_{i,l-1}^{-1})=-2$$
$$\mathcal{N}_t(Y_{i,l+1}^{-1},Y_{i,l-1}^{-1})=-2\text{ , }\mathcal{N}_t(A_{i,l+1}^{-1},Y_{i,l})=2$$
For example $\mathcal{N}_t(Y_{i,l},A_{i,l-1}^{-1})=\mathcal{N}(Y_{i,l},A_{i,l-1}^{-1})-\mathcal{N}(A_{i,l-1}^{-1},Y_{i,l})=2$ because $\tilde{Y}_{i,l}\tilde{A}_{i,l-1}^{-1}=t^2\tilde{A}_{i,l-1}^{-1}\tilde{Y}_{i,l}$.

\noindent We have $\mathcal{N}_t(m_1,m_2)=A+B+C+D$ where:

$A=\mathcal{N}_t(m_1^y,m_2^y)$

$B=\underset{i,j\in I,l,k\in\ZZ}{\sum}y_{i,l}v_{j,k}'\mathcal{N}_t(Y_{i,l},A_{j,k}^{-1})=\underset{i\in I,l\in\ZZ}{\sum}y_{i,l}v_{i,l-1}'\mathcal{N}_t(Y_{i,l},A_{i,l-1}^{-1})=2\underset{i\in I,l\in\ZZ}{\sum}y_{i,l}v_{i,l-1}'$

$C=\underset{i,j\in I,l,k\in\ZZ}{\sum}v_{i,l}y_{j,k}'\mathcal{N}_t(A_{i,l}^{-1},Y_{j,k})=\underset{i\in I,l\in\ZZ}{\sum}v_{i,l+1}y_{i,l}'\mathcal{N}_t(A_{i,l+1}^{-1},Y_{i,l})=2\underset{i\in I,l\in\ZZ}{\sum}v_{i,l+1}y_{i,l}'$

$D=\underset{i,j\in I,l,k\in\ZZ}{\sum}v_{i,l}v_{j,k}'\mathcal{N}_t(A_{i,l}^{-1},A_{j,k}^{-1})
\\=\underset{i\in I,l\in\ZZ}{\sum}v_{i,l+1}v_{i,l-1}'\mathcal{N}_t(A_{i,l+1}^{-1},A_{i,l-1}^{-1})+v_{i,l}v_{i,l}'\mathcal{N}_t(Y_{i,l+1}^{-1},Y_{i,l-1}^{-1})+\underset{C_{j,i}=-1,l\in\ZZ}{\sum}v_{i,l+1}v_{j,l}'\mathcal{N}_t(A_{i,l+1}^{-1},Y_{i,l})
\\=-2\underset{i\in I,l\in\ZZ}{\sum}(v_{i,l+1}v_{i,l-1}'+v_{i,l}v_{i,l'})+2\underset{C_{j,i}=-1,l\in\ZZ}{\sum}v_{i,l+1}v_{j,l}'$

\noindent In particular, we have:
$$B+C+D=2\underset{i\in I,l\in\ZZ}{\sum}(y_{i,l}v_{i,l-1}'+v_{i,l+1}y_{i,l}'-v_{i,l+1}v_{i,l-1}'-v_{i,l}v_{i,l}')+2\underset{C_{j,i}=-1,l\in\ZZ}{\sum}v_{i,l+1}v_{j,l}'$$\qed

\noindent The bicharacter $d$ was introduced for the $ADE$-case by Nakajima in \cite{Nab} motivated by geometry. It particular this proposition \ref{form} gives a new motivation for this deformed structure. 

\subsubsection{Presentation related to the basis $A$ and identification with \cite{Vas}}\label{varva}

\begin{lem} For $m_1,m_2\in A$, we have:
$$\mathcal{N}(m_1,m_2)=\underset{i,j\in I, l,k\in\ZZ}{\sum}u_{i,l}(m_1)u_{j,k}(m_2)((\tilde{C}_{j,i}(z))_{r_j+l-k}-(\tilde{C}_{j,i}(z))_{-r_j+l-k})$$
\end{lem}

\demo First we can compute in $\Yim_u$:
$$\tilde{Y}_{i,l}\tilde{Y}_{j,k}=\text{exp}(\underset{m>0}{\sum}h^{2m}[y_i[-m],y_j[m]]q^{m(k-l)}):\tilde{Y}_{i,l}\tilde{Y}_{j,k}:=t_{\tilde{C}_{j,i}(z)z^{k-l}(z^{-r_j}-z^{r_j})}:\tilde{Y}_{i,l}\tilde{Y}_{j,k}:$$
and as $N(\tilde{Y}_{i,l})=N(\tilde{Y}_{j,k})=0$ we have $\mathcal{N}(\tilde{Y}_{i,l},\tilde{Y}_{j,k})=(\tilde{C}_{j,i}(z))_{r_j+l-k}-(\tilde{C}_{j,i}(z))_{-r_j+l-k}$.\qed

\noindent In $sl_2$-case we have $C(z)=z+z^{-1}$ and $\tilde{C}(z)=\frac{1}{z+z^{-1}}=\underset{r\geq 0}{\sum}(-1)^rz^{-2r-1}$. So:
$$\tilde{Y}_l\tilde{Y}_k=t^s:\tilde{Y}_l\tilde{Y}_k:$$
where:

$s=0$ if $l-k=1+2r$, $r\in\ZZ$

$s=0$ if $l-k=2r$, $r>0$

$s=2(-1)^{r+1}$ if $l-k=2r$, $r<0$

$s=-1$ if $l=k$

\noindent It is analogous to the multiplication introduced for the $ADE$-case by Varagnolo-Vasserot in \cite{Vas}: we suppose we are in the $ADE$-case, denote $P=\underset{i\in I}{\bigoplus}\ZZ \omega_i$ (resp. $Q=\underset{i\in I}{\bigoplus}\ZZ \alpha_i$) the weight-lattice (resp. root-lattice) and:

$\bar{}:P\otimes\ZZ[z^{\pm}]\rightarrow P\otimes\ZZ[z^{\pm}]$ is defined by $\overline{\lambda \otimes P(z)}=\lambda\otimes P(z^{-1})$.

$(,):Q\otimes \ZZ((z^{-1}))\times P\otimes\ZZ((z^{-1}))\rightarrow \ZZ((z^{-1}))$ is the $\ZZ((z^{-1}))$-bilinear form defined by $(\alpha_i,\omega_j)=\delta_{i,j}$.

$\Omega^{-1}:P\otimes\ZZ[z^{\pm}]\rightarrow Q\otimes\ZZ((z^{-1}))$ is defined by $\Omega^{-1}(\omega_i)=\underset{k\in I}{\sum}\tilde{C}_{i,k}(z)\alpha_k$.

\noindent The map $\epsilon:P\otimes\ZZ[z^{\pm}]\times P\otimes\ZZ[z^{\pm}]\rightarrow \ZZ$ is defined by: 
$$\epsilon_{\lambda,\mu}=\pi_0((z^{-1}\Omega^{-1}(\bar{\lambda})|\mu))$$
The multiplication of \cite{Vas} is defined by:
$$Y_{i,l}Y_{j,m}=t^{2\epsilon_{z^l\omega_i,z^m\omega_j}-2\epsilon_{z^m\omega_j,z^l\omega_i}}Y_{j,m}Y_{i,l}$$
So we can compute:
$$\epsilon_{z^l\omega_i,z^m\omega_j}=\pi_0((z^{-1}\Omega^{-1}(z^{-l}\omega_i)|z^m\omega_j))=\pi_0(\underset{k\in I}{\sum}(z^{-1-l}\tilde{C}_{i,k}(z)\alpha_k|z^m\omega_j))$$
$$=\pi_0(z^{m-l-1}\tilde{C}_{i,j}(z))=(\tilde{C}_{i,j}(z))_{l+1-m}$$
If we set $\epsilon'_{\lambda,\mu}=\pi_0((z\Omega^{-1}(\bar{\lambda})|\mu))$ then we have $\epsilon'_{z^l\omega_i,z^m\omega_j}=(\tilde{C}_{i,j}(z))_{l-1-m}$ and:
$$\epsilon_{z^l\omega_i,z^m\omega_j}-\epsilon'_{z^l\omega_i,z^m\omega_j}=\mathcal{N}({Y}_{i,l},{Y}_{j,m})$$

\section{Deformed screening operators}\label{scr} Motivated by the screening currents of \cite{Freb} we give in this section a ``concrete'' approach to deformations of screening operators. In particular the $t$-analogues of screening operators defined in \cite{Her01} will appear as commutators in $\mathcal{H}_h$. Let us begin with some background about classic screening operators.

\subsection{Reminder: classic screening operators (\cite{Fre},\cite{Fre2})}\label{screclas} 

\subsubsection{Classic screening operators and symmetry property of $q$-characters}
Recall the definition of 
\\$\pi_+(\tilde{A}_{i,l}^{\pm})=A_{i,l}^{\pm}\in \Yim$ and of $u_{i,l}:A\rightarrow\ZZ$ in section \ref{defoal}.
\begin{defi} The $i^{\text{th}}$-screening operator is the $\ZZ$-linear map defined by:\label{si}
$$S_i:\Yim\rightarrow\Yim_i=\frac{\underset{l\in\ZZ}{\bigoplus}\Yim.S_{i,l}}{\underset{l\in\ZZ}{\sum}\Yim.(S_{i,l+2r_i}-A_{i,l+r_i}.S_{i,l})}$$\label{yi}
$$\forall m\in A, S_i(m)=\underset{l\in\ZZ}{\sum}u_{i,l}(m)S_{i,l}$$
\end{defi}
\noindent Note that the $i^{\text{th}}$-screening operator can also be defined as the derivation such that:
$$S_i(1)=0\text{ , }\forall j\in I,l\in\ZZ,S_i(Y_{j,l})=\delta_{i,j}Y_{i,l}.S_{i,l}$$
\begin{thm}({\bf Frenkel, Reshetikhin, Mukhin} \cite{Fre},\cite{Fre2})\label{simme}
The image of $\chi_q:\ZZ[X_{i,l}]_{i\in I,l\in\ZZ}\rightarrow\Yim$ is:
$$\text{Im}(\chi_q)=\underset{i\in I}{\bigcap}\text{Ker}(S_i)$$
\end{thm}

\noindent It is analogous to the classical symmetry property of $\chi$: $\text{Im}(\chi)=\ZZ[y_i^{\pm}]_{i\in I}^W$.

\subsubsection{Structure of the kernel of $S_i$}\label{mi} Let $\mathfrak{K}_i=\text{Ker}(S_i)$\label{ki}. It is a $\ZZ$-subalgebra of $\Yim$.

\begin{thm}\label{deux}({\bf Frenkel, Reshetikhin, Mukhin} \cite{Fre},\cite{Fre2}) The $\ZZ$-subalgebra $\mathfrak{K}_i$ of $\Yim$ is generated by the $Y_{i,l}(1+A_{i,l+r_i}^{-1}),Y_{j,l}^{\pm}$ ($j\neq i,l\in \ZZ$). 
\end{thm}

\noindent For $m\in B_i$, we denote:
$$E_i(m)=m\underset{l\in\ZZ}{\prod}(1+{{A_{i,l+r_i}}}^{-1})^{u_{i,l}(m)}\in\mathfrak{K}_i$$\label{eim}
In particular:
\begin{cor}\label{el} The $\ZZ$-module $\mathfrak{K}_i$ is freely generated by the $E_i(m)$ ($m\in B_i$):
$$\mathfrak{K}_i=\underset{m\in B_i}{\bigoplus}\ZZ E_i(m)\simeq \ZZ^{(B_i)}$$
\end{cor}

\subsubsection{Examples in the $sl_2$-case}\label{exlsl} We suppose in this section that we are in the $sl_2$-case. For $m\in B$, let $L(m)=\chi_q(V_m)$ be the $q$-character of the $\U_q(\hat{sl_2})$-irreducible representation of highest weight $m$. In particular $L(m)\in\mathfrak{K}$ and $\mathfrak{K}=\underset{m\in B}{\bigoplus}\ZZ L(m)$.

\noindent In \cite{Fre} an explicit formula for $L(m)$ is given: a $\sigma\subset\ZZ$ is called a $2$-segment if $\sigma$ is of the form $\sigma=\{l,l+2,...,l+2k\}$. Two $2$-segment are said to be in special position if their union is a $2$-segment that properly contains each of them. All finite subset of $\ZZ$ with multiplicity $(l,u_l)_{l\in\ZZ}$ ($u_l\geq 0$) can be broken in a unique way into a union of $2$-segments which are not in pairwise special position. 

\noindent For $m\in B$ we decompose $m=\underset{j}{\prod}\underset{l\in \sigma_j}{\prod}Y_l\in B$ where the $(\sigma_j)_j$ is the decomposition of the $(l,u_l(m))_{l\in\ZZ}$. We have:
$$L(m)=\underset{j}{\prod}L(\underset{l\in \sigma_j}{\prod}Y_l)$$
So it suffices to give the formula for a $2$-segments:
$$L(Y_lY_{l+2}Y_{l+4}...Y_{l+2k})=Y_lY_{l+2}Y_{l+4}...Y_{l+2k}+Y_lY_{l+2}...Y_{l+2(k-1)}Y_{l+2(k+1)}^{-1}$$
$$+Y_lY_{l+2}...Y_{l+2(k-2)}Y_{l+2k}^{-1}Y_{l+2(k+1)}^{-1}+...+Y_{l+2}^{-1}Y_{l+2}^{-1}...Y_{l+2(k+1)}^{-1}$$
We say that $m$ is irregular if there are $j_1\neq j_2$ such that 
$$\sigma_{j_1}\subset \sigma_{j_2}\text{ and }\sigma_{j_1}+2\subset\sigma_{j_2}$$
\begin{lem}({\bf Frenkel, Reshetikhin} \cite{Fre})\label{dominl} There is a dominant $\Yim$-monomial other than $m$ in $L(m)$ if and only if $m$ is irregular.\end{lem}

\subsubsection{Complements: another basis of $\mathfrak{K}_i$}\label{compl} Let us go back to the general case. Let $\Yim_{sl_2}=\ZZ[Y_l^{\pm}]_{l\in\ZZ}$ the ring $\Yim$ for the $sl_2$-case. Let $i$ be in $I$ and for $0\leq k\leq r_i-1$, let $\omega_k: A\rightarrow \Yim_{sl_2}$ be the map defined by:
$$\omega_k(m)=\underset{l\in\ZZ}{\prod}Y_l^{u_{i,k+lr_i}(m)}$$
and $\nu_k:\ZZ[(Y_{l-1}Y_{l+1})^{-1}]_{l\in\ZZ}\rightarrow\Yim$ be the ring homomorphism such that $\nu_k((Y_{l-1}Y_{l+1})^{-1})=A_{i,k+lr_i}^{-1}$.

\noindent For $m\in B_i$, $\omega_k(m)$ is dominant in $\Yim_{sl_2}$ and so we can define $L(\omega_k(m))$ (see section \ref{exlsl}). We have $L(\omega_k(m))\omega_k(m)^{-1}\in\ZZ[(Y_{l-1}Y_{l+1})^{-1}]_{l\in\ZZ}$. We introduce:
$$L_i(m)=m\underset{0\leq k\leq r_i-1}{\prod}\nu_k(L(\omega_k(m))\omega_k(m)^{-1})\in \mathfrak{K}_i$$\label{lim}
In analogy with the corollary \ref{el} the $\ZZ$-module $\mathfrak{K}_i$ is freely generated by the $L_i(m)$ ($m\in B_i$):
$$\mathfrak{K}_i=\underset{m\in B_i}{\bigoplus}\ZZ L_i(m)\simeq \ZZ^{(B_i)}$$

\subsection{Screening currents}

Following \cite{Freb}, for $i\in I,l\in\ZZ$, introduce $\tilde{S}_{i,l}\in \mathcal{H}_h$:\label{tsil}
$$\tilde{S}_{i,l}=\text{exp}(\underset{m>0}{\sum}h^m\frac{a_{i}[m]}{q_i^{m}-q_i^{-m}}q^{lm})\text{exp}(\underset{m>0}{\sum}h^m\frac{a_i[-m]}{q_i^{-m}-q_i^{m}}q^{-lm})$$

\begin{lem}\label{currents} We have the following relations in $\mathcal{H}_h$:
$$\tilde{A}_{i,l}\tilde{S}_{i,l-r_i}=t_{-z^{-2r_i}-1}\tilde{S}_{i,l+r_i}$$
$$\tilde{S}_{i,l}\tilde{A}_{j,k}=t_{C_{i,j}(z)(z^{(k-l)}+z^{(l-k)})}\tilde{A}_{j,k}\tilde{S}_{i,l}$$
$$\tilde{S}_{i,l}\tilde{Y}_{j,k}=t_{\delta_{i,j}(z^{(k-l)}+z^{(l-k)})}\tilde{Y}_{j,k}\tilde{S}_{i,l}$$
\end{lem}

\demo

As for lemma \ref{relu} we compute in $\mathcal{H}_h$:

$\tilde{A}_{i,l}\tilde{S}_{i,l-r_i}
\\=\text{exp}(\underset{m>0}{\sum}h^ma_{i}[m]q^{lm})(\text{exp}(\underset{m>0}{\sum}h^ma_i[-m]q^{-lm})\text{exp}(\underset{m>0}{\sum}h^m\frac{a_{i}[m]}{q^{mr_i}-q^{-mr_i}}q^{m(l-r_i)})
\\\text{exp}(\underset{m>0}{\sum}h^m\frac{a_i[-m]}{q^{-mr_i}-q^{mr_i}}q^{m(r_i-l)})
\\=\text{exp}(\underset{m>0}{\sum}h^{2m}\frac{-q^{2mr_i}+q^{-2mr_i}}{q^{mr_i}-q^{-mr_i}}q^{-mr_i}c_m)\text{exp}(\underset{m>0}{\sum}h^ma_{i}[m]q^{lm}+h^m\frac{a_{i}[m]}{q^{mr_i}-q^{-mr_i}}q^{m(l-r_i)})
\\\text{exp}(\underset{m>0}{\sum}h^m\frac{a_i[-m]}{q^{-mr_i}-q^{mr_i}}q^{m(r_i-l)}+h^ma_i[-m]q^{-lm})
\\=t_{-z^{-2r_i}-1}\text{exp}(\underset{m>0}{\sum}h^m a_i[m](1+\frac{q^{-mr_i}}{q^{mr_i}-q^{-mr_i}})q^{lm})\text{exp}(\underset{m>0}{\sum}h^m a_i[m](1+\frac{q^{mr_i}}{q^{-mr_i}-q^{mr_i}})q^{-lm})
\\=t_{-z^{-2r_i}-1}\tilde{S}_{i,l+r_i}$

$\tilde{S}_{i,l}\tilde{A}_{j,k}
\\=\text{exp}(\underset{m>0}{\sum}h^m\frac{a_{i}[m]}{q^{mr_i}-q^{-mr_i}}q^{lm})(\text{exp}(\underset{m>0}{\sum}h^m\frac{a_i[-m]}{q^{-mr_i}-q^{mr_i}}q^{-lm})\text{exp}(\underset{m>0}{\sum}h^ma_{j}[m]q^{km}))
\\\text{exp}(\underset{m>0}{\sum}h^ma_j[-m]q^{-km})
\\=\text{exp}(\underset{m>0}{\sum}h^{2m}\frac{(q^{-mr_i}-q^{mr_i})C_{i,j}(q^m)}{q^{-mr_i}-q^{mr_i}}q^{m(k-l)}c_m)\text{exp}(\underset{m>0}{\sum}h^m\frac{a_{i}[m]}{q^{mr_i}-q^{-mr_i}}q^{lm})
\\\text{exp}(\underset{m>0}{\sum}h^ma_{j}[m]q^{km})\text{exp}(\underset{m>0}{\sum}h^m\frac{a_i[-m]}{q^{-mr_i}-q^{mr_i}}q^{-lm})\text{exp}(\underset{m>0}{\sum}h^ma_j[-m]q^{-km})
\\=\text{exp}(\underset{m>0}{\sum}h^{2m}C_{i,j}(q^m)(q^{m(k-l)}+q^{(l-k)m})c_m)\tilde{A}_{j,k}\tilde{S}_{i,l}$

\noindent Finally:

$\tilde{S}_{i,l}\tilde{Y}_{j,k}
\\=\text{exp}(\underset{m>0}{\sum}h^m\frac{a_{i}[m]}{q^{mr_i}-q^{-mr_i}}q^{lm})(\text{exp}(\underset{m>0}{\sum}h^m\frac{a_i[-m]}{q^{-mr_i}-q^{mr_i}}q^{-lm})\text{exp}(\underset{m>0}{\sum}h^my_{j}[m]q^{km}))
\\\text{exp}(\underset{m>0}{\sum}h^my_j[-m]q^{-km})
\\=\text{exp}(\underset{m>0}{\sum}h^{2m}\delta_{i,j}\frac{(q^{-mr_i}-q^{mr_i})}{q^{-mr_i}-q^{mr_i}}q^{m(k-l)}c_m)\text{exp}(\underset{m>0}{\sum}h^m\frac{a_{i}[m]}{q^{mr_i}-q^{-mr_i}}q^{lm})
\\\text{exp}(\underset{m>0}{\sum}h^my_{j}[m]q^{km})\text{exp}(\underset{m>0}{\sum}h^m\frac{a_i[-m]}{q^{-mr_i}-q^{mr_i}}q^{-lm})\text{exp}(\underset{m>0}{\sum}h^my_j[-m]q^{-km})
\\=\text{exp}(\underset{m>0}{\sum}\delta_{i,j}h^{2m}(q^{m(k-l)}+q^{(l-k)m})c_m)\tilde{Y}_{j,k}\tilde{S}_{i,l}$
\qed

\subsection{Deformed bimodules} In this section we define and study a $t$-analogue $\Yim_{i,t}$ of the module $\Yim_i$.

\noindent For $i\in I$, let $\Yim_{i,u}$\label{yiu} be the $\Yim_u$ sub left-module of $\mathcal{H}_h$ generated by the $\tilde{S}_{i,l}$ ($l\in \ZZ$). It follows from lemma \ref{currents} that $(\tilde{S}_{i,l})_{-r_i\leq l <r_i}$ generate $\Yim_{i,u}$ and that it is also a subbimodule of $\mathcal{H}_h$. Denote by $\tilde{S}_{i,l}\in \mathcal{H}_t$\label{tsit} the image of $\tilde{S}_{i,l}\in \mathcal{H}_h$ in $\mathcal{H}_t$.

\begin{defi} $\Yim_{i,t}$\label{yit} is the sub left-module of $\mathcal{H}_t$ generated by the $\tilde{S}_{i,l}$ ($l\in\ZZ$).\end{defi}
\noindent In particular it is to say the image of $\Yim_{i,u}$ in $\mathcal{H}_t$. It follows from lemma \ref{currents} that for $l\in\ZZ$, we have in $\Yim_{i,t}$:
$$\tilde{A}_{i,l}\tilde{S}_{i,l-r_i}=t^{-1}\tilde{S}_{i,l+r_i}$$ 
It particular $\Yim_{i,t}$ is generated by the $(\tilde{S}_{i,l})_{-r_i\leq l <r_i}$.

\noindent It follows from lemma \ref{currents} that for $l\in\ZZ$, we have:
$$\tilde{S}_{i,l}.\tilde{Y}_{j,k}=t^{2\delta_{i,j}\delta_{l,k}}\tilde{Y}_{j,k}.\tilde{S}_{i,l}\text{ , }\tilde{S}_{i,l}.t=t.\tilde{S}_{i,l}$$
In particular $\Yim_{i,t}$ a subbimodule of $\mathcal{H}_t$. Moreover:
$$\tilde{S}_{i,l}.\tilde{A}_{i,k}=t^{2\delta_{l-k,r_i}+2\delta_{l-k,-r_i}}\tilde{A}_{i,k}.\tilde{S}_{i,l}$$
$$\tilde{S}_{i,l}.\tilde{A}_{j,k}=t^{-2\underset{r=C_{i,j}+1,C_{i,j}+3,...,-C_{i,j}-1}{\sum}\delta_{l-k,r}}\tilde{A}_{j,k}.\tilde{S}_{i,l}\text{ (if $i\neq j$)}$$

\begin{prop}\label{tcur} The $\Yim_t$ left module $\Yim_{i,t}$ is freely generated by $(\tilde{S}_{i,l})_{-r_i\leq l <r_i}$:
$$\Yim_{i,t}=\underset{-r_i\leq l <r_i}{\bigoplus}\Yim_t\tilde{S}_{i,l}\simeq \Yim_t^{2r_i}$$
\end{prop}

\demo We saw that $(\tilde{S}_{i,l})_{-r_i\leq l <r_i}$ generate $\Yim_{i,t}$. We prove they are $\Yim_t$-linearly independent: 
\\for $(R_1,...,R_n)\in\mathfrak{U}^n$, introduce:
$$Y_{R_1,...,R_n}=\text{exp}(\underset{m>0, j\in I}{\sum}h^my_j[m]R_j(q^m))\in \mathcal{H}_t^+$$
It follows from lemma \ref{indgene} that the $(Y_{R})_{R\in\mathfrak{U}^n}$ are $\ZZ$-linearly independent. Note that we have $\pi_+(\Yim_{i,t})\subset \underset{R\in \mathfrak{U}^n}{\bigoplus}\ZZ Y_R$ and that $\Yim=\underset{R\in\ZZ[z^{\pm}]^n}{\bigoplus}\ZZ Y_R$. Suppose we have a linear combination ($\lambda_r\in\Yim_t$):
$$\lambda_{-r_i}\tilde{S}_{i,-r_i}+...+\lambda_{r_i-1}\tilde{S}_{i,r_i-1}=0$$ 
Introduce $\mu_{k,R}\in\ZZ$ such that:
$$\pi_+(\lambda_k)=\underset{R\in\ZZ[z^{\pm}]^n}{\sum}\mu_{k,R}Y_R$$
and $R_{i,k}=(R_{i,k}^{1}(z),...,R_{i,k}^{n}(z))\in\mathfrak{U}^n$ such that $\pi_+(\tilde{S}_{i,k})=Y_{R_{i,k}}$. If we apply $\pi_+$ to the linear combination, we get:
$$\underset{R\in\ZZ[z^{\pm}]^n,-r_i\leq k\leq r_i-1}{\sum}\mu_{k,R}Y_RY_{R_{i,k}}=0$$ 
and we have for all $R'\in \mathfrak{U}$:
$$\underset{-r_i\leq k\leq r_i-1/R'-R_{i,k}\in\ZZ[z^{\pm}]^n}{\sum}\mu_{k,R'-R_{i,k}}=0$$
Suppose we have $-r_i\leq k_1\neq k_2\leq r_i -1$ such that $R'-R_{i,k_1},R'-R_{i,k_2}\in \ZZ[z^{\pm}]^n$. So $R_{i,k_1}-R_{i,k_2}\in\ZZ[z^{\pm}]^n$. But $a_i[m]=\underset{j\in I}{\sum}C_{j,i}(q^m)y_j[m]$, so for $j\in I$:
$$C_{j,i}(z)\frac{z^{k_1}-z^{k_2}}{z^{r_i}-z^{-r_i}}=(R_{i,k_1}^j(z)-R_{i,k_2}^j(z))\in \ZZ[z^{\pm}]$$ 
In particular for $j=i$ we have $C_{i,i}(z)\frac{z^{k_1}-z^{k_2}}{z^{r_i}-z^{-r_i}}=\frac{(z^{r_i}+z^{-r_i})(z^{k_1}-z^{k_2})}{z^{r_i}-z^{-r_i}}\in \ZZ[z^{\pm}]$. This is impossible because $|k_1-k_2|<2r_i$. So we have only one term in the sum and all $\mu_{k,R}=0$. So $\pi_+(\lambda_k)=0$, and $\lambda_k\in (t-1)\Yim_t$. We have by induction for all $m>0$, $\lambda_k\in (t-1)^m\Yim_t$. It is possible if and only if $\lambda_k=0$.\qed

\noindent Denote by $\Yim_i$ the $\Yim$-bimodule $\pi_+(\Yim_{i,t})$. It is consistent with the notations of section \ref{screclas}.

\subsection{$t$-analogues of screening operators}

We introduced $t$-analogues of screening operators in \cite{Her01}. The picture of the last section enables us to define them from a new point of view.

\noindent For $m$ a $\Yim_t$-monomial, we have:
$$[\tilde{S}_{i,l},m]=\tilde{S}_{i,l}m-m\tilde{S}_{i,l}=(t^{2u_{i,l}(m)}-1)m\tilde{S}_{i,l}=t^{u_{i,l}(m)}(t-t^{-1})[u_{i,l}(m)]_tm\tilde{S}_{i,l}$$
So for $\lambda\in \Yim_t$ we have $[\tilde{S}_{i,l},\lambda]\in (t^2-1)\Yim_{i,t}$, and $[\tilde{S}_{i,l},\lambda]\neq 0$ only for a finite number of $l\in\ZZ$. So we can define:

\begin{defi} The $i^{th}$ $t$-screening operator is the map $S_{i,t}:\Yim_t\rightarrow \Yim_{i,t}$ such that ($\lambda\in\Yim_t$):
$$S_{i,t}(\lambda)=\frac{1}{t^2-1}\underset{l\in\ZZ}{\sum}[\tilde{S}_{i,l},\lambda]\in \Yim_{i,t}$$\end{defi}

\noindent In particular, $S_{i,t}$ is $\ZZ[t^{\pm}]$-linear and a derivation. It is our map of \cite{Her01}.

\noindent For $m$ a $\Yim_t$-monomial, we have $\pi_+(S_{i,t}(m))=\pi_+(t^{u_{i,l}(m)-1}[u_{i,l}(m)]_t)\pi_+(m\tilde{S}_{i,l})=u_{i,l}(m)\pi_+(m\tilde{S}_{i,l})$ and the following commutative diagram:
$$\begin{array}{rcccl}
\Yim_t&\stackrel{S_{i,t}}{\longrightarrow}&\Yim_{i,t}\\
\pi_+\downarrow &&\downarrow&\pi_+\\
\Yim&\stackrel{S_i}{\longrightarrow}&\Yim_i\end{array}$$  

\subsection{Kernel of deformed screening operators}\label{kernelun}

\subsubsection{Structure of the kernel}\label{defitl} We proved in \cite{Her01} a $t$-analogue of theorem \ref{deux}:

 \begin{thm}\label{her} (\cite{Her01}) The kernel of the $i^{th}$ $t$-screening operator $S_{i,t}$ is the $\ZZ[t^{\pm}]$-subalgebra of $\Yim_t$ generated by the $\tilde{Y}_{i,l}(1+t\tilde{A}_{i,l+r_i}^{-1}),\tilde{Y}_{j,l}^{\pm}$ ($j\neq i,l\in \ZZ$).\end{thm}

\demo For the first inclusion we compute:
$$S_{i,t}(\tilde{Y}_{i,l}(1+t\tilde{A}_{i,l+r_i}^{-1}))=\tilde{Y}_{i,l}\tilde{S}_{i,l}+t\tilde{Y}_{i,l}\tilde{A}_{i,l+r_i}^{-1}(-t^{-2})\tilde{S}_{i,l+2r_i}=\tilde{Y}_{i,l}(\tilde{S}_{i,l}-t^{-1}\tilde{A}_{i,l+r_i}^{-1}\tilde{S}_{i,l+2r_i})=0$$
For the other inclusion we refer to \cite{Her01}.
\qed

\noindent Let $\mathfrak{K}_{i,t}=\text{Ker}(S_{i,t})$\label{kit}. It is a $\ZZ[t^{\pm}]$-subalgebra of $\Yim_t$. In particular we have $\pi_+(\mathfrak{K}_{i,t})=\mathfrak{K}_i$ (consequence of theorem \ref{deux} and \ref{her}).

\noindent For $m\in B_i$ introduce: (recall that $\overset{\rightarrow}{\underset{l\in\ZZ}{\prod}}U_l$ means $...U_{-1}U_0U_1U_2...$):
$$E_{i,t}(m)=\overset{\rightarrow}{\underset{l\in\ZZ}{\prod}}((\tilde{Y}_{i,l}(1+t\tilde{A}_{i,l+r_i}^{-1}))^{u_{i,l}(m)}\underset{j\neq i}{\prod}\tilde{Y}_{j,l}^{u_{j,l}(m)})$$\label{teitm}
It is well defined because it follows from theorem \ref{dessus} that for $j\neq i, l\in\ZZ$, $(\tilde{Y}_{i,l}(1+t\tilde{A}_{i,l+r_i}^{-1}))$ and $\tilde{Y}_{j,l}$ commute. For $m\in B_i$, the formula shows that the $\Yim_t$-monomials of $E_{i,t}(m)$ are the $\Yim$-monomials of $E_i(m)$ (with identification by $\pi_+$). Such elements were used in \cite{Nab} for the $ADE$ case.

\noindent The theorem \ref{her} allows us to describe $\mathfrak{K}_{i,t}$:
\begin{cor}\label{hers} For all $m\in B_i$, we have $E_{i,t}(m)\in\mathfrak{K}_{i,t}$. Moreover: 
$$\mathfrak{K}_{i,t}=\underset{m\in B_i}{\bigoplus}\ZZ[t^{\pm}]E_{i,t}(m)\simeq \ZZ[t^{\pm}]^{(B_i)}$$\end{cor}

\demo First $E_{i,t}(m)\in\mathfrak{K}_{i,t}$ as product of elements of $\mathfrak{K}_{i,t}$. We show easily that the $E_{i,t}(m)$ are $\ZZ[t^{\pm}]$-linearly independent by looking at a maximal $\Yim_t$-monomial in a linear combination.

\noindent Let us prove that the $E_{i,t}(m)$ are $\ZZ[t^{\pm}]$-generators of $\mathfrak{K}_{i,t}$: for a product $\chi$ of the algebra-generators of theorem \ref{her}, let us look at the highest weight $\Yim_t$-monomial $m$. Then $E_{i,t}(m)$ is this product up to the order in the multiplication. But for $p=1$ or $p\geq 3$, $Y_{i,l}Y_{i,l+pr_i}$ is the unique dominant $\Yim$-monomial of $E_i(Y_{i,l})E_i(Y_{i,l+pr_i})$, so:
$$\tilde{Y}_{i,l}(1+t\tilde{A}_{i,l+r_i}^{-1})\tilde{Y}_{i,l+pr_i}(1+t\tilde{A}_{i,l+pr_i+r_i}^{-1})\in t^{\ZZ}\tilde{Y}_{i,l+pr_i}(1+t\tilde{A}_{i,l+pr_i+r_i}^{-1})\tilde{Y}_{i,l}(1+t\tilde{A}_{i,l+r_i}^{-1})$$
And for $p=2$:
$$\tilde{Y}_{i,l}(1+t\tilde{A}_{i,l+r_i}^{-1})\tilde{Y}_{i,l+2r_i}(1+t\tilde{A}_{i,l+3r_i}^{-1})-\tilde{Y}_{i,l+2r_i}(1+t\tilde{A}_{i,l+3r_i}^{-1})\tilde{Y}_{i,l}(1+t\tilde{A}_{i,l+r_i}^{-1})$$
$$\in \ZZ[t^{\pm}]+t^{\ZZ}\tilde{Y}_{i,l}(1+t\tilde{A}_{i,l+r_i}^{-1})\tilde{Y}_{i,l+2r_i}(1+t\tilde{A}_{i,l+3r_i}^{-1})$$\qed

\subsubsection{Elements of $\mathfrak{K}_{i,t}$ with a unique $i$-dominant $\Yim_t$-monomial}

\begin{prop}\label{defifprem} For $m\in B_i$, there is a unique $F_{i,t}(m)\in\mathfrak{K}_{i,t}$\label{tfitm} such that $m$ is the unique $i$-dominant $\Yim_t$-monomial of $F_{i,t}(m)$. Moreover :
$$\mathfrak{K}_{i,t}=\underset{m\in B_i}{\bigoplus}F_{i,t}(m)$$
\end{prop}

\demo It follows from corollary \ref{hers} that an element of $\mathfrak{K}_{i,t}$ has at least one $i$-dominant $\Yim_t$-monomial. In particular we have the uniqueness of $F_{i,t}(m)$. 

\noindent For the existence, let us look at the $sl_2$-case. Let $m$ be in $B$. It follows from the lemma \ref{fini} that $\{MA_{i_1,l_1}^{-1}...A_{i_R,l_R}^{-1}/R\geq 0,l_1,...,l_R\geq l(M)\}\cap B$ is finite (where $l(M)=\text{min}\{l\in\ZZ/\exists i\in I,u_{i,l}(M)\neq 0\}$). We define on this set a total ordering compatible with the partial ordering: $m_L=m>m_{L-1}>...>m_1$. Let us prove by induction on $l$ the existence of $F_t(m_l)$. The unique dominant $\Yim_t$-monomial of $E_t(m_1)$ is $m_1$ so $F_t(m_1)=E_t(m_1)$. In general let $\lambda_1(t),...,\lambda_{l-1}(t)\in\ZZ[t^{\pm}]$ be the coefficient of the dominant $\Yim_t$-monomials $m_1,...,m_{l-1}$ in $E_t(m_l)$. We put:
$$F_t(m_l)=E_t(m_l)-\underset{r=1...l-1}{\sum}\lambda_r(t)F_t(m_r)$$
Notice that this construction gives $F_t(m)\in m\ZZ[\tilde{A}_l^{-1},t^{\pm}]_{l\in\ZZ}$.

\noindent For the general case, let $i$ be in $I$ and $m$ be in $B_i$. Consider $\omega_k(m)$ as in section \ref{compl}. The study of the $sl_2$-case allows us to set $\chi_k=\omega_k(m)^{-1}F_t(\omega_k(m))\in\ZZ[\tilde{A}_l^{-1},t^{\pm}]_l$. And using the $\ZZ[t^{\pm}]$-algebra homomorphism $\nu_{k,t}:\ZZ[\tilde{A}_l^{-1},t^{\pm}]_{l\in\ZZ}\rightarrow \ZZ[\tilde{A}_{i,l}^{-1},t^{\pm}]_{i\in I,l\in\ZZ}$ defined by $\nu_{k,t}(\tilde{A}_l^{-1})=\tilde{A}_{i,k+lr_i}^{-1}$, we set (the terms of the product commute): 
$$F_{i,t}(m)=m\underset{0\leq k\leq r_i-1}{\prod}\nu_{k,t}(\chi_k)\in\mathfrak{K}_{i,t}$$ 
For the last assertion, we have $E_{i,t}(m)=\underset{l=1...L}{\sum}\lambda_l(t)F_{i,t}(m_l)$ where $m_1,...,m_L$ are the $i$-dominant $\Yim_t$-monomials of $E_{i,t}(m)$ with coefficients $\lambda_1(t),...,\lambda_L(t)\in\ZZ[t^{\pm}]$.\qed

\noindent In the same way there is a unique $F_i(m)\in\mathfrak{K}_i$\label{fim} such that $m$ is the unique $i$-dominant $\Yim$-monomial of $F_i(m)$. Moreover $F_i(m)=\pi_+(F_{i,t}(m))$.

\subsubsection{Examples in the $sl_2$-case} In this section we suppose that $\Glie=sl_2$ and we compute $F_t(m)=F_{1,t}(m)$ in some examples with the help of section \ref{exlsl}. 

\begin{lem}\label{fexp} Let $\sigma=\{l,l+2,...,l+2k\}$ be a $2$-segment and $m_{\sigma}=\tilde{Y}_l\tilde{Y}_{l+2}...\tilde{Y}_{l+2k}\in B$. Then we have the formula:
$$F_t(m_{\sigma})=m_{\sigma}(1+t\tilde{A}_{l+2k+1}^{-1}+t^2\tilde{A}_{l+(2k+1)}^{-1}\tilde{A}_{l+(2k-1)}^{-1}+...+t^k\tilde{A}_{l+(2k+1)}^{-1}\tilde{A}_{l+(2k-1)}^{-1}...\tilde{A}_{l+1}^{-1})$$
If $\sigma_1,\sigma_2$ are $2$-segments not in special position, we have:
$$F_t(m_{\sigma_1})F_t(m_{\sigma_2})=t^{\mathcal{N}(m_{\sigma_1},m_{\sigma_2})-\mathcal{N}(m_{\sigma_2},m_{\sigma_1})}F_t(m_{\sigma_2})F_t(m_{\sigma_1})$$
If $\sigma_1,...,\sigma_R$ are $2$-segments such that $m_{\sigma_1}...m_{\sigma_r}$ is regular, we have:
$$F_t(m_{\sigma_1}...m_{\sigma_R})=F_t(m_{\sigma_1})...F_t(m_{\sigma_R})$$
\end{lem}

\noindent In particular if $m\in B$ verifies $\forall l\in\ZZ, u_l(m)\leq 1$ then it is of the form $m=m_{\sigma_1}...m_{\sigma_R}$ where the $\sigma_r$ are $2$-segments such that $\text{max}(\sigma_r)+2<\text{min}(\sigma_{r+1})$. So the lemma \ref{fexp} gives an explicit formula $F_t(m)=F_t(m_{\sigma_1})...F_t(m_{\sigma_R})$.

\demo First we need some relations in $\Yim_{1,t}$ : we know that for $l\in\ZZ$ we have $t\tilde{S}_{l-1}=\tilde{A}_l^{-1}\tilde{S}_{l+1}=t^2\tilde{S}_{l+1}\tilde{A}_l^{-1}$, so $t^{-1}\tilde{S}_{l-1}=\tilde{S}_{l+1}\tilde{A}_l^{-1}$. So we get by induction that for $r\geq 0$:
$$t^{-r}\tilde{S}_{l+1-2r}=\tilde{S}_{l+1}\tilde{A}_l^{-1}\tilde{A}_{l-2}^{-1}...\tilde{A}_{l-2(r-1)}^{-1}$$
As $u_{i,l+1}(\tilde{A}_l^{-1}\tilde{A}_{l-2}^{-1}...\tilde{A}_{l-2(r-1)}^{-1})=u_{i,l+1}(\tilde{A}_l^{-1})=-1$, we get:
$$t^{-r}\tilde{S}_{l+1-2r}=t^{-2}\tilde{A}_l^{-1}\tilde{A}_{l-2}^{-1}...\tilde{A}_{l-2(r-1)}^{-1}\tilde{S}_{l+1}$$
For $r'\geq 0$, by multiplying on the left by $\tilde{A}_{l+2r'}^{-1}\tilde{A}_{l+2(r'-1)}^{-1}...\tilde{A}_{l+2}^{-1}$, we get:
$$t^{-r}\tilde{A}_{l+2r'}^{-1}\tilde{A}_{l+2(r'-1)}^{-1}...\tilde{A}_{l+2}^{-1}\tilde{S}_{l+1-2r}=t^{-2}\tilde{A}_{l+2r'}^{-1}\tilde{A}_{l+2(r'-1)}^{-1}...\tilde{A}_{l-2(r-1)}^{-1}\tilde{S}_{l+1}$$
If we put $r'=1+R',r=R-R',l=L-1-2R'$, we get for $0\leq R'\leq R$:
$$t^{R'}\tilde{A}_{L+1}^{-1}\tilde{A}_{L-1}^{-1}...\tilde{A}_{L+1-2R'}^{-1}\tilde{S}_{L-2R}=t^{R-2}\tilde{A}_{L+1}^{-1}\tilde{A}_{L-1}^{-1}...\tilde{A}_{L+1-2R}^{-1}\tilde{S}_{L-2R'}$$
Now let be $m=\tilde{Y}_0\tilde{Y}_2...\tilde{Y}_l$ and $\chi\in\Yim_t$ given by the formula in the lemma. Let us compute $\tilde{S}_t(\chi)$:
$$\tilde{S}_t(\chi)=m(\tilde{S}_0+\tilde{S}_2+...+\tilde{S}_l)$$
$$+tm\tilde{A}^{-1}_{l+1}(\tilde{S}_0+\tilde{S}_2+...+\tilde{S}_{l-2}-t^{-2}\tilde{S}_{l+2})$$
$$+t^2m\tilde{A}^{-1}_{l+1}\tilde{A}^{-1}_{l-1}(\tilde{S}_0+\tilde{S}_2+...+\tilde{S}_{l-4}-t^{-2}\tilde{S}_{l}-t^{-2}\tilde{S}_{l+2})$$
$$+...$$
$$+t^lm\tilde{A}^{-1}_{l+1}\tilde{A}^{-1}_{l-1}...\tilde{A}^{-1}_{1}(-t^{-2}\tilde{S}_2+...-t^{-2}\tilde{S}_{l}-t^{-2}\tilde{S}_{l+2})$$
$$=m(\tilde{S}_0+\tilde{S}_2+...+\tilde{S}_l)$$
$$+tm\tilde{A}^{-1}_{l+1}(\tilde{S}_0+\tilde{S}_2+...+\tilde{S}_{l-2})-m\tilde{S}_l$$
$$+t^2m\tilde{A}^{-1}_{l+1}\tilde{A}^{-1}_{l-1}(\tilde{S}_0+\tilde{S}_2+...+\tilde{S}_{l-4})-tm\tilde{A}^{-1}_{l+1}\tilde{S}_{l-2}-m\tilde{S}_{l-2}$$
$$+...$$
$$-mt^{l-1}\tilde{A}^{-1}_{l+1}\tilde{A}^{-1}_{l-1}...\tilde{A}^{-1}_{3}-...-t^{-2}\tilde{S}_{l}-m\tilde{S}_0$$
$$=0$$
So $\chi\in\mathfrak{K}_t$. But we see on the formula that $m$ is the unique dominant monomial of $\chi$. So $\chi=F_t(m)$.

\noindent For the second point, we have two cases:

if $m_{\sigma_1}m_{\sigma_2}$ is regular, it follows from lemma \ref{dominl} that $L(m_{\sigma_1})L(m_{\sigma_2})=L(m_{\sigma_2})L(m_{\sigma_1})$ has no dominant monomial other than $m_{\sigma_1}m_{\sigma_2}$. But our formula shows that $F_t(m_{\sigma_1})$ (resp. $F_t(m_{\sigma_2}$)) has the same monomials than $L(m_{\sigma_1})$ (resp. $L(m_{\sigma_2})$). So 
$$F_t(m_{\sigma_1})F_t(m_{\sigma_2})-t^{\mathcal{N}(m_{\sigma_1},m_{\sigma_2})-\mathcal{N}(m_{\sigma_2},m_{\sigma_1})}F_t(m_{\sigma_2})F_t(m_{\sigma_1})$$ 
has no dominant $\Yim_t$-monomial because $m_{\sigma_1}m_{\sigma_2}-t^{\mathcal{N}(m_{\sigma_1},m_{\sigma_2})-\mathcal{N}(m_{\sigma_2},m_{\sigma_1})}m_{\sigma_2}m_{\sigma_1}=0$.

if $m_{\sigma_1}m_{\sigma_2}$ is irregular, we have for example $\sigma_{j_1}\subset \sigma_{j_2}$ and $\sigma_{j_1}+2\subset\sigma_{j_2}$. Let us write $\sigma_{j_1}=\{l_1,l_1+2,...,,p_1\}$ and $\sigma_2=\{l_2,l_2+2,...,,p_2\}$. So we have $l_2\leq l_1$ and $p_1\leq p_2-2$. Let $m=m_1m_2$ be a dominant $\Yim$-monomial of $L(m_{\sigma_1}m_{\sigma_2})=L(m_{\sigma_1})L(m_{\sigma_2})$ where $m_1$ (resp. $m_2$) is a $\Yim$-monomial of $L(m_{\sigma_1})$ (resp. $L(m_{\sigma_2})$). If $m_2$ is not $m_{\sigma_2}$, we have $Y_{p_2}^{-1}$ in $m_2$ which can not be canceled by $m_1$. So $m=m_1m_{\sigma_2}$. Let us write $m_1=m_{\sigma_1}A_{p_1+1}^{-1}...A_{p_1+1-2r}^{-1}$. So we just have to prove:
$$\tilde{A}_{p_1+1}^{-1}...\tilde{A}_{p_1+1-2r}^{-1}m_{\sigma_2}=m_{\sigma_2}\tilde{A}_{p_1+1}^{-1}...\tilde{A}_{p_1+1-2r}^{-1}$$
This follows from ($l\in\ZZ$):
$$\tilde{A}_l^{-1}\tilde{Y}_{l-1}\tilde{Y}_{l+1}=\tilde{Y}_{l-1}\tilde{Y}_{l+1}\tilde{A}_l^{-1}$$

\noindent For the last assertion it suffices to show that $F_t(m_{\sigma_1})...F_t(m_{\sigma_R})$ has no other dominant $\Yim_t$-monomial than $m_{\sigma_1}...m_{\sigma_R}$. But $F_t(m_{\sigma_1})...F_t(m_{\sigma_R})$ has the same monomials than $L(m_{\sigma_1})...L(m_{\sigma_R})=L(m_{\sigma_1}...m_{\sigma_R})$. As $m_{\sigma_1}...m_{\sigma_R}$ is regular we get the result.\qed

\subsubsection{Technical complements} Let us go back to the general case. We give some technical results which will be used in the following to compute $F_{i,t}(m)$ in some cases (see proposition \ref{cpfacile} and section \ref{fin}).

\begin{lem}\label{calcn} Let $i$ be in $I$, $l\in\ZZ$, $M\in A$ such that $u_{i,l}(M)=1$ and $u_{i,l+2r_i}=0$. Then we have $\mathcal{N}(M,\tilde{A}_{i,l+r_i}^{-1})=-1$. In particular $\pi^{-1}(MA_{i,l+r_i}^{-1})=tM\tilde{A}_{i,l+r_i}^{-1}$.\end{lem}

\demo We can suppose $M=:M:$ and we compute in $\Yim_u$:
$$M\tilde{A}_{i,l+r_i}^{-1}=\pi_+(m)\text{exp}(\underset{m>0,r\in\ZZ,j\in I}{\sum}u_{j,r}(M)h^mq^{-rm}y_j[-m])$$
$$\text{exp}(\underset{m>0}{\sum}-h^mq^{-(l+r_i)m}a_i[-m])\text{exp}(\underset{m>0}{\sum}-h^mq^{(l+r_i)m}a_i[m])$$
$$=:M\tilde{A}_{i,l+r_i}^{-1}:\text{exp}(\underset{m>0}{\sum}h^{2m}([a_i[-m],a_i[m]]-\underset{r\in\ZZ}{\sum}u_{i,r}(m)[y_i[-m],a_i[m]]q^{(l+r_i-r)m}c_m)=t_R:\tilde{Y}_{i,l}\tilde{A}_{i,l+r_i}^{-1}:$$
where:
$$R(z)=-(z^{2r_i}-z^{-2r_i})+\underset{r\in\ZZ}{\sum}u_{i,r}(M)z^{(l+r_i-r)}(z^{r_i}-z^{-r_i})$$
So:
$$\mathcal{N}(\tilde{Y}_{i,l},M\tilde{A}_{i,l+r_i}^{-1})=\underset{r\in\ZZ}{\sum}u_{i,r}(M)(z^{2r_i+l-r}-z^{l-r})_0=-u_{i,l}(M)+u_{i,l+2r_i}(M)=-1$$
\qed

\begin{lem}\label{pidonne} Let $m$ be in $B_i$ such that $\forall l\in\ZZ, u_{i,l}(m)\leq 1$ and for $1\leq r\leq 2r_i$ the set 
\\$\{l\in\ZZ/u_{i,r+2lr_i}(m)=1\}$ is a $1$-segment. Then we have $F_{i,t}(m)=\pi^{-1}(F_i(m))$.\end{lem}

\demo Let us look at the $sl_2$-case : $m=m_1m_2=m_{\sigma_1}m_{\sigma_2}$ where $\sigma_1,\sigma_2$ are $2$-segment. So the lemma \ref{fexp} gives an explicit formula for $F_t(m)$ and it follows from lemma \ref{calcn} that $F_t(m)=\pi^{-1}(F(m))$.

\noindent We go back to the general case : let us write $m=m'm_1...m_{2r_i}$ where $m'=\underset{j\neq i,l\in\ZZ}{\prod}Y_{j,l}^{u_{j,l}(m)}$ and $m_r=\underset{l\in\ZZ}{\prod}Y_{i,r+2lr_i}^{u_{i,r+2lr_i}(m)}$. We have $m_r$ of the form $m_r=Y_{i,l_r}Y_{i,l_r+2r_i}...Y_{i,l_r+2n_ir_i}$. We have $F_{i,t}(m)=t^{-N(m'm_1...m_r)}m'F_{i,t}(m_1)...F_{i,t}(m_{2r_i})$. The study of the $sl_2$-case gives $F_{i,t}(m_r)=\pi^{-1}(F_i(m_r))$. It follows from lemma \ref{calcn} that:
$$t^{-N(m'm_1...m_r)}m'\pi^{-1}(F_i(m_1))...\pi^{-1}(F_i(m_r))=\pi^{-1}(m'F_i(m_1)...F_i(m_r))=\pi^{-1}(F_i(m))$$
\qed

\section{Intersection of kernels of deformed screening operators}\label{rests} Motivated by theorem \ref{simme} we study the structure of a completion of $\mathfrak{K}_t=\underset{i\in I}{\bigcap}\text{Ker}(S_{i,t})$ in order to construct $\chi_{q,t}$ in section \ref{conschi}. Note that in the $sl_2$-case we have $\mathfrak{K}_t=\text{Ker}(S_{1,t})$ that was studied in section \ref{scr}.

\subsection{Reminder: classic case (\cite{Fre}, \cite{Fre2})}
\subsubsection{The elements $E(m)$ and $q$-characters}

For $J\subset I$, denote the $\ZZ$-subalgebra $\mathfrak{K}_J=\underset{i\in J}{\bigcap}\mathfrak{K}_i\subset\Yim$ and $\mathfrak{K}=\mathfrak{K}_I$.

\begin{lem}\label{least} (\cite{Fre}, \cite{Fre2}) A non zero element of $\mathfrak{K}_J$ has at least one $J$-dominant $\Yim$-monomial.\end{lem}

\demo It suffices to look at a maximal weight $\Yim$-monomial $m$ of $\chi\in\mathfrak{K}_J$: for $i\in J$ we have $m\in B_i$ because $\chi\in\mathfrak{K}_i$.\qed

\begin{thm}(\cite{Fre}, \cite{Fre2}) For $i\in I$ there is a unique $E(Y_{i,0})\in \mathfrak{K}$\label{ei} such that ${Y}_{i,0}$ is the unique dominant $\Yim$-monomial in $E(Y_{i,0})$.
\end{thm}

\noindent The uniqueness follows from lemma \ref{least}. For the existence we have $E(Y_{i,0})=\chi_q(V_{\omega_i}(1))$ (theorem \ref{simme}).

\noindent Note that the existence of $E(Y_{i,0})\in\mathfrak{K}$ suffices to characterize $\chi_q: \text{Rep}\rightarrow \mathfrak{K}$. It is the ring homomorphism such that $\chi_q(X_{i,l})=s_l(E(Y_{i,0}))$ where $s_l:\Yim\rightarrow \Yim$ is given by $s_l({Y}_{j,k})={Y}_{j,k+l}$. 

\noindent For $m\in B$, we defined the standard module $M_m$ in section \ref{back}. We set: 
$$E(m)=\underset{m\in B}{\prod}s_l(E(Y_{i,0}))^{u_{i,l}(m)}=\chi_q(M_m)\in\mathfrak{K}$$\label{em}
\noindent We defined the simple module $V_m$ in section \ref{back}. We set $L(m)=\chi_q(V_m)\in\mathfrak{K}$\label{lm}. We have:
$$\mathfrak{K}=\underset{m\in B}{\bigoplus}\ZZ E(m)=\underset{m\in B}{\bigoplus}\ZZ L(m)\simeq \ZZ^{(B)}$$
For $m\in B$, we can also define a unique $F(m)\in\mathfrak{K}$\label{fm} such that $m$ is the unique dominant $\Yim$-monomial which appears in $F(m)$ (see for example the proof of proposition \ref{defifprem}).

\subsubsection{Technical complements} For $J\subset I$, let $\Glie_J$ be the semi-simple Lie algebra of Cartan Matrix $(C_{i,j})_{i,j\in J}$ and $\U_q(\hat{\Glie})_J$ the associated quantum affine algebra with coefficient $(r_i)_{i\in J}$. In analogy with the definition of $E_i(m),L_i(m)$ using the $sl_2$-case (section \ref{compl}), we define for $m\in B_J$: $E_J(m)$, $L_J(m)$, $F_J(m)\in \mathfrak{K}_J$ using $\U_q(\hat{\Glie})_J$. We have:
$$\mathfrak{K}_J=\underset{m\in B_J}{\bigoplus}\ZZ E_J(m)=\underset{m\in B_J}{\bigoplus}\ZZ L_J(m)=\underset{m\in B_J}{\bigoplus}\ZZ F_J(m)\simeq \ZZ^{(B_J)}$$

\noindent As a direct consequence of proposition \ref{aidafm} we have :

\begin{lem}\label{aidahfm} For $m\in B$, we have $E(m)\in\ZZ[Y_{i,l}]_{i\in I,l\geq l(m)}$ where $l(m)=\text{min}\{l\in\ZZ/\exists i\in I,u_{i,l}(m)\neq 0\}$.\end{lem}

\subsection{Completion of the deformed algebras}\label{complesection} In this section we introduce completions of $\Yim_t$ and of $\mathfrak{K}_{J,t}=\underset{i\in J}{\bigcap}\mathfrak{K}_{i,t}\subset \Yim_t$ ($J\subset I$). We have the following motivation: we have seen $\pi_+(\mathfrak{K}_{J,t})\subset \mathfrak{K}_J$ (section \ref{scr}). In order to prove an analogue of the other inclusion (theorem \ref{con}) we have to introduce completions where infinite sums are allowed.

\subsubsection{The completion $\Yim_t^{\infty}$ of $\Yim_t$}

Let $\overset{\infty}{A}_t$\label{infat} be the $\ZZ[t^{\pm}]$-module $\overset{\infty}{A}_t=\underset{m\in A}{\prod}\ZZ[t^{\pm}].m\simeq \ZZ[t^{\pm}]^A$. An element $(\lambda_m(t)m)_{m\in A}\in \overset{\infty}{A}_t$ is noted $\underset{m\in A}{\sum}\lambda_m(t)m$. We have $\underset{m\in A}{\bigoplus}\ZZ[t^{\pm}].m=\Yim_t\subset \overset{\infty}{A}_t$. The algebra structure of $\Yim_t$ gives a $\ZZ[t^{\pm}]$-bilinear morphisms $\Yim_t\otimes\overset{\infty}{A}_t\rightarrow\overset{\infty}{A}_t$ and $\overset{\infty}{A}_t\otimes\Yim_t\rightarrow\overset{\infty}{A}_t$ such that $\overset{\infty}{A}_t$ is a $\Yim_t$-bimodule. But the $\ZZ[t^{\pm}]$-algebra structure of $\Yim_t$ can not be naturally extended to $\overset{\infty}{A}_t$. We define a $\ZZ[t^{\pm}]$-submodule $\Yim_t^{\infty}$\label{tytinf} with $\Yim_t\subset\Yim_t^{\infty}\subset\overset{\infty}{A}_t$, for which it is the case:

\noindent Let $\Yim_t^A$\label{tyta} be the $\ZZ[t^{\pm}]$-subalgebra of $\Yim_t$ generated by the $(\tilde{A}_{i,l}^{-1})_{i\in I,l\in\ZZ}$. We gave in proposition \ref{yenga} the structure of $\Yim_t^A$. In particular we have $\Yim_t^A=\underset{K\geq 0}{\bigoplus}\Yim_t^{A,K}$ where for $K\geq 0$:
$$\Yim_t^{A,K}=\underset{m=:\tilde{A}_{i_1,l_1}^{-1}...\tilde{A}_{i_K,l_K}^{-1}:}{\bigoplus}\ZZ[t^{\pm}].m\subset \Yim_t^A$$
Note that for $K_1,K_2\geq 0$, $\Yim_t^{A,K_1}\Yim_t^{A,K_2}\subset\Yim_t^{A,K_1+K_2}$ for the multiplication of $\Yim_t$. So $\Yim_t^A$ is a graded algebra if we set $\text{deg}(x)=K$ for $x\in\Yim_t^{A,K}$. Denote by $\Yim_t^{A,\infty}$ the completion of $\Yim_t^A$ for this gradation. It is a sub-$\ZZ[t^{\pm}]$-module of $\overset{\infty}{A}_t$. 

\begin{defi} We define $\Yim_t^{\infty}$ as the sub $\Yim_t$-leftmodule of $\overset{\infty}{A}_t$ generated by $\Yim_t^{A,\infty}$.\end{defi}
\noindent In particular, we have: $\Yim_t^{\infty}=\underset{M\in A}{\sum}M.\Yim_t^{A,\infty}\subset \overset{\infty}{A}_t$.

\begin{lem} There is a unique algebra structure on $\Yim_t^{\infty}$ compatible with the structure of $\Yim_t\subset \Yim_t^{\infty}$.\end{lem}

\demo The structure is unique because the elements of $\Yim_t^{\infty}$ are infinite sums of elements of $\Yim_t$. For $M\in A$, we have $\Yim_t^{A,\infty}.M\subset M.\Yim_t^{A,\infty}$, so $\Yim_t^{\infty}$ is a sub $\Yim_t$-bimodule of $\overset{\infty}{A}_t$. For $M\in A$ and $\lambda\in\Yim_t^{A, \infty}$ denote $\lambda^M\in\Yim_t^{A, \infty}$ such that $\lambda.M=M.\lambda^M$. We define the $\ZZ[t^{\pm}]$-algebra structure on $\Yim_t^{\infty}$ by ($M,M'\in A,\lambda,\lambda'\in\Yim_t^{A,\infty}$): 
$$(M.\lambda)(M'.\lambda')=MM'.(\lambda^{M'}\lambda')$$
It is well defined because for $M_1,M_2,M\in A,\lambda,\lambda_2\in\Yim_t^A$ we have $M_1\lambda_1=M_2\lambda_2\Rightarrow M_1M\lambda_1^M=M_2M\lambda_2^M$.\qed

\subsubsection{The completion $\mathfrak{K}_{i,t}^{\infty}$ of $\mathfrak{K}_{i,t}$} We define a completion of $\mathfrak{K}_{i,t}$ analog to the completed algebra $\Yim_t^{\infty}$. 

\noindent For $M\in A$, we define a $\ZZ[t^{\pm}]$-linear endomorphism $E_{i,t}^M:M\Yim_t^{A,\infty}\rightarrow M\Yim_t^{A,\infty}$\label{teitmap} such that ($m$ $\Yim_t^A$-monomial):
$$E_{i,t}^M(Mm)=0\text{ if $:Mm:\notin B_i$}$$
$$E_{i,t}^M(Mm)=E_{i,t}(Mm)\text{ if $:Mm:\in B_i$}$$
It is well-defined because if $m\in\Yim_t^{A,K}$ and $:Mm:\in B_i$ we have $E_{i,t}(Mm)\in M\underset{K'\geq K}{\bigoplus}\Yim^{A,K'}$.

\begin{defi}We define $\mathfrak{K}_{i,t}^{\infty}=\underset{M\in A}{\sum}\text{Im}(E_{i,t}^M)\subset\Yim_t^{\infty}$\label{kitinf}.\end{defi}
\noindent For $J\subset I$, we set $\mathfrak{K}_{J,t}^{\infty}=\underset{i\in J}{\bigcap}\mathfrak{K}_{i,t}^{\infty}$ and $\mathfrak{K}_t^{\infty}=\mathfrak{K}_{I,t}^{\infty}$.

\begin{lem}\label{leasto} A non zero element of $\mathfrak{K}_{J,t}^{\infty}$ has at least one $J$-dominant $\Yim_t$-monomial.\end{lem}

\demo Analog to the proof of lemma \ref{least}.

\begin{lem}\label{alginf} For $J\subset I$, we have $\mathfrak{K}_{J,t}^{\infty}\cap \Yim_t=\mathfrak{K}_{J,t}$. Moreover $\mathfrak{K}_{J,t}^{\infty}$ is a $\ZZ[t^{\pm}]$-subalgebra of $\Yim_t^{\infty}$.\end{lem}

\demo It suffices to prove the results for $J=\{i\}$. First for $m\in B_i$ we have $E_{i,t}(m)=E_{i,t}^m(m)\in\mathfrak{K}_{i,t}^{\infty}$ and so $\mathfrak{K}_{i,t}=\underset{m\in B_i}{\bigoplus}\ZZ[t^{\pm}]E_{i,t}(m)\subset \mathfrak{K}_{i,t}^{\infty}\cap\Yim_t$. Now let $\chi$ be in $\mathfrak{K}_{i,t}^{\infty}$ such that $\chi$ has only a finite number of $\Yim_t$-monomials. In particular it has only a finite number of $i$-dominant $\Yim_t$-monomials $m_1,...,m_r$ with coefficients $\lambda_1(t),...,\lambda_r(t)$. In particular it follows from lemma \ref{leasto} that $\chi=\lambda_1(t)F_{i,t}(m_1)+...+\lambda_r(t)F_{i,t}(m_r)\in\mathfrak{K}_{i,t}$ (see proposition \ref{defifprem} for the definition of $F_{i,t}(m)$).

\noindent For the last assertion, consider $M_1$, $M_2\in A$ and $m_1$, $m_2$ $\Yim_t^A$-monomials such that $:M_1m_1:,:M_2m_2:\in B_i$. Then $E_{i,t}(M_1m_1)E_{i,t}(M_2m_2)$ is in the the sub-algebra $\mathfrak{K}_{i,t}\subset\Yim_t$ and in $\text{Im}(E_{i,t}^{M_1M_2})$.\qed

\noindent In the same way for $t=1$ we define the $\ZZ$-algebra $\mathfrak{\Yim}^{\infty}$ and the $\ZZ$-subalgebras $\mathfrak{K}_J^{\infty}\subset\Yim^{\infty}$.

\noindent The surjective map $\pi_+:\Yim_t\rightarrow\Yim$ is naturally extended to a surjective map $\pi_+:\Yim_t^{\infty}\rightarrow\Yim^{\infty}$. For $i\in I$, we have $\pi_+(\mathfrak{K}_{i,t}^{\infty})=\mathfrak{K}_i^{\infty}$ and for $J\subset I$, $\pi_+(\mathfrak{K}_{J,t}^{\infty})\subset\mathfrak{K}_J^{\infty}$. The other inclusion is equivalent to theorem \ref{con}.

\subsubsection{Special submodules of $\Yim_t^{\infty}$}\label{infsum} For $m\in A$, $K\geq 0$ we construct a subset $D_{m,K}\subset m\{\tilde{A}_{i_1,l_1}^{-1}...\tilde{A}_{i_K,l_K}^{-1}\}$\label{dmk} stable by the maps $E_{i,t}^{m}$ such that $\underset{K\geq 0}{\bigcup}D_{m,K}$ is countable: we say that $m'\in D_{m,K}$ if and only if there is a finite sequence $(m_0=m,m_1,...,m_R=m')$ of length $R\leq K$, such that for all $1\leq r\leq R$, there is $r'<r$, $J\subset I$ such that $m_{r'}\in B_J$ and for $r'<r''\leq r$, $m_{r''}$ is a $\Yim$-monomial of $E_J(m_{r'})$ and $m_{r''}m_{r''-1}^{-1}\in\{A_{j,l}^{-1}/l\in\ZZ,j\in J\}$.

\noindent The definition means that ``there is chain of monomials of some $E_J(m'')$ from $m$ to $m'$''.

\begin{lem} The set $D_{m,K}$ is finite. In particular, the set $D_m$ is countable.\end{lem}

\demo Let us prove by induction on $K\geq 0$ that $D_{m,K}$ is finite: we have $D_{m,0}=\{m\}$ and: 
$$D_{m,K+1}\subset \underset{J\subset I, m'\in D_{m,K}\cap B_J}{\bigcup}\{\text{$\Yim$-monomials of }E_J(m')\}$$\qed

\begin{lem}\label{ordrep} For $m,m'\in A$ such that $m'\in D_{m}$ we have $D_{m'}\subseteq D_{m}$. For $M\in A$, the set $B\cap D_M$ is finite.\end{lem}

\demo Consider $(m_0=m,m_1,...,m_R=m')$ a sequence adapted to the definition of $D_m$. Let $m''$ be in $D_{m'}$ and $(m_R=m',m_{R+1},...,m_{R'}=m'')$ a sequence adapted to the definition of $D_{m'}$. So $(m_0,m_1,...,m_{R'})$ is adapted to the definition of $D_m$, and $m''\in D_m$.

\noindent Let us look at $m\in B\bigcap D_M$: we can see by induction on the length of a sequence $(m_0=M,m_1,...,m_R=m)$ adapted to the definition of $D_M$ that $m$ is of the form $m=MM'$ where $M'=\underset{i\in I,l\geq l_1}{\prod}A_{i,l}^{-v_{i,l}}$ ($v_{i,l}\geq 0$). So the last assertion follows from lemma \ref{fini}.\qed

\begin{defi} $\tilde{D}_m$\label{tdm} is the $\ZZ[t^{\pm}]$-submodule of $\Yim_t^{\infty}$ whose elements are of the form $(\lambda_m(t)m)_{m\in D_m}$.\end{defi}

\noindent For $m\in A$ introduce $m_0=m>m_1>m_2>...$ the countable set $D_m$ with a total ordering compatible with the partial ordering. For $k\geq 0$ consider an element $F_k\in \tilde{D}_{m_k}$. 

\noindent Note that some infinite sums make sense in $\tilde{D}_m$: for $k\geq 0$, we have $D_{m_k}\subset \{m_k,m_{k+1},...\}$. So $m_k$ appears only in the $F_{k'}$ with $k'\leq k$ and the infinite sum $\underset{k\geq 0}{\sum}F_k$ makes sense in $\tilde{D}_m$.

\subsection{Crucial result for our construction}

\noindent Our construction of $q,t$-characters is based on theorem \ref{con} proved in this section.

\subsubsection{Statement} 

\begin{defi}\label{tftm} For $n\geq 1$ denote $P(n)$\label{pn} the property ``for all semi-simple Lie-algebras $\Glie$ of rank $\text{rk}(\Glie)=n$, for all $m\in B$ there is a unique $F_t(m)\in\mathfrak{K}_t^{\infty}\cap \tilde{D}_m$ such that $m$ is the unique dominant $\Yim_t$-monomial of $F_t(m)$.''.\end{defi}

\begin{thm}\label{con} For all $n\geq 1$, the property $P(n)$ is true.\end{thm}

\noindent Note that for $n=1$, that is to say $\Glie=sl_2$, the result follows from section \ref{scr}.

\noindent The uniqueness follows from lemma \ref{leasto} : if $\chi_1,\chi_2\in \mathfrak{K}_t^{\infty}$ are solutions, then $\chi_1-\chi_2$ has no dominant $\Yim_t$-monomial, so $\chi_1=\chi_2$.

\noindent Remark: in the simply-laced case the existence is a consequence of the geometric theory of quivers \cite{Naa}, \cite{Nab}, and in $A_n,D_n$-cases of algebraic explicit constructions \cite{Nac}. In the rest of this section \ref{rests} we give an algebraic proof of this theorem in the general case.

\subsubsection{Outline of the proof} First we give some preliminary technical results (section \ref{conste}) in which we construct $t$-analogues of the $E(m)$. Next we prove $P(n)$ by induction on $n$. Our proof has 3 steps:

\noindent Step 1 (section \ref{petitdeux}): we prove $P(1)$ and $P(2)$ using a more precise property $Q(n)$ such that $Q(n)\Rightarrow P(n)$. The property $Q(n)$ has the following advantage: it can be verified by computation in elementary cases $n=1,2$. 

\noindent Step 2 (section \ref{consp}): we give some consequences of $P(n)$ which will be used in the proof of $P(r)$ ($r>n$): we give the structure of $\mathfrak{K}_t^{\infty}$ (proposition \ref{thth}) for $\text{rk}(\Glie)=n$ and the structure of $\mathfrak{K}_{J,t}^{\infty}$ where $J\subset I$, $|J|=n$ and $|I|>n$ (corollary \ref{recufond}). 

\noindent Step 3 (section \ref{proof}): we prove $P(n)$ ($n \geq 3$) assuming $P(r)$, $r\leq n$ are true. We give an algorithm (section \ref{defialgo}) to construct explicitly $F_t(m)$. It is called $t$-algorithm and is a $t$-analogue of Frenkel-Mukhin algorithm \cite{Fre2} (a deformed algorithm was also used by Nakajima in the $ADE$-case \cite{Naa}). As we do not know {\it a priori} the algorithm is well defined the general case, we have to show that it never fails (lemma \ref{nfail}) and gives a convenient element (lemma \ref{conv}).

\subsection{Preliminary: Construction of the $E_t(m)$}\label{conste}

\begin{lem}\label{copieun} We suppose that for $i\in I$, there is $F_t(\tilde{Y}_{i,0})\in\mathfrak{K}_t^{\infty}\cap \tilde{D}_{\tilde{Y}_{i,0}}$ such that $\tilde{Y}_{i,0}$ is the unique dominant $\Yim_t$-monomial of $F_t(\tilde{Y}_{i,0})$. Then:

i) All $\Yim_t$-monomials of $F_t(\tilde{Y}_{i,0})$, except the highest weight $\Yim_t$-monomial, are right negative.

ii) All $\Yim_t$-monomials of $F_t(\tilde{Y}_{i,0})$ are products of $\tilde{Y}_{j,l}^{\pm}$ with $l\geq 0$. 

iii) The only $\Yim_t$-monomial of $F_t(\tilde{Y}_{i,0})$ which contains a $\tilde{Y}_{j,0}^{\pm}$ ($j\in I$) is the highest weight monomial $\tilde{Y}_{i,0}$.

iv) The $F_t(\tilde{Y}_{i,0})$ ($i\in I$) commute.\end{lem}
\noindent Note that (i),(ii) and (iii) appeared in \cite{Fre2}.

\demo 

i) It suffices to prove that all $\Yim_t$-monomials $m_0=Y_{i,0},m_1,...$ of $D_{Y_{i,0}}$ except $Y_{i,0}$ are right negative. But $m_1$ is the monomial $Y_{i,0}A_{i,1}^{-1}$ of $E_i(Y_{i,0})$ and it is right negative. We can now prove the statement by induction: suppose that $m_r$ is a monomial of $E_J(m_{r'})$, where $m_{r'}$ is right negative. So $m_r$ is a product of $m_{r'}$ by some $A_{j,l}^{-1}$ ($l\in\ZZ$).Those monomials are right negative because a product of right negative monomial is right negative. 

ii) Suppose that $m\in A$ is product of $Y_{k,l}^{\pm}$ with $l\geq 0$. It follows from lemma \ref{aidahfm} that all monomials of $D_{m}$ are product of $Y_{k,l}^{\pm}$ with $l\geq 0$.

iii) All $\Yim$-monomials of $D_{Y_{i,0}}$ except $\tilde{Y}_{i,0}$ are in $D_{Y_{i,0}A_{i,r_i}^{-1}}$. But $l(Y_{i,0}A_{i,r_i}^{-1})\geq 1$ and we can conclude with the help of lemma \ref{aidahfm}.

iv) Let $i\neq j$ be in $I$ and look at $F_t(\tilde{Y}_{i,0})F_t(\tilde{Y}_{j,0})$. Suppose we have a dominant $\Yim_t$-monomial $m_0=m_1m_2$ in $F_t(\tilde{Y}_{i,0})F_t(\tilde{Y}_{j,0})$ different from the highest weight $\Yim_t$-monomial $\tilde{Y}_{i,0}\tilde{Y}_{j,0}$. We have for example $m_1\neq \tilde{Y}_{i,0}$, so $m_1$ is right negative. Let $l_1$ be the maximal $l$ such that a $\tilde{Y}_{k,l}$ appears in $m_1$. We have $u_{k,l}(m_1)<0$ and $l>0$. As $u_{k,l}(m_0)\geq 0$ we have $u_{k,l}(m_2)>0$ and $m_2\neq {Y}_{j,0}$. So $m_2$ is right negative and there is $k'\in I$ and $l'>l$ such that $u_{k',l'}(m_2)<0$. So $u_{k',l'}(m_1)>0$, contradiction. So the highest weight $\Yim_t$-monomial of $F_t(\tilde{Y}_{i,0})F_t(\tilde{Y}_{j,0})$ is the unique dominant $\Yim_t$-monomial. In the same way the highest weight $\Yim_t$-monomial of $F_t(\tilde{Y}_{j,0})F_t(\tilde{Y}_{i,0})$ is the unique dominant $\Yim_t$-monomial. But we have $\tilde{Y}_{i,0}\tilde{Y}_{j,0}=\tilde{Y}_{j,0}\tilde{Y}_{i,0}$, so $F_t(\tilde{Y}_{i,0})F_t(\tilde{Y}_{j,0})-F_t(\tilde{Y}_{j,0})F_t(\tilde{Y}_{i,0})\in\mathfrak{K}_t^{\infty}$ has no dominant $\Yim_t$-monomial, so is equal to $0$.\qed

\noindent Denote, for $l\in\ZZ$, by $s_l:\Yim_t^{\infty}\rightarrow \Yim_t^{\infty}$ the endomorphism of $\ZZ[t^{\pm}]$-algebra such that $s_l(\tilde{Y}_{j,k})=\tilde{Y}_{j,k+l}$ (it is well-defined because the defining relations of $\Yim_t$ are invariant for $k\mapsto k+l$). If the hypothesis of the lemma \ref{copieun} are verified, we can define for $m\in t^{\ZZ}B$ :
$$E_t(m)=m (\overset{\rightarrow}{\underset{l\in\ZZ}{\prod}}\underset{i\in I}{\prod}\tilde{Y}_{i,l}^{u_{i,l}(m)})^{-1}\overset{\rightarrow}{\underset{l\in\ZZ}{\prod}}\underset{i\in I}{\prod}s_l(F_t(\tilde{Y}_{i,0}))^{u_{i,l}(m)}\in\mathfrak{K}_t^{\infty}$$\label{tetm}
because for $l\in\ZZ$ the product $\underset{i\in I}{\prod}s_l(F_t(\tilde{Y}_{i,0}))^{u_{i,l}(m)}$ is commutative (lemma \ref{copieun}).

\subsection{Step 1: Proof of $P(1)$ and $P(2)$}\label{petitdeux} The aim of this section is to prove $P(1)$ and $P(2)$. First we define a more precise property $Q(n)$ such that $Q(n)\Rightarrow P(n)$.

\subsubsection{The property $Q(n)$}\label{qn}

\begin{defi} For $n\geq 1$ denote $Q(n)$ the property ``for all semi-simple Lie-algebras $\Glie$ of rank $\text{rk}(\Glie)=n$, for all $i\in I$ there is a unique $F_t(\tilde{Y}_{i,0})\in\mathfrak{K}_t\cap \tilde{D}_{\tilde{Y}_{i,0}}$ such that $\tilde{Y}_{i,0}$ is the unique dominant $\Yim_t$-monomial of $F_t(\tilde{Y}_{i,0})$. Moreover $F_t(\tilde{Y}_{i,0})$ has the same monomials as $E(Y_{i,0})$''.\end{defi}

\noindent The property $Q(n)$ is more precise than $P(n)$ because it asks that $F_t(\tilde{Y}_{i,0})$ has only a finite number of monomials.

\begin{lem} For $n\geq 1$, the property $Q(n)$ implies the property $P(n)$.\end{lem}

\demo We suppose $Q(n)$ is true. In particular the section \ref{conste} enables us to construct $E_t(m)\in\mathfrak{K}_t^{\infty}$ for $m\in B$. The defining formula of $E_t(m)$ shows that it has the same monomials as $E(m)$. So $E_t(m)\in\tilde{D}_m$ and $E_t(m)\in\mathfrak{K}_t$.

\noindent Let us prove $P(n)$: let $m$ be in $B$. The uniqueness of $F_t(m)$ follows from lemma \ref{leasto}. Let $m_L=m>m_{L-1}>...>m_1$ be the dominant monomials of $D_m$ with a total ordering compatible with the partial ordering (it follows from lemma \ref{fini} that $D_m\cap B$ is finite). Let us prove by induction on $l$ the existence of $F_t(m_l)$. The unique dominant of $D_{m_1}$ is $m_1$ so $F_t(m_1)=E_t(m_1)\in\tilde{D}_{m_1}$. In general let $\lambda_1(t),...,\lambda_{l-1}(t)\in\ZZ[t^{\pm}]$ be the coefficient of the dominant $\Yim_t$-monomials $m_1,...,m_{l-1}$ in $E_t(m_l)$. We put:
$$F_t(m_l)=E_t(m_l)-\underset{r=1...l-1}{\sum}\lambda_r(t)F_t(m_r)$$
We see in the construction that $F_t(m)\in\tilde{D}_m$ because for $m'\in D_m$ we have $E_t(m')\in \tilde{D}_{m'}\subseteq\tilde{D}_m$ (lemma \ref{ordrep}).
\qed

\subsubsection{Cases $n=1,n=2$}\label{petit} We need the following general technical result:

\begin{prop}\label{cpfacile} Let $m$ be in $B$ such that all monomial $m'$ of $F(m)$ verifies : $\forall i\in I, m'\in B_i$ implies $\forall l\in\ZZ, u_{i,l}(m')\leq 1$ and for $1\leq r\leq 2r_i$ the set $\{l\in\ZZ/u_{i,r+2lr_i}(m')=1\}$ is a $1$-segment. Then $\pi^{-1}(F(m))\in\Yim_t$ is in $\mathfrak{K}_t$ and has a unique dominant monomial $m$.\end{prop}

\demo Let us write $F(m)=\underset{m'\in A}{\sum}\mu(m')m'$ ($\mu(m')\in\ZZ$). Let $i$ be in $I$ and consider the decomposition of $F(m)$ in $\mathfrak{K}_i$:
$$F(m)=\underset{m'\in B_i}{\sum}\mu(m')F_i(m')$$
But $\mu(m')\neq 0$ implies the hypothesis of lemma \ref{pidonne} is verified for $m'\in B_i$. So $\pi^{-1}(F_i(m'))=F_{i,t}(m')$. And:
$$\pi^{-1}(F(m))=\underset{m'\in B_i}{\sum}\mu(m')F_{i,t}(m')\in \mathfrak{K}_{i,t}$$\qed

\noindent For $n=1$ (section \ref{kernelun}), $n=2$ (section \ref{fin}), we can give explicit formula for the $E(Y_{i,0})=F(Y_{i,0})$. In particular we see that the hypothesis of proposition \ref{cpfacile} are verified, so: 
\begin{cor} The properties $Q(1)$, $Q(2)$ and so $P(1)$, $P(2)$ are true.\end{cor}

\noindent This allow us to start our induction in the proof of theorem \ref{con}.

\noindent In section \ref{acase} we will see other applications of proposition \ref{cpfacile}.

\noindent Note that the hypothesis of proposition \ref{cpfacile} are not verified for fundamental monomials $m=Y_{i,0}$ in general: for example for the $D_5$-case we have in $F(Y_{2,0})$ the monomial $Y_{3,3}^2Y_{5,4}^{-1}Y_{2,4}^{-1}Y_{4,4}^{-1}$.

\subsection{Step 2: consequences of the property $P(n)$}\label{consp} Let be $n\geq 1$. We suppose in this section that $P(n)$ is proved. We give some consequences of $P(n)$ which will be used in the proof of $P(r)$ ($r>n$).

\noindent Let $\mathfrak{K}_t^{\infty,f}$\label{ktinff} be the $\ZZ[t^{\pm}]$-submodule of $\mathfrak{K}_t^{\infty}$ generated by elements with a finite number of dominant $\Yim_t$-monomials. 
\begin{prop}\label{thth} We suppose $\text{rk}(\Glie)=n$. We have:
$$\mathfrak{K}_t^{\infty, f}=\underset{m\in B}{\bigoplus}\ZZ[t^{\pm}]F_t(m)\simeq\ZZ[t^{\pm}]^{(B)}$$
Moreover for $M\in A$, we have:
$$\mathfrak{K}_t^\infty\cap\tilde{D}_M=\underset{m\in B\cap D_M}{\bigoplus}\ZZ[t^{\pm}]F_t(m)\simeq\ZZ[t^{\pm}]^{B\cap D_M}$$
\end{prop}

\demo Let $\chi$ be in $\mathfrak{K}_t^{\infty, f}$ and $m_1,...,m_L\in B$ the dominant $\Yim_t$-monomials of $\chi$ and $\lambda_1(t),...,\lambda_L(t)\in\ZZ[t^{\pm}]$ their coefficients. It follows from lemma \ref{leasto} that $\chi=\underset{l=1...L}{\sum}\lambda_l(t)F_t(m_l)$.

\noindent Let us look at the second point: lemma \ref{ordrep} shows that $m\in B\cap D_M\Rightarrow F_t(m)\in\tilde{D}_M$. In particular the inclusion $\supseteq$ is clear. For the other inclusion we prove as in the first point that $\mathfrak{K}_t^\infty\cap\tilde{D}_M=\underset{m\in B\cap D_M}{\sum}\ZZ[t^{\pm}]F_t(m)$. We can conclude because it follows from lemma \ref{fini} that $D_M\cap B$ is finite.\qed

\noindent We recall that have seen in section \ref{infsum} that some infinite sum make sense in $\tilde{D}_M$.

\begin{cor}\label{recufond} We suppose $\text{rk}(\Glie)>n$ and let $J$ be a subset of $I$ such that $|J|=n$. For $m\in B_J$, there is a unique $F_{J,t}(m)\in\mathfrak{K}_{J,t}^{\infty}$ such that $m$ is the unique $J$-dominant $\Yim_t$-monomial of $F_{J,t}(m)$. Moreover $F_{J,t}(m)\in\tilde{D}_{m}$. 

\noindent For $M\in A$, the elements of $\mathfrak{K}_{J,t}^\infty\cap\tilde{D}_M$ are infinite sums $\underset{m\in B_J\cap D_M}{\sum}\lambda_m(t)F_{J,t}(m)$. In particular:
$$\mathfrak{K}_{J,t}^\infty\cap\tilde{D}_M\simeq\ZZ[t^{\pm}]^{B_J\cap D_M}$$
\end{cor}

\demo The uniqueness of $F_{J,t}(m)$ follows from lemma \ref{leasto}. Let us write $m=m_Jm'$ where 
\\$m_J=\underset{i\in J,l\in\ZZ}{\prod}Y_{i,l}^{u_{i,l}(m)}$. So $m_J$ is a dominant $\Yim_t$-monomial of $\ZZ[Y_{i,l}^{\pm}]_{i\in J,l\in\ZZ}$. In particular the proposition \ref{thth} with the algebra $\U_q(\hat{\Glie})_J$ of rank $n$ gives $m_J\chi$ where $\chi\in\ZZ[\tilde{A}_{i,l}^{\U_q(\hat{\Glie})_J,-1},t^{\pm}]_{i\in J,l\in\ZZ}$ (where for $i\in I,l\in\ZZ$, $\tilde{A}_{i,l}^{\U_q(\hat{\Glie})_J,\pm}=\beta_{I,J}(\tilde{A}_{i,l}^{\pm})$ where $\beta_{I,J}(\tilde{Y}_{i,l}^{\pm})=\delta_{i\in J}\tilde{Y}_{i,l}^{\pm}$). So we can put $F_t(m)=m\nu_{J,t}(\chi)$ where $\nu_{J,t}:\ZZ[\tilde{A}_{i,l}^{\U_q(\hat{\Glie})_J,-1},t^{\pm}]_{i\in J,l\in\ZZ}\rightarrow\Yim_t$ is the ring homomorphism such that $\nu_{J,t}(\tilde{A}_{i,l}^{\U_q(\hat{\Glie})_J,-1})=\tilde{A}_{i,l}^{-1}$.

\noindent The last assertion is proved as in proposition \ref{thth}.\qed

\subsection{Step 3: $t$-algorithm and end of the proof of theorem \ref{con}}\label{proof} In this section we explain why the $P(r)$ ($r<n$) imply $P(n)$. In particular we define the $t$-algorithm which constructs explicitly the $F_t(m)$.

\subsubsection{The induction} 

We prove the property $P(n)$ by induction on $n\geq 1$. It follows from section \ref{petitdeux} that $P(1)$ and $P(2)$ are true. Let be $n\geq 3$ and suppose that $P(r)$ is proved for $r<n$.

\noindent Let $m_+$ be in $B$ and $m_0=m_+>m_1>m_2>...$ the countable set $D_{m_+}$ with a total ordering compatible with the partial ordering.

\noindent For $J\varsubsetneq I$ and $m\in B_J$, it follows from $P(r)$ and corollary \ref{recufond} that there is a unique $F_{J,t}(m)\in\tilde{D}_m\cap\mathfrak{K}_{J,t}^{\infty}$ such that $m$ is the unique $J$-dominant monomial of $F_{J,t}(m)$ and that the elements of $\tilde{D}_{m_+}\cap \mathfrak{K}_{J,t}^{\infty}$ are the infinite sums of $\Yim_t^{\infty}$:  $\underset{m\in D_{m_+}\cap B_J}{\sum}\lambda_m(t)F_{J,t}(m)$ where $\lambda_m(t)\in\ZZ[t^{\pm}]$.

\noindent If $m\in A-B_J$, denote $F_{J,t}(m)=0$.

\subsubsection{Definition of the $t$-algorithm}\label{defialgo} For $r,r'\geq 0$ and $J\subsetneq I$ denote $[F_{J,t}(m_{r'})]_{m_r}\in\ZZ[t^{\pm}]$ the coefficient of $m_r$ in $F_{J,t}(m_{r'})$. 

\begin{defi} We call $t$-algorithm the following inductive definition of the sequences $(s(m_r)(t))_{r\geq 0}\in\ZZ[t^{\pm}]^{\NN}$, $(s_J(m_r)(t))_{r\geq 0}\in\ZZ[t^{\pm}]^{\NN}$ ($J\varsubsetneq I$)\label{smrt}:
$$s(m_0)(t)=1\text{ , }s_J(m_0)(t)=0$$
and for $r\geq 1, J\subsetneq I$:
$$s_J(m_r)(t)=\underset{r'<r}{\sum}(s(m_{r'})(t)-s_J(m_{r'})(t))[F_{J,t}(m_{r'})]_{m_r} $$
$$\text{ if }m_r\notin B_J, s(m_r)(t)=s_J(m_r)(t)$$
$$\text{ if }m_r\in B, s(m_r)(t)=0$$
\end{defi}

\noindent We have to prove that the $t$-algorithm defines the sequences in a unique way. We see that if $s(m_r),s_J(m_r)$ are defined for $r\leq R$ so are $s_J(m_{R+1})$ for $J\subsetneq I$. The $s_J(m_R)$ impose the value of $s(m_{R+1})$ and by induction the uniqueness is clear. We say that the $t$-algorithm is well defined to step $R$ if there exist $s(m_{r}), s_J(m_r)$ such that the formulas of the $t$-algorithm are verified for $r\leq R$.

\begin{lem} The $t$-algorithm is well defined to step $r$ if and only if:
$$\forall J_1,J_2\varsubsetneq I, \forall r'\leq r, m_{r'}\notin B_{J_1} \text{ and }m_{r'}\notin B_{J_2}\Rightarrow s_{J_1}(m_{r'})(t)=s_{J_2}(m_{r'})(t)$$\end{lem}

\demo If for $r'<r$ the $s(m_{r'})(t),s_J(m_{r'})(t)$ are well defined, so is $s_J(m_r)(t)$. If $m_r\in B$, $s(m_r)(t)=0$ is well defined. If $m_r\notin B$, it is well defined if and only if $\{s_J(m_r)(t)/m_r\notin B_J\}$ has one unique element.\qed

\subsubsection{The $t$-algorithm never fails} If the $t$-algorithm is well defined to all steps, we say that the $t$-algorithm never fails. In this section we show that the $t$-algorithm never fails.

\noindent If the $t$-algorithm is well defined to step $r$, for $J\varsubsetneq I$ we set:
$$\mu_J(m_r)(t)=s(m_r)(t)-s_J(m_r)(t)$$
$$\chi_J^r=\underset{r'\leq r}{\sum}\mu_{J}(m_{r'}(t))F_{J,t}(m_{r'})\in\mathfrak{K}_{J,t}^{\infty}$$

\begin{lem} If the $t$-algorithm is well defined to step $r$, for $J\subset I$ we have:
$$\chi_J^r\in (\underset{r'\leq r}{\sum} s(m_{r'})(t)m_{r'})+s_J(m_{r+1})(t)m_{r+1}+\underset{r'>r+1}{\sum}\ZZ[t^{\pm}]m_{r'}$$
For $J_1\subset J_2\subsetneq I$, we have:
$$\chi_{J_2}^r=\chi_{J_1}^r+\underset{r'>r}{\sum}\lambda_{r'}(t)F_{J_1,t}(m_{r'})$$
where $\lambda_{r'}(t)\in\ZZ[t^{\pm}]$. In particular, if $m_{r+1}\notin B_{J_1}$, we have $s_{J_1,t}(m_{r+1})=s_{J_2,t}(m_{r+1})$.
\end{lem}

\demo For $r'\leq r$ let us compute the coefficient $(\chi_J^r)_{m_{r'}}\in\ZZ[t^{\pm}]$ of $m_{r'}$ in $\chi_J^r$:
$$(\chi_J^r)_{m_{r'}}=\underset{r''\leq r'}{\sum}(s(m_{r''})(t)-s_J(m_{r''})(t))[F_{J,t}(m_{r''})]_{m_{r'}}$$
$$=(s(m_{r'})(t)-s_J(m_{r'})(t))[F_{J,t}(m_{r'})]_{m_{r'}}+\underset{r''<r'}{\sum}(s(m_{r''})(t)-s_J(m_{r''})(t))[F_{J,t}(m_{r''})]_{m_{r'}}$$
$$=(s(m_{r'})(t)-s_J(m_{r'})(t))+s_J(m_{r'})(t)=s(m_{r'})(t)$$
Let us compute the coefficient $(\chi_J^r)_{m_{r+1}}\in\ZZ[t^{\pm}]$ of $m_{r+1}$ in $\chi_J^r$:
$$(\chi_J^r)_{m_{r+1}}=\underset{r''<r+1}{\sum}(s(m_{r''})(t)-s_J(m_{r''})(t))[F_{J,t}(m_{r''})]_{m_{r+1}}=s_J(m_{r+1})$$
For the second point let $J_1\subset J_2\subsetneq I$. We have $\chi_{J_2}^r\in \mathfrak{K}_{J_1,t}^{\infty}\cap\tilde{D}_{m+}$ and it follows from $P(|J_1|)$ and corollary \ref{recufond} (or section \ref{petit} if $|J_1|\leq 2$) that we can introduce $\lambda_{m_{r'}}(t)\in\ZZ[t^{\pm}]$ such that :
$$\chi_{J_2}^r=\underset{r'\geq 0}{\sum}\lambda_{m_{r'}}(t)F_{J_1,t}(m_{r'})$$
We show by induction on $r'$ that for $r'\leq r$, $m_{r'}\in B_{J_1}\Rightarrow \lambda_{m_{r'}}(t)=\mu_{J_1}(m_{r'})(t)$. First we have $\lambda_{m_0}(t)=(\chi_{J_2}^r)_{m_0}=s(m_0)(t)=1=\mu_{J_1}(m_0)$. For $r'\leq r$: 
$$s(m_{r'})(t)=\lambda_{m_{r'}}(t)+\underset{r''<r'}{\sum}\lambda_{m_{r''}}(t)[F_{J_1,t}(m_{r''})]_{m_{r'}}$$
$$\lambda_{m_{r'}}(t)=s(m_{r'})(t)-\underset{r''<r'}{\sum}\mu_{J_1}(m_{r'})(t)[F_{J_1,t}(m_{r''})]_{m_{r'}}=s(m_{r'})(t)-s_{J_1}(m_{r'})(t)=\mu_{J_1}(m_{r'})(t)$$
For the last assertion if $m_{r+1}\notin B_{J_1}$, the coefficient of $m_{r+1}$ in $\underset{r'>r}{\sum}\ZZ[t^{\pm}]F_{J_1,t}(m_{r'})$ is 0, and $(\chi_{J_2}^r)_{m_{r+1}}=(\chi_{J_1}^r)_{m_{r+1}}$. It follows from the first point that $s_{J_1,t}(m_{r+1})=s_{J_2,t}(m_{r+1})$.\qed

\begin{lem}\label{nfail} The $t$-algorithm never fails.\end{lem}

\demo Suppose the sequence is well defined until the step $r-1$ and let $J_1,J_2\varsubsetneq I$ such that $m_r\notin B_{J_1}$ and $m_r\notin B_{J_2}$. Let $i$ be in $J_1$, $j$ in $J_2$ such that $m_r\notin B_i$ and $m_r\notin B_j$. Consider $J=\{i,j\}\varsubsetneq I$. The $\chi_J^{r-1},\chi_i^{r-1},\chi_j^{r-1}\in\Yim_t$ have the same coefficient $s(m_{r'})(t)$ on $m_{r'}$ for $r'\leq r-1$. Moreover:
$$s_i(m_r)(t)=(\chi_i^{r-1})_{m_r}\text{ , }s_j(m_r)(t)=(\chi_j^{r-1})_{m_r}\text{ , }s_J(m_r)(t)=(\chi_J^{r-1})_{m_r}$$
But $m_r\notin B_J$, so:
$$\chi_J^{r-1}=\underset{r'\leq r-1}{\sum}\mu_i(m_{r'})(t)F_{i,t}(m_{r'})+\underset{r'\geq r+1}{\sum}\lambda_{m_{r'}}(t)F_{i,t}(m_{r'})$$
So $(\chi_J^{r-1})_{m_r}=(\chi_i^{r-1})_{m_r}$ and we have $s_i(m_r)(t)=s_J(m_r)(t)$. In the same way we have $s_i(m_r)(t)=s_{J_1}(m_r)(t)$, $s_j(m_r)(t)=s_J(m_r)(t)$ and $s_j(m_r)(t)=s_{J_2}(m_r)(t)$. So we can conclude $s_{J_1}(m_r)(t)=s_{J_2}(m_r)(t)$.\qed

\subsubsection{Proof of $P(n)$} It follows from lemma \ref{nfail} that $\chi=\underset{r\geq 0}{\sum}s(m_r)(t)m_r\in\Yim_t^{\infty}$ is well defined.

\begin{lem}\label{conv} We have $\chi\in \mathfrak{K}_t^{\infty}\bigcap\tilde{D}_{m_+}$. Moreover the only dominant $\Yim_t$-monomial in $\chi$ is $m_0=m_+$.\end{lem}

\demo The defining formula of $\chi$ gives $\chi\in\tilde{D}_{m_+}$. Let $i$ be in $I$ and:
$$\chi_i=\underset{r\geq 0}{\sum}\mu_i(m_r)(t)F_{i,t}(m_r)\in\mathfrak{K}_{i,t}^{\infty}$$
Let us compute for $r\geq 0$ the coefficient of $m_r$ in $\chi-\chi_i$:
$$(\chi-\chi_i)_{m_r}=s(m_r)(t)-\underset{r'\leq r}{\sum}\mu_i(m_{r'})(t)[F_{i,t}(m_{r'})]_{m_r}$$
$$=s(m_r)(t)-s_i(m_r)(t)-\mu_i(m_r)(t)[F_{i,t}(m_{r})]_{m_r}=(s(m_r)(t)-s_i(m_r)(t))(1-[F_{i,t}(m_{r})]_{m_r})$$
We have two cases: 

if $m_r\in B_i$, we have $1-[F_{i,t}(m_{r})]_{m_r}=0$.

if $m_r\notin B_i$, we have $s(m_r)(t)-s_i(m_r)(t)=0$.

\noindent So $\chi=\chi_i\in\mathfrak{K}_{i,t}^{\infty}$, and $\chi\in \mathfrak{K}_t^{\infty}$.

\noindent The last assertion follows from the definition of the algorithm: for $r>0$, $m_r\in B\Rightarrow s(m_r)(t)=0$.\qed

\noindent This lemma implies:

\begin{cor} For $n\geq 3$, if the $P(r)$ ($r<n$) are true, then $P(n)$ is true.\end{cor}

\noindent In particular the theorem \ref{con} is proved by induction on $n$. 

\section{Morphism of $q,t$-characters and applications}\label{conschi}

\subsection{Morphism of $q,t$-characters}\label{aida}

\subsubsection{Definition of the morphism}

We set $\text{Rep}_t=\text{Rep}\otimes_{\ZZ}\ZZ[t^{\pm}]=\ZZ[X_{i,l},t^{\pm}]_{i\in I,l\in\ZZ}$\label{rept}. We say that $M\in\text{Rep}_t$ is a $\text{Rep}_t$-monomial if it is of the form $M=\underset{i\in I,l\in\ZZ}{\prod}X_{i,l}^{x_{i,l}}$ ($x_{i,l}\geq 0$). In this case denote $x_{i,l}(M)=x_{i,l}$. Recall the definition of the $E_t(m)$ (section \ref{conste}).

\begin{defi}\label{mqt} The morphism of $q,t$-characters is the $\ZZ[t^{\pm}]$-linear map $\chi_{q,t}:\text{Rep}_t\rightarrow \Yim_t^{\infty}$\label{chiqt} such that ($u_{i,l}\geq 0$):
$$\chi_{q,t}(\underset{i\in I,l\in\ZZ}{\prod}X_{i,l}^{u_{i,l}})=E_t(\underset{i\in I,l\in\ZZ}{\prod}Y_{i,l}^{u_{i,l}})$$
\end{defi}

\subsubsection{Properties of $\chi_{q,t}$}

\begin{thm}\label{axiomes} We have $\pi_+(\text{Im}(\chi_{q,t}))\subset\Yim$ and the following diagram is commutative:
$$\begin{array}{rcccl}
\text{Rep}&\stackrel{\chi_{q,t}}{\longrightarrow}&\text{Im}(\chi_{q,t})\\
\text{id}\downarrow &&\downarrow&\pi_+\\
\text{Rep}&\stackrel{\chi_q}{\longrightarrow}&\Yim
\end{array}$$ 
In particular the map $\chi_{q,t}$ is injective. The $\ZZ[t^{\pm}]$-linear map $\chi_{q,t}:\text{Rep}_t\rightarrow\Yim_t^{\infty}$ is characterized by the three following properties:

1) For a $\text{Rep}_t$-monomial $M$ define $m=\pi^{-1}(\underset{i\in I,l\in\ZZ}{\prod}Y_{i,l}^{x_{i,l}(M)})\in A$ and $\tilde{m}\in A_t$ as in section \ref{tildeat}. Then we have :
$$\chi_{q,t}(M)=\tilde{m}+\underset{m'<m}{\sum}a_{m'}(t)m'\text{ (where $a_{m'}(t)\in\ZZ[t^{\pm}]$)}$$

2) The image of $\text{Im}(\chi_{q,t})$ is contained in $\mathfrak{K}_t^{\infty}$.

3) Let $M_1,M_2$ be $\text{Rep}_t$-monomials such that $\text{max}\{l/\underset{i\in I}{\sum}x_{i,l}(M_1)>0\}\leq \text{min}\{l/\underset{i\in I}{\sum}x_{i,l}(M_2)>0\}$. We have :
$$\chi_{q,t}(M_1M_2)=\chi_{q,t}(M_1)\chi_{q,t}(M_2)$$
\end{thm}

\noindent Note that the properties $1,2,3$ are generalizations of the defining axioms introduced by Nakajima in \cite{Nab} for the $ADE$-case; in particular in the $ADE$-case $\chi_{q,t}$ is the morphism of $q,t$-characters constructed in \cite{Nab}.

\demo $\pi_+(\text{Im}(\chi_{q,t}))\subset\Yim$ means that only a finite number of $\Yim_t$-monomials of $E_t(m)$ have coefficient $\lambda(t)\notin (t-1)\ZZ[t^{\pm}]$. As $F_t(\tilde{Y}_{i,0})$ has no dominant $\Yim_t$-monomial other than $\tilde{Y}_{i,0}$, we have the same property for $\pi_+(F_t(\tilde{Y}_{i,0}))\in\mathfrak{K}^{\infty}$ and $\pi_+(F_t(\tilde{Y}_{i,0}))=E(Y_{i,0})\in \Yim$. As $\Yim$ is a subalgebra of $\Yim^{\infty}$ we get $\pi_+(E_t(m))\in\Yim$ with the help of the defining formula.

\noindent The diagram is commutative because $\pi_+\circ s_l=s_l\circ \pi_+$ and $\pi_+(F_t(\tilde{Y}_{i,0}))=E(Y_{i,0})$. It is proved by Frenkel, Reshetikhin in \cite{Fre} that $\chi_q$ is injective, so $\chi_{q,t}$ is injective.

\noindent Let us show that $\chi_{q,t}$ verifies the three properties:

1) By definition we have $\chi_{q,t}(M)=E_t(m)$. But $s_l(F_t(\tilde{Y}_{i,0}))=F_t(\tilde{Y}_{i,l})\in\tilde{D}(\tilde{Y}_{i,l})$. In particular $s_l(F_t(\tilde{Y}_{i,0}))$ is of the form $\tilde{Y}_{i,l}+\underset{m'<Y_{i,l}}{\sum}\lambda_{m'}(t)m'$ and we get the property for $E_t(m)$ by multiplication. 

2) We have $s_l(F_t(\tilde{Y}_{i,0}))=E_t(\tilde{Y}_{i,l})\in\mathfrak{K}_t^{\infty}$ and $\mathfrak{K}_t^{\infty}$ is a subalgebra of $\Yim_t^{\infty}$, so $\text{Im}(\chi_{q,t})\subset\mathfrak{K}_t^{\infty}$.

3) If we set $L=\text{max}\{l/\underset{i\in I}{\sum}x_{i,l}(M_1)>0\}$, $m_1=\underset{i\in I,l\in\ZZ}{\prod}Y_{i,l}^{x_{i,l}(M_1)}$, $m_2=\underset{i\in I,l\in\ZZ}{\prod}Y_{i,l}^{x_{i,l}(M_2)}$, we have:
$$E_t(m_1)=\overset{\rightarrow}{\underset{l\leq L}{\prod}}\underset{i\in I}{\prod}s_l(F_t(\tilde{Y}_{i,0}))^{x_{i,l}(M_1)}\text{ , }E_t(m_2)=\overset{\rightarrow}{\underset{l\geq L}{\prod}}\underset{i\in I}{\prod}s_l(F_t(\tilde{Y}_{i,0}))^{x_{i,l}(M_2)}$$
and in particular:
$$E_t(m_1m_2)=E_t(m_1)E_t(m_2)$$
\noindent Finally let $f:\text{Rep}_t\rightarrow\Yim_t^{\infty}$ be a $\ZZ[t^{\pm}]$-linear homomorphism which verifies properties 1,2,3. We saw that the only element of $\mathfrak{K}_t^{\infty}$ with highest weight monomial $\tilde{Y}_{i,l}$ is $s_l(F_t(\tilde{Y}_{i,0}))$. In particular we have $f(X_{i,l})=E_t(Y_{i,l})$. Using property 3, we get for $M\in\text{Rep}_t$ a monomial :
$$f(M)=\overset{\rightarrow}{\underset{l\in\ZZ}{\prod}}\underset{i\in I}{\prod}f(X_{i,l})^{u_{i,l}(m)}=\overset{\rightarrow}{\underset{l\in\ZZ}{\prod}}\underset{i\in I}{\prod}s_l(F_t(\tilde{Y}_{i,0}))^{u_{i,l}(m)}=\chi_{q,t}(M)$$
\qed

\subsection{Quantization of the Grothendieck Ring}\label{quanta} In this section we see that $\chi_{q,t}$ allows us to define a deformed algebra structure on $\text{Rep}_t$ generalizing the quantization of \cite{Nab}. The point is to show that $\text{Im}(\chi_{q,t})$ is a subalgebra of $\mathfrak{K}_t^{\infty}$.

\subsubsection{Generators of $\mathfrak{K}_t^{\infty,f}$} Recall the definition of $\mathfrak{K}_t^{\infty,f}$ in section \ref{consp}. For $m\in B$, all monomials of $E_t(m)$ are in $\{mA_{i_1,l_1}^{-1}...A_{i_K,l_K}^{-1}/k\geq 0,l_k\geq L\}$ where $L=\text{min}\{l\in\ZZ,\exists i\in I,u_{i,l}(m)>0\}$. So it follows from lemma \ref{fini} that $E_t(m)\in\mathfrak{K}_t^{\infty}$ has only a finite number of dominant $\Yim_t$-monomials, that is to say $E_t(m)\in\mathfrak{K}_t^{\infty,f}$.

\begin{prop}\label{diago} The $\ZZ[t^{\pm}]$-module $\mathfrak{K}_t^{\infty,f}$ is freely generated by the $E_t(m)$:
$$\mathfrak{K}_t^{\infty, f}=\underset{m\in B}{\bigoplus}\ZZ[t^{\pm}]E_t(m)\simeq\ZZ[t^{\pm}]^{(B)}$$
\end{prop}

\demo The $E_t(m)$ are $\ZZ[t^{\pm}]$-linearly independent and we saw $E_t(m)\in\mathfrak{K}_t^{\infty, f}$. It suffices to prove that the $E_t(m)$ generate the $F_t(m)$: let us look at $m_0\in B$ and consider $L=\text{min}\{l\in\ZZ,\exists i\in I,u_{i,l}(m_0)>0\}$. In the proof of lemma \ref{fini} we saw there is only a finite dominant monomials in $\{m_0A_{i_1,l_1}^{-v_{i_1,l_1}}...A_{i_R,l_R}^{-v_{i_R,l_R}}/R\geq 0,i_r\in I,l_r\geq L\}$. Let $m_0>m_1>...>m_D\in B$ be those monomials with a total ordering compatible with the partial ordering. In particular, for $0\leq d\leq D$ the dominant monomials of $E_t(m_d)$ are in $\{m_d,m_{d+1},...,m_D\}$. So there are elements $(\lambda_{d,d'}(t))_{0\leq d,d'\leq D}$ of $\ZZ[t^{\pm}]$ such that:
$$E_t(m_d)=\underset{d\leq d'\leq D}{\sum}\lambda_{d,d'}(t)F_t(m_{d'})$$
We have $\lambda_{d,d'}(t)=0$ if $d'<d$ and $\lambda_{d,d}(t)=1$. We have a triangular system with $1$ on the diagonal, so it is invertible in $\ZZ[t^{\pm}]$.\qed

\subsubsection{Construction of the quantization}

\begin{lem} $\mathfrak{K}_t^{\infty,f}$ is a subalgebra of $\mathfrak{K}_t^{\infty}$.\end{lem}

\demo It suffices to prove that for $m_1,m_2\in B$, $E_t(m_1)E_t(m_2)$ has only a finite number of dominant $\Yim_t$-monomials. But $E_t(m_1)E_t(m_2)$ has the same monomials as $E_t(m_1m_2)$.\qed

\noindent It follows from proposition \ref{diago} that $\chi_{q,t}$ is a $\ZZ[t^{\pm}]$-linear isomorphism between $\text{Rep}_t$ and $\mathfrak{K}_t^{\infty,f}$. So we can define:
\begin{defi} The associative deformed $\ZZ[t^{\pm}]$-algebra structure on $\text{Rep}_t$ is defined by:
$$\forall \lambda_1,\lambda_2\in\text{Rep}_t,\lambda_1*\lambda_2=\chi_{q,t}^{-1}(\chi_{q,t}(\lambda_1)\chi_{q,t}(\lambda_2))$$\label{star}\end{defi}
\subsubsection{Examples: $sl_2$-case} We make explicit computation of the deformed multiplication in the $sl_2$-case:

\begin{prop} In the $sl_2$-case, the deformed algebra structure on $\text{Rep}_t=\ZZ[X_l,t^{\pm}]_{l\in\ZZ}$ is given by:
$$X_{l_1}*X_{l_2}*...*X_{l_m}=X_{l_1}X_{l_2}...X_{l_m}\text{ if $l_1\leq l_2\leq ...\leq l_m$}$$
$$X_{l}*X_{l'}=t^{\gamma}X_{l}X_{l'}=t^{\gamma}X_{l'}*X_l\text{ if $l>l'$ and $l\neq l'+2$}$$
$$X_{l}*X_{l-2}=t^{-2}X_{l}X_{l-2}+t^{\gamma}(1-t^{-2})=t^{-2}X_{l-2}*X_{l}+(1-t^{-2})$$
where $\gamma\in\ZZ$ is defined by $\tilde{Y}_l\tilde{Y}_{l'}=t^{\gamma}\tilde{Y}_{l'}\tilde{Y}_l$.
\end{prop}

\demo For $l\in\ZZ$ we have the $q,t$-character of the fundamental representation $X_l$:
$$\chi_{q,t}(X_l)=\tilde{Y}_l+\tilde{Y}_{l+2}^{-1}=\tilde{Y}_l(1+t\tilde{A}_{l+1}^{-1})$$
The first point of the proposition follows immediately from the definition of $\chi_{q,t}$. For example, for $l,l'\in\ZZ$ we have:
$$\chi_{q,t}(X_lX_{l'})=\chi_{q,t}(X_{\text{min}(l,l')})\chi_{q,t}(X_{\text{max}(l,l')})$$
In particular if $l\leq l'$, we have $X_l*X_{l'}=X_lX_{l'}$. Suppose now that $l>l'$ and introduce $\gamma\in\ZZ$ such that $\tilde{Y}_l\tilde{Y}_{l'}=t^{\gamma}\tilde{Y}_{l'}\tilde{Y}_l$. We have:
$$\chi_{q,t}(X_l)\chi_{q,t}(X_l')=\tilde{Y}_l(1+t\tilde{A}_{l+1}^{-1})\tilde{Y}_{l'}(1+t\tilde{A}_{l'+1}^{-1})$$
$$=t^{\gamma}\tilde{Y}_{l'}\tilde{Y}_l+t^{\gamma +1}\tilde{Y}_{l'}\tilde{Y}_l\tilde{A}_{l+1}^{-1}+t^{\gamma +1+2\delta_{l,l'+2}}\tilde{Y}_{l'}\tilde{A}_{l'+1}^{-1}\tilde{Y}_l+t^{\gamma+ 2}\tilde{Y}_{l'}\tilde{A}_{l'+1}^{-1}\tilde{Y}_{l}\tilde{A}_{l+1}^{-1}$$
$$=t^{\gamma}\chi_{q,t}(X_{l'}X_l)+t^{\gamma +1}(t^{2\delta_{l,l'+2}}-1)\tilde{Y}_{l'}\tilde{A}_{l'+1}^{-1}\tilde{Y}_l$$
If $l\neq l'+2$ we get $X_l*X_{l'}=t^{\gamma}X_{l'}*X_l$. If $l=l'+2$, we have:
$$\tilde{Y}_{l'}\tilde{A}_{l'+1}^{-1}\tilde{Y}_{l'+2}=t^{-1}\tilde{Y}_{l'+2}^{-1}\tilde{Y}_{l'+2}=t^{-1}$$
But $t^2\tilde{Y}_l\tilde{Y}_{l-2}=\tilde{Y}_{l-2}\tilde{Y}_l$, so $X_l*X_{l-2}=t^{-2}X_{l-2}*X_l+t^{-2}(t^2-1)$.
\qed

\noindent Note that $\gamma$ were computed in section \ref{varva}.

\noindent We see that the new $\ZZ[t^{\pm}]$-algebra structure is not commutative and not even twisted polynomial.

\subsection{An involution of the Grothendieck ring}\label{invo} In this section we construct an antimultiplicative involution of the Grothendieck ring $\text{Rep}_t$. The construction is motivated by the point view adopted in this article : it is just replacing $c_{|l|}$ by $-c_{|l|}$. In the $ADE$-case such an involution were introduced Nakajima \cite{Nab} with different motivations.

\subsubsection{An antihomomorphism of $\mathcal{H}$}

\begin{lem} There is a unique $\CC$-linear isomorphism of $\mathcal{H}$ which is antimultiplicative and such that:
$$\overline{c_m}=-c_m\text{ , }\overline{a_i[r]}=a_i[r]\text{ ($m>0,i\in I,r\in\ZZ-\{0\}$)}$$
Moreover it is an involution.
\end{lem}

\demo It suffices to show it is compatible with the defining relations of $\mathcal{H}$ ($i,j\in I,m,r\in\ZZ -\{0\}$):
$$\overline{[a_i[m],a_j[r]]}=\overline{a_i[m]a_j[r]}-\overline{a_j[r]a_i[m]}=-[a_i[m],a_j[r]]$$
$$\overline{\delta_{m,-r}(q^m-q^{-m})B_{i,j}(q^m)c_{|m|}}=-\delta_{m,-r}(q^m-q^{-m})B_{i,j}(q^m)c_{|m|}$$
For the last assertion, we have $\overline{\overline{c_m}}=c_m$ and $\overline{\overline{a_i[r]}}=a_i[r]$, and an algebra morphism which fixes the generators is the identity.\qed

\noindent It can be naturally extended to an antimultiplicative $\CC$-isomorphism of $\mathcal{H}_h$. 

\begin{lem}\label{stbl} The $\ZZ$-subalgebra $\Yim_u\subset \mathcal{H}_h$ verifies $\overline{\Yim_u}\subset\Yim_u$.\end{lem}

\demo It suffices to check on the generators of $\Yim_u$ ($R\in\mathfrak{U},i\in I,l\in\ZZ$):
$$\overline{t_R}=\text{exp}(\underset{m>0}{\sum}h^{2m}R(q^m)(-c_m))=t_{-R}$$ 
$$\overline{\tilde{Y}_{i,l}}=\text{exp}(\underset{m>0}{\sum}h^m y_i[-m]q^{-lm})\text{exp}(\underset{m>0}{\sum}h^m y_i[m]q^{lm})$$
$$=\text{exp}(\underset{m>0}{\sum}h^{2m}[y_i[-m],y_i[m]])\tilde{Y}_{i,l}
=t_{-\tilde{C}_{i,i}(q)(q_i-q_i^{-1})}\tilde{Y}_{i,l}\in\Yim_u$$
$$\overline{\tilde{Y}_{i,l}}^{-1}=(\overline{\tilde{Y}_{i,l}})^{-1}=t_{\tilde{C}_{i,i}(q)(q_i-q_i^{-1})}\tilde{Y}_{i,l}^{-1}\in\Yim_u$$
\qed

\subsubsection{Involution of $\Yim_t$} As for $R,R'\in\mathfrak{U}$, we have $\pi_0(R)=\pi_0(R')\Leftrightarrow\pi_0(-R)=\pi_0(-R')$, the involution of $\Yim_u$ (resp. of $\mathcal{H}_h$) is compatible with the defining relations of $\Yim_t$ (resp. $\mathcal{H}_t$). We get a $\ZZ$-linear involution of $\Yim_t$ (resp. of $\mathcal{H}_t$). For $\lambda,\lambda'\in\Yim_t,\alpha\in\ZZ$, we have:
$$\overline{\lambda.\lambda'}=\overline{\lambda'}.\overline{\lambda}\text{ , }\overline{t^{\alpha}\lambda}=t^{-\alpha}\overline{\lambda}$$

\noindent Note that in $\Yim_u$ for $i\in I,l\in\ZZ$:
$$\overline{\tilde{A}_{i,l}}=\text{exp}(\underset{m>0}{\sum}h^m a_i[-m]q^{-lm})\text{exp}(\underset{m>0}{\sum}h^m a_i[m]q^{lm})$$
$$=\text{exp}(\underset{m>0}{\sum}h^{2m}[a_i[-m],a_i[m]]c_m)\tilde{A}_{i,l}=t_{(-q_i^2+q_i^{-2})}\tilde{A}_{i,l}$$
So in $\Yim_t$ we have $\overline{\tilde{A}_{i,l}}=\tilde{A}_{i,l}$ and $\overline{\tilde{A}_{i,l}^{-1}}=\tilde{A}_{i,l}^{-1}$.

\subsubsection{The involution of deformed bimodules}

\begin{lem} For $i\in I$, the $\Yim_{i,u}\subset \mathcal{H}_h$ verifies $\overline{\Yim_{i,u}}\subset\Yim_{i,u}$.\end{lem}

\demo First we compute for $i\in I,l\in\ZZ$:
$$\overline{\tilde{S}_{i,l}}=\text{exp}(\underset{m>0}{\sum}h^m\frac{a_i[-m]}{q_i^{-m}-q_i^{m}}q^{-lm})\text{exp}(\underset{m>0}{\sum}h^m\frac{a_{i}[m]}{q_i^{m}-q_i^{-m}}q^{lm})$$
$$=\text{exp}(\underset{m>0}{\sum}h^{2m}\frac{[a_i[-m],a_i[m]]}{-(q_i^{-m}-q_i^{m})^2}c_m)\tilde{S}_{i,l}=t_{\frac{q_i+q_i^{-1}}{q_i-q_i^{-1}}}\tilde{S}_{i,l}\in\Yim_{i,u}$$
Now for $\lambda\in\Yim_u$, we have $\overline{\lambda.\tilde{S}_{i,l}}=t_{\frac{q_i+q_i^{-1}}{q_i-q_i^{-1}}}\tilde{S}_{i,l}\overline{\lambda}$. But it is in $\Yim_{i,u}$ because $\overline{\lambda}\in\Yim_u$ (lemma \ref{stbl}) and $\Yim_{i,u}$ is a $\Yim_u$-subbimodule of $\mathcal{H}_h$ (lemma \ref{currents}).\qed

\noindent In $\mathcal{H}_t$ we have  $\overline{\tilde{S}_{i,l}}=t\tilde{S}_{i,l}$ because $\pi_0(\frac{q_i+q_i^{-1}}{q_i-q_i^{-1}})=1$. As said before we get a $\ZZ$-linear involution of $\Yim_{i,t}$ such that:
$$\overline{\lambda \tilde{S}_{i,l}}=t\tilde{S}_{i,l}\overline{\lambda}$$

\noindent We introduced such an involution in \cite{Her01}. With this new point of view, the compatibility with the relation $\tilde{A}_{i,l-r_i}\tilde{S}_{i,l}=t^{-1}\tilde{S}_{i,l+r_i}$ is a direct consequence of lemma \ref{currents} and needs no computation; for example:
$$\overline{\tilde{A}_{i,l-r_i}\tilde{S}_{i,l}}=t\tilde{S}_{i,l}\tilde{A}_{i,l-r_i}=t^3\tilde{A}_{i,l-r_i}\tilde{S}_{i,l}=t^2\tilde{S}_{i,l+r_i}$$
$$\overline{t^{-1}\tilde{S}_{i,l+r_i}}=t\overline{\tilde{S}_{i,l+r_i}}=t^2\tilde{S}_{i,l+r_i}$$

\subsubsection{The induced involution of $\text{Rep}_t$}

\begin{lem} For $i\in I$, the subalgebra $\mathfrak{K}_{i,t}\subset\Yim_t$ verifies $\overline{\mathfrak{K}_{i,t}}\subset\mathfrak{K}_{i,t}$.\end{lem}

\demo Suppose $\lambda\in\mathfrak{K}_{i,t}$, that is to say $S_{i,t}(\lambda)=0$. So $\overline{(t^2-1)S_{i,t}(\lambda)}=0$ and:
$$\underset{l\in\ZZ}{\sum}(\overline{\tilde{S}_{i,l}\lambda}-\overline{\lambda\tilde{S}_{i,l}})=0\Rightarrow t\underset{l\in\ZZ}{\sum}(\overline{\lambda}\tilde{S}_{i,l}-\tilde{S}_{i,l}\overline{\lambda})=0$$
So $t(1-t^2)S_{i,t}(\overline{\lambda})=0$ and $\overline{\lambda}\in\mathfrak{K}_{i,t}$.\qed

\noindent Note that $\chi\in\Yim_t$ has the same monomials as $\overline{\chi}$, that is to say if $\chi=\underset{m\in A}{\sum}\lambda(t)m$ and $\overline{\chi}=\underset{m\in A}{\sum}\mu(t)m$, we have $\lambda(t)\neq 0\Leftrightarrow\mu(t)\neq 0$. In particular we can naturally extend our involution to an antimultiplicative involution on $\Yim_t^{\infty}$. Moreover we have $\overline{\mathfrak{K}_t^{\infty}}\subset \mathfrak{K}_t^{\infty}$ and $\overline{\mathfrak{K}_t^{\infty,f}}=\overline{\text{Im}(\chi_{q,t})}\subset \text{Im}(\chi_{q,t})$. So we can define:

\begin{defi} The $\ZZ$-linear involution of $\text{Rep}_t$ is defined by:
$$\forall \lambda\in\text{Rep}_t\text{ , }\overline{\lambda}=\chi_{q,t}^{-1}(\overline{\chi_{q,t}(\lambda)})$$\end{defi}

\subsection{Analogues of Kazhdan-Lusztig polynomials} In this section we define analogues of Kazhdan-Lusztig polynomials (see \cite{kalu}) with the help of the antimultiplicative involution of section \ref{invo} in the same spirit Nakajima did for the $ADE$-case \cite{Nab}. Let us begin we some technical properties of the action of the involution on monomials.

\subsubsection{Invariance of monomials}\label{invmon} We recall that the $\Yim_t^A$-monomials are products of the $\tilde{A}_{i,l}^{-1}$ ($i\in I, l\in\ZZ$).

\begin{lem} For $M$ a $\Yim_t$-monomial and $m$ a $\Yim_t^A$-monomial there is a unique $\alpha(M,m)\in\ZZ$\label{alpham} such that $\overline{t^{\alpha(M,m)}Mm}=t^{\alpha(M,m)}\overline{M}m$.\end{lem}

\demo Let $\beta\in\ZZ$ such that $\overline{m}=t^{\beta}m$. We have $\overline{Mm}=\overline{m}\overline{M}=t^{\beta+\gamma}\overline{M}m$ where $\gamma\in 2\ZZ$ (section \ref{notuil}). So it suffices to prove that $\beta\in 2\ZZ$. 

\noindent Let us compute $\beta$. Let $\pi_+(m)=\underset{i\in I,l\in\ZZ}{\prod}A_{i,l}^{-v_{i,l}}$. In $\Yim_u$ we have $\pi_+(m)\pi_-(m)=t_R\pi_-(m)\pi_+(m)$ where $\pi_0(R)=\beta$ and:
$$R(q)=\underset{i,j\in I,r,r'\in\ZZ}{\sum}v_{i,r}v_{j,r'}\underset{l>0}{\sum}q^{lr-lr'}\frac{[a_i[l],a_j[-l]]}{c_l}$$
where for $l>0$ we set $\frac{[a_i[l],a_j[-l]]}{c_l}=B_{i,j}(q^l)(q^{l}-q^{-l})\in\ZZ[q^{\pm}]$ which is antisymmetric. For $i=j$, we have the term:
$$\underset{r,r'\in\ZZ}{\sum}v_{i,r}v_{i,r'}\underset{l>0}{\sum}q^{lr-lr'}\frac{[a_i[l],a_i[-l]]}{c_l}$$
$$=\underset{l>0}{\sum}(\underset{\{r,r'\}\subset\ZZ, r\neq r'}{\sum}v_{i,r}(m)v_{i,r'}(m)(q^{l(r-r')}+q^{l(r'-r)})+\underset{r\in\ZZ}{\sum}v_{i,r}(m)^2)\frac{[a_i[l],a_i[-l]]}{c_l}$$
It is antisymmetric, so it has no term in $q^0$. So $\pi_0(R)=\pi_0(R')$ where $R'$ is the sum of the contributions for $i\neq j$: 
$$\underset{r,r'\in\ZZ}{\sum}v_{i,r}(m)v_{j,r'}(m)\underset{l>0}{\sum}q^{lr-lr'}(\frac{[a_i[l],a_j[-l]]}{c_l}+\frac{[a_j[l],a_i[-l]]}{c_l})$$
$$=2\underset{r,r'\in\ZZ}{\sum}v_{i,r}(m)v_{j,r'}(m)\underset{l>0}{\sum}q^{lr-lr'}\frac{[a_i[l],a_j[-l]]}{c_l}$$
In particular $\pi_0(R')\in 2\ZZ$.\qed

\noindent For $M$ a $\Yim_t$-monomial denote $A^{\text{inv}}_M=\{t^{\alpha(m,M)}Mm/\text{$m$ $\Yim_t^A$-monomial}\}$\label{ainv}. In particular for $m'\in A^{\text{inv}}_M$ we have $\overline{m'}{m'}^{-1}=\overline{M}M^{-1}$.

\subsubsection{The polynomials} For $M$ a $\Yim_t$-monomial, denote $B^{\text{inv}}_M=t^{\ZZ}B\cap A^{\text{inv}}_M\label{binv}$. 

\begin{thm}\label{expol} For $m\in t^{\ZZ}B$ there is a unique $L_t(m)\in\mathfrak{K}_t^{\infty}$\label{tltm} such that:
$$\overline{L_t(m)}=(\overline{m}m^{-1})L_t(m)$$
$$E_t(m)=L_t(m)+\underset{m'<m, m' \in B^{\text{inv}}_m}{\sum}P_{m',m}(t)L_t(m')$$
where $P_{m',m}(t)\in t^{-1}\ZZ[t^{-1}]$.
\end{thm}

\noindent Those polynomials $P_{m',m}(t)$ are called analogues to Kazhdan-Lusztig polynomials and the $L_t(m)$ ($m\in B$) for a canonical basis of $\mathfrak{K}_t^{f, \infty}$. Such polynomials were introduced by Nakajima \cite{Nab} for the $ADE$-case.

\demo First consider $\overline{F_t(m)}$: it is in $\mathfrak{K}_t^{\infty}$ and has only one dominant $\Yim_t$-monomial $\overline{m}$, so $\overline{F_t(m)}=\overline{m}m^{-1}F_t(m)$.

\noindent Let be $m=m_L>m_{L-1}>...>m_0$ the finite set $t^{\ZZ}D(m)\cap B^{\text{inv}}_m$ (see lemma \ref{ordrep}) with a total ordering compatible with the partial ordering. Note that it follows from section \ref{invmon} that for $L\geq l\geq 0$, we have $\overline{m_l}m_l^{-1}=\overline{m}m^{-1}$.

\noindent We have $E_t(m_0)=F_t(m_0)$ and so $\overline{E_t(m_0)}=\overline{m_0}m_0^{-1}E_t(m_0)$. As $B_{m_0}^{\text{inv}}=\{m_0\}$, we have $L_t(m_0)=E_t(m_0)$. We suppose by induction that the $L_t(m_l)$ ($L-1\geq l\geq 0$) are uniquely and well defined. In particular $m_l$ is of highest weight in $L_t(m_l)$, $\overline{L_t(m_l)}=\overline{m_l}m_l^{-1}L_t(m_l)=\overline{m}m^{-1}L_t(m_l)$, and we can write:
$$\tilde{D}_t(m_L)\cap \mathfrak{K}_t^{\infty}=\ZZ[t^{\pm}]F_t(m_L)\oplus\underset{0\leq l\leq L-1}{\bigoplus}\ZZ[t^{\pm}]L_t(m_l)$$
In particular consider $\alpha_{l,L}(t)\in\ZZ[t^{\pm}]$ such that:
$$E_t(m)=F_t(m)+\underset{l<L}{\sum}\alpha_{l,L}(t)L_t(m_l)$$
We want $L_t(m)$ of the form :
$$L_t(m)=F_t(m)+\underset{l<L}{\sum}\beta_{l,L}(t)L_t(m_l)$$
The condition $\overline{L_t(m)}=\overline{m}m^{-1}mL_t(m)$ means that the $\beta_{l,L}(t)$ are symmetric. The condition $P_{m',m}(t)\in t^{-1}\ZZ[t^{-1}]$ means $\alpha_{l,L}(t)-\beta_{l,L}(t)\in t^{-1}\ZZ[t^{-1}]$. So it suffices to prove that those two conditions uniquely define the $\beta_{l,L}(t)$: let us write $\alpha_{l,L}(t)=\alpha_{l,L}^+(t)+\alpha_{l,L}^0(t)+\alpha_{l,L}^-(t)$ (resp. $\beta_{l,L}(t)=\beta_{l,L}^+(t)+\beta_{l,L}^0(t)+\beta_{l,L}^-(t)$) where $\alpha_{l,L}^{\pm}(t)\in t^{\pm}\ZZ[t^{\pm}]$ and $\alpha_{l,L}^0(t)\in\ZZ$ (resp. for $\beta$). The condition $\alpha_{l,L}(t)-\beta_{l,L}(t)\in t^{-1}\ZZ[t^{-1}]$ means $\beta_{l,L}^0(t)=\alpha_{l,L}^0(t)$ and $\beta_{l,L}^-(t)=\alpha_{l,L}^-(t)$. The symmetry of $\beta_{l,L}(t)$ means $\beta_{l,L}^+(t)=\beta_{l,L}^-(t^{-1})=\alpha_{l,L}^-(t^{-1})$.\qed

\subsubsection{Examples for $\Glie=sl_2$}\label{exltcalc} In this section we suppose that $\Glie=sl_2$.

\begin{prop}\label{calcex} Let $m\in t^{\ZZ}B$ such that $\forall l\in\ZZ, u_l(m)\leq 1$. Then $L_t(m)=F_t(m)$. Moreover:
$$E_t(m)=L_t(m)+\underset{m'<m/m'\in B_m^{\text{inv}}}{\sum}t^{-R(m')}L_t(m')$$
where $R(m')\geq 1$ is given by $\pi_+(m'm^{-1})=A_{i_1,l_1}^{-1}...A_{i_R,l_R}^{-1}$. In particular for $m'\in B^{\text{inv}}_m$ such that $m'<m$ we have $P_{m',m}(t)=t^{-R(m')}$.\end{prop}

\demo Note that a dominant monomial $m'<m$ verifies $\forall l\in\ZZ, u_l(m')\leq 1$ and appears in $E_t(m)$. We know that $\tilde{D}_m\cap\mathfrak{K}_t=\underset{m'\in t^{\ZZ}D_m\cap B^{\text{inv}}_m}{\bigoplus}\ZZ[t^{\pm}]F_t(m')$. We can introduce $P_{m',m}(t)\in\ZZ[t^{\pm}]$ such that:
$$E_t(m)=F_t(m)+\underset{m'\in t^{\ZZ}D_m\cap B^{\text{inv}}-\{m\} }{\sum}P_{m',m}(t)F_t(m')$$
So by induction it suffices to show that $P_{m',m}(t)\in t^{-1}\ZZ[t^{-1}]$.  

\noindent $P_{m',m}(t)$ is the coefficient of $m'$ in $E_t(m)$. A dominant $\Yim_t$-monomial $M$ which appears in $E_t(m)$ is of the form:
$$M=m(m_1...m_{R+1})^{-1}m_1t\tilde{A}_{l_1}^{-1}m_2t\tilde{A}_{l_2}^{-1}m_3...t\tilde{A}_{l_R}^{-1}m_{R+1}$$
where $l_1<...<l_R\in\ZZ$ verify $\{l_r+2,l_r-2\}\cap\{l_1,...,l_{r-1},l_{r+1},...,l_R\}$ is empty, $u_{l_r-1}(m)=u_{l_r+1}(m)=1$ and we have set $m_r=\underset{l_{r-1}<l\leq l_r}{\overset{\rightarrow}{\prod}}\tilde{Y}_l^{u_l(m)}$. Such a monomial appears one time in $E_t(m)$. In particular $P_{m',m}(t)=t^{\alpha}$ where $\alpha\in\ZZ$ is given by $M=t^{\alpha}m'$ that is to say $\overline{M}M^{-1}=t^{-2\alpha}{m'}^{-1}m'=t^{-2\alpha}m^{-1}m$. So we compute:

\noindent $\overline{M}M^{-1}=t^{-2R}\overline{m_{R+1}}\tilde{A}_{l_R}^{-1}\overline{m_R}...\tilde{A}_{l_1}^{-1}\overline{m_1}(\overline{m}_1^{-1}...\overline{m}_{R+1}^{-1})\overline{m}m_{R+1}^{-1}\tilde{A}_{l_R}m_R^{-1}...\tilde{A}_{l_1}m_1^{-1}(m_1...m_{R+1})m^{-1}
\\=t^{-2R}t^{4R}\tilde{A}_{l_R}^{-1}...\tilde{A}_{l_1}^{-1}\overline{m}\tilde{A}_{l_R}...\tilde{A}_{l_1}m^{-1}
\\=t^{2R}\tilde{A}_{l_R}^{-1}...\tilde{A}_{l_1}^{-1}\tilde{A}_{l_R}...\tilde{A}_{l_1}\overline{m}m^{-1}=t^{2R}\overline{m}m^{-1}$\qed

\noindent Let us look at another example $m=\tilde{Y}_0^2\tilde{Y}_2$. We have:
$$E_t(m)=L_t(m)+t^{-2}L_t(m')$$
where $m'=t\tilde{Y}_0^2\tilde{Y}_2\tilde{A}_1^{-1}\in B_{m}^{\text{inv}}$ and:
$$L_t(m)=F_t(\tilde{Y}_0)F_t(\tilde{Y}_0\tilde{Y}_2)=\tilde{Y}_0(1+t\tilde{A}_1^{-1})\tilde{Y}_0\tilde{Y}_2(1+t\tilde{A}_3^{-1}(1+t\tilde{A}_1^{-1}))$$
$$L_t(m')=F_t(m')=t\tilde{Y}_0^2\tilde{Y}_2\tilde{A}_1^{-1}(1+t\tilde{A}_1^{-1})$$
Indeed the dominant monomials appearing in $E_t(m)$ are $m$ and $\tilde{Y}_0t\tilde{A}_1^{-1}\tilde{Y}_0\tilde{Y}_2+\tilde{Y}_0^2t\tilde{A}_1^{-1}\tilde{Y}_2=(1+t^{-2})m'$.

\noindent In particular: $P_{m',m}(t)=t^{-2}$.

\subsubsection{Example in non-simply laced case}\label{exnons} We suppose that $C=\begin{pmatrix}2 & -2\\-1 & 2\end{pmatrix}$ and $m=\tilde{Y}_{2,0}\tilde{Y}_{1,5}$. The formulas for $E_t(\tilde{Y}_{2,0})$ and $E_t(\tilde{Y}_{1,5})$ are given is section \ref{fin}. We have:
$$E_t(m)=L_t(m)+t^{-1}L_t(m')$$
where $m'=t\tilde{Y}_{2,0}\tilde{Y}_{1,5}\tilde{A}_{2,2}^{-1}\tilde{A}_{1,4}^{-1}\in B_m^{\text{inv}}$ and:
$$L_t(m')=F_t(m')=t\tilde{Y}_{2,0}\tilde{Y}_{1,5}\tilde{A}_{2,2}^{-1}\tilde{A}_{1,4}^{-1}(1+t\tilde{A}_{1,2}^{-1}(1+t\tilde{A}_{2,4}^{-1}(1+t\tilde{A}_{1,6}^{-1})))$$
Indeed the dominant monomials appearing in $E_t(m)$ are $m=\tilde{Y}_{2,0}\tilde{Y}_{1,5}$ and $\tilde{Y}_{2,0}t\tilde{A}_{2,2}^{-1}t\tilde{A}_{1,4}^{-1}\tilde{Y}_{1,5}=t^{-1}m'$. 

\noindent In particular $P_{m',m}(t)=t^{-1}$.

\section{Questions and conjectures}\label{quest}

\subsection{Positivity of coefficients}\label{acase}

\begin{prop}\label{cpaconj} If $\Glie$ is of type $A_n$ ($n\geq 1$), the coefficients of $\chi_{q,t}(Y_{i,0})$ are in $\NN[t^{\pm}]$.\end{prop}

\demo We show that for all $i\in I$ the hypothesis of proposition \ref{cpfacile} for $m=Y_{i,0}$ are verified; in particular the property $Q$ of section \ref{qn} will be verified.

\noindent Let $i$ be in $I$. For $j\in I$, let us write $E(Y_{i,0})=\underset{m\in B_j}{\sum}\lambda_j(m)E_j(m)\in\mathfrak{K}_j$ where $\lambda_j(m)\in\ZZ$. Let $D$ be the set $D=\{\text{monomials of $E_j(m)$ }/j\in I,m\in B_j,\lambda_j(m)\neq 0\}$. It suffices to prove that for $j\in I$, $m\in B_j\cap D\Rightarrow u_j(m)\leq 1$ (because proposition \ref{cpfacile} implies that for all $i\in I$, $F_t(\tilde{Y}_{i,0})=\pi^{-1}(E(Y_{i,0}))$).

\noindent As $E(Y_{i,0})=F(Y_{i,0})$, $Y_{i,0}$ is the unique dominant $\Yim$-monomial in $E(Y_{i,0})$. So for a monomial $m\in D$ there is a finite sequence $\{m_0=Y_{i,0},m_1,...,m_R=m\}$ such that for all $1\leq r\leq R$, there is $r'<r$ and $j\in I$ such that $m_{r'}\in B_j$ and for $r'<r''\leq r$, $m_{r''}$ is a monomial of $E_j(m_{r'})$ and $m_{r''}m_{r''-1}^{-1}\in\{A_{j,l}^{-1}/l\in\ZZ\}$. Such a sequence is said to be adapted to $m$. Suppose there is $j\in I$ and $m\in B_j\cap D$ such that $u_j(m)\geq 2$. So there is $m'\leq m$ in $D\cap B_j$ such that $u_j(m)=2$. So we can consider $m_0\in D$ such that there is $j_0\in I$, $m_0\in B_{j_0}$, $u_{j_0}(m)\geq 2$ and for all $m'<m_0$ in $D$ we have $\forall j\in I,m'\in B_j\Rightarrow u_j(m')\leq 1$. Let us write:
$$m_0=Y_{j_0,q^{l}}Y_{j_0,q^{m}}\underset{j\neq j_0}{\prod}m_0^{(j)}$$
where for $j\neq j_0$, $m_0^{(j)}=\underset{l\in\ZZ}{\prod}Y_{j,l}^{u_{j,l}(m_0)}$. In a finite sequence adapted to $m_0$, a term $Y_{j_0,q^{l}}$ or $Y_{j_0,q^{m}}$ must come from a $E_{j_0+1}(m_1)$ or a $E_{j_0-1}(m_1)$. So for example we have $m_1<m_0$ in $D$ of the form $m_1=Y_{j_0,q^{m}}Y_{j_0+1,q^{l-1}}\underset{j\neq j_0,j_0+1}{\prod}m_1^{(j)}$. In all cases we get a monomial $m_1<m_0$ in $D$ of the form:
$$m_1=Y_{j_1,q^{m_1}}Y_{j_1+1,q^{l_1}}\underset{j\neq j_1,j_1+1}{\prod}m_1^{(j)}$$
But the term $Y_{j_1+1,q^{l-1}}$ can not come from a $E_{j_1}(m_2)$ because we would have $u_{j_1}(m_2)\geq 2$. So we have $m_2<m_1$ in $D$ of the form:  
$$m_2=Y_{j_2,q^{m_2}}Y_{j_2+2,q^{l_2}}\underset{j\neq j_2,j_2+1,j_2+2}{\prod}m_2^{(j)}$$
This term must come from a $E_{j_2-1}, E_{j_2+3}$. By induction, we get $m_N<m_0$ in $D$ of the form :
$$m_N=Y_{1,q^{m_N}}Y_{n,q^{l_N}}\underset{j\neq 1,..,n}{\prod}m_N^{(j)}=Y_{1,q^{m-N}}Y_{n,q^{l_N}}$$
It is a dominant monomial of $D\subset D_{Y_{i,0}}$ which is not $Y_{i,0}$. It is impossible (proof of lemma \ref{copieun}).\qed

\noindent An analog result is also geometrically proved by Nakajima for the $ADE$-case in \cite{Nab} (it is also algebraically  for $AD$-cases proved in \cite{Nac}). Those results and the explicit formulas in $n=1,2$-cases (see section \ref{fin}) suggest:

\begin{conj}\label{conun} The coefficients of $F_t(\tilde{Y}_{i,0})=\chi_{q,t}(Y_{i,0})$ are in $\NN[t^{\pm}]$.\end{conj}

\noindent In particular for $m\in B$, the coefficients of $E_t(m)$ would be in $\NN[t^{\pm}]$; moreover $\chi_{q,t}(Y_{i,0})$ and $\chi_q(Y_{i,0})$ would have the same monomials, the $t$-algorithm would stop and $\text{Im}(\chi_{q,t})\subset\Yim_t$.

\noindent At the time he wrote this paper the author does not know a general proof of the conjecture. However a case by case investigation seems possible: the cases $G_2, B_2, C_2$ are checked in section \ref{fin} and the cases $F_4, B_n, C_n$ ($n\leq 10$) have been checked on a computer. So a combinatorial proof for series $B_n, C_n$ ($n\geq 2$) analog to the proof of proposition \ref{cpaconj} would complete the picture.

\subsection{Decomposition in irreducible modules}

The proposition \ref{calcex} suggests:

\begin{conj}\label{condeux} For $m\in B$ we have $\pi_+(L_t(m))=L(m)$.\end{conj}

\noindent In the $ADE$-case the conjecture \ref{condeux} is proved by Nakajima with the help of geometry (\cite{Nab}). In particular this conjecture implies that the coefficients of $\pi_+(L_t(m))$ are non negative. It gives a way to compute explicitly the decomposition of a standard module in irreducible modules, because the conjecture \ref{condeux} implies:
$$E(m)=L(m)+\underset{m'<m}{\sum}P_{m',m}(1)L(m')$$
In particular we would have $P_{m',m}(1)\geq 0$.

\noindent In section \ref{exltcalc} we have studied some examples:

-In proposition \ref{calcex} for $\Glie=sl_2$ and $m\in B$ such that $\forall l\in\ZZ, u_l(m)\leq 1$: we have $\pi_+(L_t(m))=F(m)=L(m)$ and:
$$E(m)=\underset{m'\in B/m'\leq m}{\sum}L(m')$$
-For $\Glie=sl_2$ and $m=\tilde{Y}_0^2\tilde{Y}_2$: we have $\pi_+(L_t(m))=F(Y_0)F(Y_0Y_2)=L(m)$ and:
$$E(Y_0^2Y_2)=L(Y_0^2Y_2)+L(Y_0)$$
Note that $L(Y_0^2Y_2)$ has two dominant monomials $Y_0^2Y_2$ and $Y_0$ because $Y_0^2Y_2$ is irregular (lemma \ref{dominl}).

-For $C=B_2$ and $m=\tilde{Y}_{2,0}\tilde{Y}_{1,5}$. The $\pi_+(L_t(\tilde{Y}_{2,0}\tilde{Y}_{1,5}))$ has non negative coefficients and the conjecture implies $E(Y_{2,0}Y_{1,5})=L(Y_{2,0}Y_{1,5})+L(Y_{1,1})$.

\subsection{Further applications and generalizations}

We hope to address the following questions in the future:

\subsubsection{Iterated deformed screening operators}

Our presentation of deformed screening operators as commutators leads to the definition of iterated deformed screening operators. For example in order 2 we set:
$$\tilde{S}_{j,i,t}(m)=[\underset{l\in\ZZ}{\sum}\tilde{S}_{j,l},S_{i,t}(m)]$$

\subsubsection{Possible generalizations}

\noindent Some generalizations of the approach used in this article will be studied: 

a) the theory of $q$-characters at roots of unity (\cite{Fre3}) suggests a generalization to the case $q^N=1$. 

b) in this article we decided to work with $\Yim_t$ which is a quotient of $\Yim_u$. The same construction with $\Yim_u$ will give characters with an infinity of parameters of deformation $t_r=\text{exp}(\underset{l>0}{\sum}h^{2l}q^{lr}c_l)$ ($r\in\ZZ$).

c) our construction is independent of representation theory and could be established for other generalized Cartan matrices (in particular for twisted affine cases). 

\clearpage

\section{Appendix}\label{fin}

There are 5 types of semi-simple Lie algebra of rank 2: $A_1\times A_1$, $A_2$, $C_2$, $B_2$, $G_2$ (see for example \cite{Kac}). In each case we give the formula for $E(1),E(2)\in\mathfrak{K}$ and we see that the hypothesis of proposition \ref{cpfacile} is verified. In particular we have $E_t(\tilde{Y}_{1,0})=\pi^{-1}(E(1)),E_t(\tilde{Y}_{2,0})=\pi^{-1}(E(2))\in\mathfrak{K}_t$. 

\noindent Following \cite{Fre}, we represent the $E(1),E(2)\in\mathfrak{K}$ as a $I\times\ZZ$-oriented colored tree. For $\chi\in\mathfrak{K}$ the tree $\Gamma_{\chi}$ is defined as follows: the set of vertices is the set of $\Yim$-monomials of $\chi$. We draw an arrow of color $(i,l)$ from $m_1$ to $m_2$ if $m_2=A_{i,l}^{-1}m_1$ and if in the decomposition $\chi=\underset{m\in B_i}{\sum}\mu_m L_i(m)$ there is $M\in B_i$ such that $\mu_M\neq 0$ and $m_1,m_2$ appear in $L_i(M)$.

\noindent Then we give a formula for $E_t(\tilde{Y}_{1,0}),E_t(\tilde{Y}_{2,0})$ and we write it in $\mathfrak{K}_{1,t}$ and in $\mathfrak{K}_{2,t}$. 

\subsection{$A_1\times A_1$-case} The Cartan matrix is $C=\begin{pmatrix}2 & 0\\0 & 2 \end{pmatrix}$ and $r_1=r_2=1$ (note that in this case the computations keep unchanged for all $r_1,r_2$).

$$\xymatrix{Y_{1,0} \ar[d]^{1,1}
\\Y_{1,2}^{-1}}\text{ and }\xymatrix{Y_{2,0} \ar[d]^{2,1}
\\Y_{2,2}^{-1}}$$

$$E_t(\tilde{Y}_{1,0})=\pi^{-1}(Y_{1,0}+Y_{1,2}^{-1})=\tilde{Y}_{1,0}(1+t\tilde{A}_{1,1}^{-1})\in\mathfrak{K}_{1,t}$$
$$=\tilde{Y}_{1,0}+\tilde{Y}_{1,2}^{-1}\in\mathfrak{K}_{2,t}$$
$$E_t(\tilde{Y}_{2,0})=\pi^{-1}(Y_{2,0}+Y_{2,2}^{-1})=\tilde{Y}_{2,0}(1+t\tilde{A}_{2,1}^{-1})\in\mathfrak{K}_{2,t}$$
$$=\tilde{Y}_{2,0}+\tilde{Y}_{2,2}^{-1}\in\mathfrak{K}_{1,t}$$

\subsection{$A_2$-case} The Cartan matrix is $C=\begin{pmatrix}2 & -1\\-1 & 2 \end{pmatrix}$. It is symmetric, $r_1=r_2=1$:

$$\xymatrix{Y_{1,0} \ar[d]^{1,1}
\\Y_{1,2}^{-1}Y_{2,1}\ar[d]^{2,2}
\\Y_{2,3}^{-1}}
\text{ and }
\xymatrix{Y_{2,0} \ar[d]^{2,1}
\\Y_{2,2}^{-1}Y_{1,1}\ar[d]^{1,2}
\\Y_{1,3}^{-1}}
$$

$$E_t(\tilde{Y}_{1,0})=\pi^{-1}(Y_{1,0}+Y_{1,2}^{-1}Y_{2,1}+Y_{2,3}^{-1})=\tilde{Y}_{1,0}(1+t\tilde{A}_{1,1}^{-1})+\tilde{Y}_{2,3}^{-1}\in\mathfrak{K}_{1,t}$$
$$=\tilde{Y}_{1,0}+:\tilde{Y}_{1,2}^{-1}\tilde{Y}_{2,1}:(1+t\tilde{A}_{2,2}^{-1})\in\mathfrak{K}_{2,t}$$
$$E_t(\tilde{Y}_{2,0})=\pi^{-1}(Y_{2,0}+Y_{2,2}^{-1}Y_{1,1}+Y_{1,3}^{-1})=\tilde{Y}_{2,0}(1+t\tilde{A}_{2,1}^{-1})+\tilde{Y}_{1,3}^{-1}\in\mathfrak{K}_{2,t}$$
$$=\tilde{Y}_{2,0}+:\tilde{Y}_{2,2}^{-1}\tilde{Y}_{1,1}:(1+t\tilde{A}_{1,2}^{-1})\in\mathfrak{K}_{1,t}$$

\clearpage

\subsection{$C_2,B_2$-case} The two cases are dual so it suffices to compute for the Cartan matrix $C=\begin{pmatrix}2 & -2\\-1 & 2 \end{pmatrix}$ and $r_1=1$, $r_2=2$.

$$\xymatrix{Y_{1,0} \ar[d]^{1,1}
\\Y_{1,2}^{-1}Y_{2,1}\ar[d]^{2,3}
\\Y_{2,5}^{-1}Y_{1,4}\ar[d]^{1,5}
\\Y_{1,6}^{-1}}
\text{ and }
\xymatrix{Y_{2,0}\ar[d]^{2,2}
\\Y_{2,4}^{-1}Y_{1,1}Y_{1,3}\ar[d]^{1,4}
\\Y_{1,1}Y_{1,5}^{-1}\ar[d]^{1,2}
\\Y_{1,3}^{-1}Y_{1,5}^{-1}Y_{2,2}\ar[d]^{2,4}
\\Y_{2,6}^{-1}}$$

$$E_t(\tilde{Y}_{1,0})=\pi^{-1}(Y_{1,0}+Y_{1,2}^{-1}Y_{2,1}+Y_{2,5}^{-1}Y_{1,4}+Y_{1,6}^{-1})$$
$$=\tilde{Y}_{1,0}(1+t\tilde{A}_{1,1}^{-1})+:\tilde{Y}_{2,5}^{-1}\tilde{Y}_{1,4}:(1+t\tilde{A}_{1,5}^{-1})\in\mathfrak{K}_{1,t}$$
$$=\tilde{Y}_{1,0}+:\tilde{Y}_{1,2}^{-1}\tilde{Y}_{2,1}:(1+t\tilde{A}_{2,3}^{-1})+\tilde{Y}_{1,6}^{-1}\in\mathfrak{K}_{2,t}$$

$$E_t(\tilde{Y}_{2,0})=\pi^{-1}(Y_{2,0}+Y_{2,4}^{-1}Y_{1,1}Y_{1,3}+Y_{1,1}Y_{1,5}^{-1}+Y_{1,3}^{-1}Y_{1,5}^{-1}Y_{2,2}+Y_{2,6}^{-1})$$
$$=\tilde{Y}_{2,0}+:\tilde{Y}_{2,4}^{-1}\tilde{Y}_{1,1}\tilde{Y}_{1,3}:(1+t\tilde{A}_{1,4}^{-1}+t^2\tilde{A}_{1,4}^{-1}\tilde{A}_{1,2}^{-1})+\tilde{Y}_{2,6}^{-1}\in\mathfrak{K}_{1,t}$$
$$=\tilde{Y}_{2,0}(1+t\tilde{A}_{2,2}^{-1})+:\tilde{Y}_{1,1}\tilde{Y}_{1,5}^{-1}:+:\tilde{Y}_{1,3}^{-1}\tilde{Y}_{1,5}^{-1}\tilde{Y}_{2,2}:(1+t\tilde{A}_{2,4}^{-1})\in\mathfrak{K}_{2,t}$$

\clearpage

\subsection{$G_2$-case} The Cartan matrix is $C=\begin{pmatrix}2 & -3\\-1 & 2 \end{pmatrix}$ and $r_1=1$, $r_2=3$.

\subsubsection{First fundamental representation}

$$\xymatrix{Y_{1,0} \ar[d]^{1,1}
\\Y_{1,2}^{-1}Y_{2,1}\ar[d]^{2,4}
\\Y_{2,7}^{-1}Y_{1,4}Y_{1,6}\ar[d]^{1,7}
\\Y_{1,4}Y_{1,8}^{-1}\ar[d]^{1,5}
\\Y_{1,6}^{-1}Y_{1,8}^{-1}Y_{2,5}\ar[d]^{2,8}
\\Y_{2,11}^{-1}Y_{1,10}\ar[d]^{1,11}
\\Y_{1,12}^{-1}}$$

$$E_t(\tilde{Y}_{1,0})=\pi^{-1}(Y_{1,0}+Y_{1,2}^{-1}Y_{2,1}+Y_{2,7}^{-1}Y_{1,4}Y_{1,6}+Y_{1,4}Y_{1,6}+Y_{1,6}^{-1}Y_{1,8}^{-1}Y_{2,5}+Y_{2,11}^{-1}Y_{1,10}+Y_{1,12}^{-1})$$
$$=\tilde{Y}_{1,0}(1+\tilde{A}_{1,1}^{-1})
+:\tilde{Y}_{2,7}^{-1}\tilde{Y}_{1,4}\tilde{Y}_{1,6}:(1+t\tilde{A}_{1,7}^{-1}+t^2\tilde{A}_{1,7}^{-1}\tilde{A}_{1,5}^{-1})
+:\tilde{Y}_{2,11}^{-1}\tilde{Y}_{1,10}:(1+t\tilde{A}_{1,11}^{-1})\in\mathfrak{K}_{1,t}$$
$$=\tilde{Y}_{1,0}
+:\tilde{Y}_{1,2}^{-1}\tilde{Y}_{2,1}:(1+t\tilde{A}_{2,4}^{-1})
+:\tilde{Y}_{1,4}\tilde{Y}_{1,6}:
+:\tilde{Y}_{1,6}^{-1}\tilde{Y}_{1,8}^{-1}\tilde{Y}_{2,5}:(1+t\tilde{A}_{2,8}^{-1})+:\tilde{Y}_{1,12}^{-1}:\in\mathfrak{K}_{2,t}$$

\clearpage

\subsubsection{Second fundamental representation}

$$
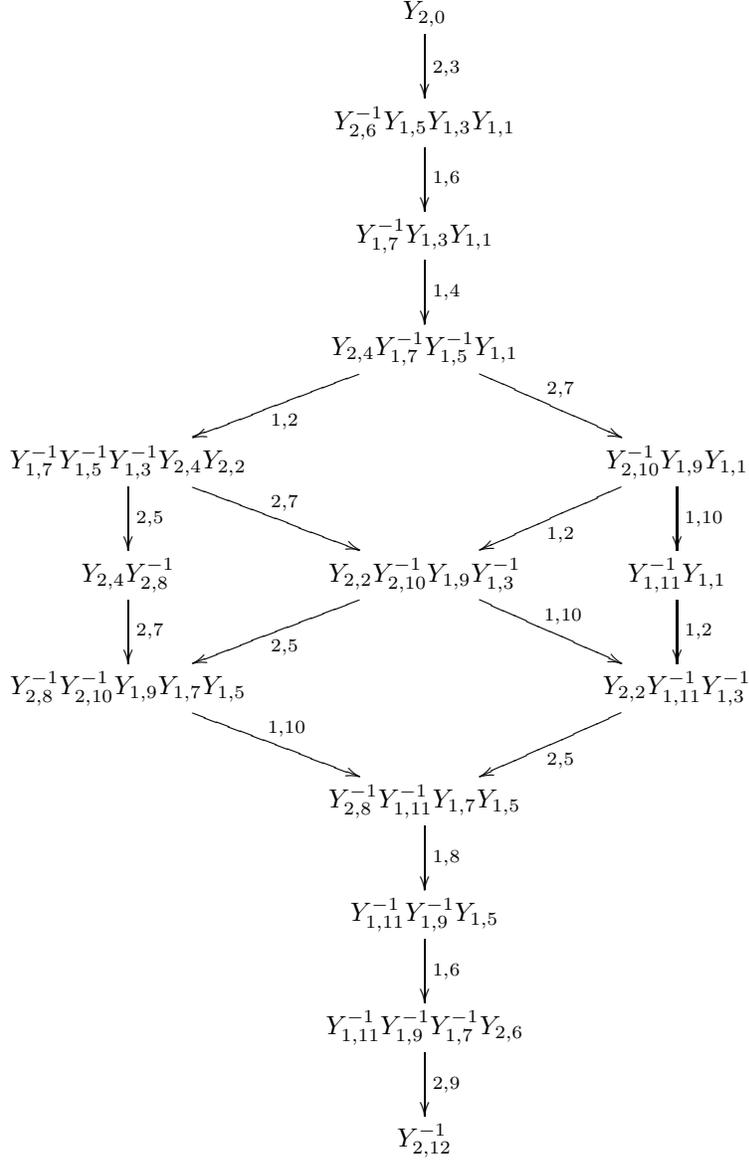
\xymatrix{&Y_{2,0} \ar[d]^{2,3}&
\\&Y_{2,6}^{-1}Y_{1,5}Y_{1,3}Y_{1,1}\ar[d]^{1,6}&
\\&Y_{1,7}^{-1}Y_{1,3}Y_{1,1}\ar[d]^{1,4}&
\\&Y_{2,4}Y_{1,7}^{-1}Y_{1,5}^{-1}Y_{1,1}\ar[ld]^{1,2}\ar[rd]^{2,7}&
\\Y_{1,7}^{-1}Y_{1,5}^{-1}Y_{1,3}^{-1}Y_{2,4}Y_{2,2}\ar[d]^{2,5}\ar[rd]^{2,7}&&Y_{2,10}^{-1}Y_{1,9}Y_{1,1}\ar[ld]^{1,2}\ar[d]^{1,10}
\\Y_{2,4}Y_{2,8}^{-1}\ar[d]^{2,7}&Y_{2,2}Y_{2,10}^{-1}Y_{1,9}Y_{1,3}^{-1}\ar[ld]^{2,5}\ar[rd]^{1,10}&Y_{1,11}^{-1}Y_{1,1}\ar[d]^{1,2}
\\Y_{2,8}^{-1}Y_{2,10}^{-1}Y_{1,9}Y_{1,7}Y_{1,5}\ar[rd]^{1,10}&&Y_{2,2}Y_{1,11}^{-1}Y_{1,3}^{-1}\ar[ld]^{2,5}
\\&Y_{2,8}^{-1}Y_{1,11}^{-1}Y_{1,7}Y_{1,5}\ar[d]^{1,8}&
\\&Y_{1,11}^{-1}Y_{1,9}^{-1}Y_{1,5}\ar[d]^{1,6}&
\\&Y_{1,11}^{-1}Y_{1,9}^{-1}Y_{1,7}^{-1}Y_{2,6}\ar[d]^{2,9}&
\\&Y_{2,12}^{-1}&}$$

$$E_t(\tilde{Y}_{2,0})=\pi^{-1}(Y_{2,0}
+Y_{2,6}^{-1}Y_{1,5}Y_{1,3}Y_{1,1}
+Y_{1,7}^{-1}Y_{1,3}Y_{1,1}0
+Y_{2,4}Y_{1,7}^{-1}Y_{1,5}^{-1}Y_{1,1}
+Y_{1,7}^{-1}Y_{1,5}^{-1}Y_{1,3}^{-1}Y_{2,4}Y_{2,2}$$
$$+Y_{2,10}^{-1}Y_{1,9}Y_{1,1}
+Y_{2,4}Y_{2,8}^{-1}+Y_{2,2}Y_{2,10}^{-1}Y_{1,9}Y_{1,3}^{-1}
+Y_{1,11}^{-1}Y_{1,1}
+Y_{2,8}^{-1}Y_{2,10}^{-1}Y_{1,9}Y_{1,7}Y_{1,5}
+Y_{2,2}Y_{1,11}^{-1}Y_{1,3}^{-1}$$
$$+Y_{2,8}^{-1}Y_{1,11}^{-1}Y_{1,7}Y_{1,5}
+Y_{1,11}^{-1}Y_{1,9}^{-1}Y_{1,5}+Y_{1,11}^{-1}Y_{1,9}^{-1}Y_{1,7}^{-1}Y_{2,6}+Y_{2,12}^{-1})$$
We use the following relations to write $E_t(\tilde{Y}_{2,0})$ in $\mathfrak{K}_{1,t}$ and in $\mathfrak{K}_{2,t}$: $\tilde{A}_{1,2}\tilde{A}_{2,7}=\tilde{A}_{2,7}\tilde{A}_{1,2}$, $\tilde{A}_{2,5}\tilde{A}_{2,7}=\tilde{A}_{2,7}\tilde{A}_{2,5}$, $\tilde{A}_{1,2}\tilde{A}_{1,10}=\tilde{A}_{1,10}\tilde{A}_{1,2}$, $\tilde{A}_{2,5}\tilde{A}_{1,10}=\tilde{A}_{1,10}\tilde{A}_{2,5}$.
$$E_t(\tilde{Y}_{2,0})=\tilde{Y}_{2,0}
+:\tilde{Y}_{2,6}^{-1}\tilde{Y}_{1,5}\tilde{Y}_{1,3}\tilde{Y}_{1,1}:(1+t\tilde{A}_{1,6}^{-1}(1+t\tilde{A}_{1,4}^{-1}(1+t\tilde{A}_{1,2}^{-1})))
+:\tilde{Y}_{2,10}^{-1}\tilde{Y}_{1,9}\tilde{Y}_{1,1}:(1+t\tilde{A}_{1,2}^{-1})(1+t\tilde{A}_{1,10}^{-1})$$
$$+:\tilde{Y}_{2,4}\tilde{Y}_{2,8}^{-1}:
+:\tilde{Y}_{2,8}^{-1}\tilde{Y}_{2,10}^{-1}\tilde{Y}_{1,9}\tilde{Y}_{1,7}\tilde{Y}_{1,5}:(1+t\tilde{A}_{1,10}^{-1}(1+t\tilde{A}_{1,8}^{-1}(1+t\tilde{A}_{1,6}^{-1})))
+\tilde{Y}_{2,12}^{-1}\in\mathfrak{K}_{1,t}$$
$$=\tilde{Y}_{2,0}(1+t\tilde{A}_{2,3}^{-1})
+:\tilde{Y}_{1,7}^{-1}\tilde{Y}_{1,3}\tilde{Y}_{1,1}:
+:\tilde{Y}_{2,4}\tilde{Y}_{1,7}^{-1}\tilde{Y}_{1,5}^{-1}\tilde{Y}_{1,1}:(1+t\tilde{A}_{2,7}^{-1})$$
$$+:\tilde{Y}_{1,7}^{-1}\tilde{Y}_{1,5}^{-1}\tilde{Y}_{1,3}^{-1}\tilde{Y}_{2,4}\tilde{Y}_{2,2}:(1+t\tilde{A}_{2,7}^{-1})(1+t\tilde{A}_{2,5}^{-1})
+:\tilde{Y}_{1,11}^{-1}\tilde{Y}_{1,1}:
+:\tilde{Y}_{2,2}\tilde{Y}_{1,11}^{-1}\tilde{Y}_{1,3}^{-1}:(1+t\tilde{A}_{2,5}^{-1})$$
$$+:\tilde{Y}_{1,11}^{-1}\tilde{Y}_{1,9}^{-1}\tilde{Y}_{1,5}:
+:\tilde{Y}_{1,11}^{-1}\tilde{Y}_{1,9}^{-1}\tilde{Y}_{1,7}^{-1}\tilde{Y}_{2,6}:(1+t\tilde{A}_{2,9}^{-1})\in\mathfrak{K}_{2,t}$$

\twocolumn

\section*{Notations}

\begin{tabular}{lll}

$A$ & set of $\Yim$-monomials & p \pageref{a}

\\ $A_t$ & set of $\Yim_t$-monomials &  p \pageref{tat}

\\ $A_m^{\text{inv}}, B_m^{\text{inv}}$ & set of $\Yim_t$-monomials & p \pageref{ainv}

\\$\overset{\infty}{A}_t$ & product module & p \pageref{infat}

\\$\alpha$ & map $(I\times \ZZ)^2\rightarrow \ZZ$ & p \pageref{alpha}

\\$\alpha(m)$ & character & p \pageref{alpham}

\\$a_i[m]$ & element of $\mathcal{H}$ & p \pageref{aim}

\\$\tilde{A}_{i,l}, \tilde{A}_{i,l}^{-1}$ & elements of $\Yim_u$ or $\Yim_t$ & p \pageref{tail}

\\$A_{i,l}, A_{i,l}^{-1}$ & elements of $\Yim$ & p \pageref{ail}

\\$B$ & a set of $\Yim$-monomials & p \pageref{b}

\\$B_i$, $B_J$ & a set of $\Yim$-monomials & p \pageref{bi}

\\$(B_{i,j})$ & symmetrized &

\\ &Cartan matrix & p \pageref{symcar}

\\$\beta$ & map $(I\times \ZZ)^2\rightarrow \ZZ$ & p \pageref{beta}

\\$(C_{i,j})$ & Cartan matrix & p \pageref{carmat}

\\$(\tilde{C}_{i,j})$ & inverse of $C$ & p \pageref{invcar}

\\$c_r$ & central element of $\mathcal{H}$ & p \pageref{cr}

\\$d$ & bicharacter & p \pageref{d}

\\$D_{m,K}, D_m$ & set of monomials & p \pageref{dmk}

\\$\tilde{D}_m$ & submodule of $\Yim_t^{\infty}$ & p \pageref{tdm}

\\$E_i(m)$ & element of $\mathfrak{K}_i$ & p \pageref{eim}

\\$E_{i,t}(m)$ & element of $\mathfrak{K}_{i,t}$ & p \pageref{teitm}

\\$E_{i,t}^m$ &map & p \pageref{teitmap}

\\$E(m)$ & element of $\mathfrak{K}$ & p \pageref{em}

\\$E_t(m)$ & element of $\mathfrak{K}_{t}^{\infty}$ & p \pageref{tetm}

\\$F_i(m)$ & element of $\mathfrak{K}_i$ & p \pageref{fim}

\\$F_{i,t}(m)$ & element of $\mathfrak{K}_{i,t}$ & p \pageref{tfitm}

\\$F(m)$ & element of $\mathfrak{K}$ & p \pageref{fm}

\\$F_t(m)$ & element of $\mathfrak{K}_{t}^{\infty}$ & p \pageref{tftm}

\\$\gamma$ & map $(I\times \ZZ)^2\rightarrow \ZZ$ & p \pageref{gamma}

\\$\mathcal{H}$ & Heisenberg algebra & p \pageref{zq}

\\$\mathcal{H}^+, \mathcal{H}^-$ & subalgebras of $\mathcal{H}$ & p \pageref{zqplus}

\\$\mathcal{H}_h$ & formal series in $\mathcal{H}$ & p \pageref{zqh}

\\$\mathcal{H}_t$ & quotient of $\mathcal{H}_h$ & p \pageref{zqt}

\\$\mathcal{H}_t^+, \mathcal{H}_t^-$ & subalgebras of $\mathcal{H}_t$ & p \pageref{zqtplus}

\\$\mathfrak{K}_i, \mathfrak{K}_J, \mathfrak{K}$ & subrings of $\Yim$ & p \pageref{ki}

\\$\mathfrak{K}_{i,t}, \mathfrak{K}_{J,t}, \mathfrak{K}_t$ & subrings of $\Yim_t$ & p \pageref{kit}

\\$\mathfrak{K}_{i,t}^{\infty}, \mathfrak{K}_{J,t}^{\infty}, \mathfrak{K}_t^{\infty}$ & subrings of $\Yim_t^{\infty}$ & p \pageref{kitinf}

\\$\chi_q$ & morphism &

\\ & of $q$-characters &p \pageref{chiqdefi}

\\$\chi_{q,t}$ & morphism &

\\ & of $q,t$-characters &p \pageref{chiqt}

\\$L_i(m)$ & element of $\mathfrak{K}_i$ & p \pageref{lim}

\\$L_t(m)$ & element of $\mathfrak{K}_t^{\infty}$ & p \pageref{tltm}

\\$:m:$ & monomial in $A$ & p \pageref{ptpt}

\\$\tilde{m}$ & monomial in $A_t$ & p \pageref{tm}

\\$N,N_t,\mathcal{N},\mathcal{N}_t$ & characters, bicharacters & p \pageref{n}

\\$P(n)$ & property of $n\in\NN$ & p \pageref{pn}

\end{tabular}

\begin{tabular}{lll}

 &&

\\ &&

\\$\pi$ & map & p \pageref{pi}

\\$\pi_r$ & map to $\ZZ$ & p \pageref{pir}

\\$\pi_+,\pi_-$ & endomorphisms of &

\\ &$\mathcal{H}_h$, $\mathcal{H}_t$ & p \pageref{piplus}

\\$q$ & complex number & p \pageref{qz}

\\$Q(n)$ & property of $n\in\NN$ & p \pageref{qn}

\\$\text{Rep}$ & Grothendieck ring & p \pageref{rep}

\\$\text{Rep}_t$ & deformed &

\\ &Grothendieck ring & p \pageref{rept}

\\$s(m_r)_J,s(m_r)$ & sequences of $\ZZ[t^{\pm}]$ & p \pageref{smrt}

\\$S_i$ & screening operator & p \pageref{si}

\\$\tilde{S}_{i,l}$ & screening current & p \pageref{tsil}

\\$S_{i,t}$ & $t$-screening operator & p \pageref{tsit}

\\$t$ & central element of $\Yim_t$ & p \pageref{t}

\\$t_R$ & central element of $\Yim_u$ & p \pageref{tr}

\\$u_{i,l}$ & multiplicity of $Y_{i,l}$ & p \pageref{uil}

\\$u_i$ & sum of the $u_{i,l}$ & p \pageref{ui}

\\$\mathfrak{U}$ & subring of $\QQ(q)$ & p \pageref{mathu}

\\$\U_q(\hat{\Glie})$ & quantum &

\\ & affine algebra &p \pageref{qadefi}

\\$\U_q(\hat{\Hlie})$ & Cartan algebra & p \pageref{qhdefi}

\\$X_{i,l}$ & element of $\text{Rep}$ & p \pageref{xil}

\\$y_i[m]$ & element of $\mathcal{H}$ & p \pageref{yim}

\\$Y_{i,l}, Y_{i,l}^{-1}$ & elements of $\Yim$ & p \pageref{ail}

\\$\tilde{Y}_{i,l},\tilde{Y}_{i,l}^{-1}$ & elements of $\Yim_u$ or $\Yim_t$ & p \pageref{tail}

\\$\Yim$ & subalgebra of $\mathcal{H}_h$ & p \pageref{y}

\\$\Yim_t$ & quotient of $\Yim_u$ & p \pageref{tyt}

\\$\Yim_t^+, \Yim_t^-$ & subalgebras of $\mathcal{H}_t$ & p \pageref{tytplus}

\\$\Yim_u$ & subalgebra of $\mathcal{H}_h$ & p \pageref{yu}

\\$\Yim_{i,t}$ & $\Yim_t$-module & p \pageref{yit}

\\$\Yim_{i,u}$ & $\Yim_u$-module & p \pageref{yiu}

\\$\Yim_t^{\infty}, \Yim_t^{A, \infty}$ & submodules of $\overset{\infty}{A}_t$ & p \pageref{tytinf}

\\$\Yim_t^A, \Yim_t^{A,K}$ & submodules of $\Yim_t$ & p \pageref{tyta}

\\$z$ & indeterminate & p \pageref{qz}

\\$::$ & endomorphism of &

\\ &$\mathcal{H}, \mathcal{H}_h, \Yim_u, \Yim_t$ & p \pageref{ptpt}

\\$*$ & deformed &

\\ &multiplication & p \pageref{star}

\end{tabular}

\onecolumn

\end{document}